\documentclass[a4paper]{amsart}
\let\oldtocsection=\tocsection
\let\oldtocsubsection=\tocsubsection
\let\oldtocsubsubsection=\tocsubsubsection
\renewcommand{\tocsection}[2]{\hspace{0em}\oldtocsection{#1}{#2}}
\renewcommand{\tocsubsection}[2]{\hspace{2em}\oldtocsubsection{#1}{#2}}
\renewcommand{\tocsubsubsection}[2]{\hspace{4em}\oldtocsubsubsection{#1}{#2}}
\usepackage{amssymb}
\usepackage{verbatim}
\usepackage[dvips]{graphicx}
\theoremstyle{plain}
  \newtheorem{thm}{Theorem}[section]
  \newtheorem*{thm*}{Theorem}

  \newtheorem{conj}[thm]{Conjecture}

  \newtheorem*{obs*}{Observation}
\theoremstyle{definition}
  \newtheorem{defn}[thm]{Definition}

\theoremstyle{remark}
  \newtheorem{rem}[thm]{Remark}
  
  \newtheorem*{ack}{Acknowledgments}
\newcommand{\Z}{\mathbb{Z}}
\newcommand{\C}{\mathbb{C}}
\newcommand{\R}{\mathbb{R}}

\newcommand{\Vol}{\operatorname{Vol}}
\newcommand{\CS}{\operatorname{CS}}
\newcommand{\cs}{\operatorname{cs}}
\newcommand{\Li}{\operatorname{Li}}
\newcommand{\Int}{\operatorname{Int}}
\newcommand{\arccosh}{\operatorname{arccosh}}
\newcommand{\Res}{\operatorname{Res}}
\newcommand{\Hom}{\operatorname{Hom}}

\newcommand{\tr}{\operatorname{tr}}
\newcommand{\SL}{\mathrm{SL}}
\renewcommand{\sl}{\mathfrak{sl}}
\newcommand{\Ad}[1]{\operatorname{Ad}_{#1}}
\newcommand{\Ker}{\operatorname{Ker}}

\newcommand{\FigEight}{\mathcal{E}}

\renewcommand{\Re}{\operatorname{Re}}
\renewcommand{\Im}{\operatorname{Im}}
\numberwithin{equation}{section}
\allowdisplaybreaks
\begin{document}
\title[Twice-iterated torus knot]
{The colored Jones polynomial, the Chern--Simons invariant, and the Reidemeister torsion of a twice-iterated torus knot}
\author{Hitoshi Murakami}
\address{
Graduate School of Information Sciences,
Tohoku University,
Aramaki-aza-Aoba 6-3-09, Aoba-ku,
Sendai 980-8579, Japan
}
\email{starshea@tky3.3web.ne.jp}
\date{\today}
\begin{abstract}
A generalization of the volume conjecture relates the asymptotic behavior of the colored Jones polynomial of a knot to the Chern--Simons invariant and the Reidemeister torsion of the knot complement associated with a representation of the fundamental group to the special linear group of degree two over complex numbers.
If the knot is hyperbolic, the representation can be regarded as a deformation of the holonomy representation that determines the complete hyperbolic structure.
In this article we study a similar phenomenon when the knot is a twice-iterated torus knot.
In this case, the asymptotic expansion of the colored Jones polynomial splits into sums and each summand is related to the Chern--Simons invariant and the Reidemeister torsion associated with a representation.
\end{abstract}
\thanks{This work was supported by JSPS KAKENHI Grant Numbers 23340115, 24654041.}
\maketitle
\setcounter{tocdepth}{3}
Let $J_N(K;q)\in\Z[q,q^{-1}]$ be the colored Jones polynomial of a knot $K$ in the three-sphere $S^3$ associated with the irreducible $N$-dimensional representation of the Lie algebra $\sl_2(\C)$ \cite{Jones:BULAM385,Kirillov/Reshetikhin:1989}.
We normalize it so that $J_N(\text{unknot};q)=1$.
Note that $J_2(K;q)$ is the original Jones polynomial.
R.~Kashaev conjectured \cite{Kashaev:LETMP97} that his knot invariant $\langle{K}\rangle_N\in\C$ introduced in \cite{Kashaev:MODPLA95} would grow exponentially with growth rate the volume of the knot complement $S^3\setminus{K}$ when the integer parameter $N$ goes to the infinity if the knot $K$ is hyperbolic.
J.~Murakami and the author \cite{Murakami/Murakami:ACTAM12001} proved that Kashaev's invariant coincides with $J_N\bigl(K;\exp(2\pi\sqrt{-1}/N)\bigr)$ and generalized Kashaev's conjecture to general knots.
\begin{conj}[Volume Conjecture,\cite{Kashaev:LETMP97,Murakami/Murakami:ACTAM12001}]\label{conj:VC}
For any knot, we have
\begin{equation}\label{eq:VC}
  2\pi
  \lim_{N\to\infty}
  \frac{\log\left|J_N\bigl(K;\exp(2\pi\sqrt{-1}/N)\bigr)\right|}{N}
  =
  \Vol(S^3\setminus{K}).
\end{equation}
Here $\Vol$ is the simplicial volume {\rm(}or the Gromov norm{\rm)} \cite{Gromov:INSHE82}, which is the sum of the hyperbolic volumes of the hyperbolic pieces in the JSJ decomposition \cite{Jaco/Shalen:MEMAM79,Johannson:1979} of the knot complement.
\end{conj}
Note that when the knot $K$ is hyperbolic, that is, its complement possesses a (unique) complete hyperbolic structure with finite volume, $\Vol(S^3\setminus{K})$ is just the hyperbolic volume associated with the complete hyperbolic structure.
So the volume conjecture states that the colored Jones polynomial would know the volume corresponding to the complete hyperbolic structure.
It is well known that we can always perturb the complete structure to get incomplete structures \cite{Thurston:GT3M}.
When we perturb $2\pi\sqrt{-1}$ in the left hand side of \eqref{eq:VC}, it is expected that we can also have the volume of the perturbed hyperbolic manifold:
\begin{conj}[\cite{Gukov:COMMP2005,Murakami:ADVAM22007,Gukov/Murakami:FIC2008}]\label{conj:Gukov/Murakami}
Let $K$ be a hyperbolic knot.
For a small complex number $u$, the following limit exists:
\begin{equation*}
  \lim_{N\to\infty}
  \frac{\log{J_N\Bigl(K;\exp\bigl((2\pi\sqrt{-1}+u)/N\bigr)\Bigr)}}{N}.
\end{equation*}
Put
\begin{equation*}
  S(u)
  :=
  (2\pi\sqrt{-1}+u)
  \lim_{N\to\infty}
  \frac{\log{J_N\Bigl(K;\exp\bigl((2\pi\sqrt{-1}+u)/N\bigr)\Bigr)}}{N}
\end{equation*}
and define
\begin{equation*}
  v(u)
  :=
  2\frac{d\,S(u)}{d\,u}-2\pi\sqrt{-1}.
\end{equation*}
Then the volume $\Vol(S^3\setminus{K};u)$ of the knot complement with the hyperbolic structure associated with the representation parametrized by $u$ is given as follows:
\begin{equation*}
  \Vol(S^3\setminus{K};u)
  =
  \Im{S(u)}
  -\pi\Re(u)
  -\frac{1}{2}\Re(u)\Im\bigl(v(u)\bigr).
\end{equation*}
Here the representation sends the meridian of the knot to $\begin{pmatrix}e^{u/2}&\ast\\0&e^{-u/2}\end{pmatrix}$ and the longitude to $\begin{pmatrix}e^{v(u)/2}&\ast\\0&e^{-v(u)/2}\end{pmatrix}$.
\end{conj}
The conjecture is true for the figure-eight knot \cite{Murakami/Yokota:JREIA2007}.
\par
This means that for large $N$ we can write
\begin{multline*}
  J_N\Bigl(K;\exp\bigl((2\pi\sqrt{-1}+u)/N\bigr)\Bigr)
  \\
  =
  \exp\left[\frac{S(u)N}{2\pi\sqrt{-1}+u}\right]
  \times(\text{a function of $N$ with polynomial growth})
\end{multline*}
for a complex function $S(u)$.
It is also conjectured that $S(u)$ also determines the $\SL(2;\C)$ Chern--Simons invariant (see Section~\ref{sec:CS} for details).
\par
The polynomial growth term in the above equation is expected to determine the twisted Reidemeister torsion.
\begin{conj}[\cite{Gukov/Murakami:FIC2008,Dimofte/Gukov:Columbia}]\label{conj:Dimofte/Gukov}
For a hyperbolic knot $K$, we have
\begin{multline*}
  J_N\Bigl(K;\exp\bigl((2\pi\sqrt{-1}+u)/N\bigr)\Bigr)
  \\
  =
  \frac{\sqrt{-\pi}}{2\sinh(u/2)}
  \mathbb{T}^{K}_{\mu}(u)^{-1/2}
  \left(\frac{N}{2\pi\sqrt{-1}+u}\right)^{1/2}
  \exp\left[\frac{S(u)N}{2\pi\sqrt{-1}+u}\right]
\end{multline*}
for small $u\ne0$, where $\mathbb{T}^{K}_\mu(u)$ is the twisted Reidemeister torsion of the representation parametrized by $u$ associated with the meridian $\mu$.
\end{conj}
This conjecture is proved for the figure-eight knot with real $u$ \cite{Murakami:JTOP2013}.
\par
When $K$ is not hyperbolic, especially when $K$ has no hyperbolic pieces, that is when $K$ is an iterated torus knot, then the volume of its complement is zero.
Kashaev and O.~Tirkkonen proved the volume conjecture for torus knots \cite{Kashaev/Tirkkonen:ZAPNS2000} and R.~Van der Veen proved the conjecture for general iterated torus knots \cite{van_der_Veen:2008}.
In \cite{Hikami/Murakami:COMCM2008,Hikami/Murakami:Bonn} K.~Hikami and the author showed that for a torus knot the asymptotic expansion of the colored Jones polynomial splits into sums each of which corresponds to a representation of the fundamental group into $\SL(2;\C)$.
Moreover we gave a topological interpretation for each summand.
\begin{thm}[\cite{Hikami/Murakami:Bonn}]\label{thm:Hikami/Murakami_intro}
Let $T(m,n)$ be the torus knot of type $(m,n)$ for coprime positive integers $m$ and $n$.
If $\xi$ is a complex number which is not purely imaginary and with $\Im{\xi}\ge0$, then we have the following asymptotic equivalence:
\begin{multline*}
  J_{N}\bigl(T(m,n);\exp(\xi/N)\bigr)
  \\
  \underset{N\to\infty}{\sim}
  \frac{1}{\Delta\bigl(T(m,n);\exp(\xi)\bigr)}
  +
  \frac{\sqrt{-\pi}}{2\sinh(\xi/2)}\sqrt{\frac{N}{\xi}}
  \sum_{k}\tau_{k}\exp\left[S_k(\xi)\frac{N}{\xi}\right],
\end{multline*}
where
$f(N)\underset{N\to\infty}{\sim}g(N)$ means that $f(N)=\bigl(1+o(1)\bigr)g(N)$ for $N\to\infty$,
\begin{align*}
  S_k(\xi)
  &:=
  \frac{-(2k\pi\sqrt{-1}-mn\xi)^2}{4mn},
  \\
  \intertext{and}
  \tau_{k}
  &:=
  (-1)^{k+1}
  \frac{4\sin(k\pi/m)\sin(k\pi/n)}{\sqrt{mn}}.
\end{align*}
Moreover $\tau_k^{-2}$ is the homological twisted Reidemeister torsion $\mathbb{T}^{T(a,b)}_{\mu}(\rho_k)$ of an irreducible representation $\rho_k\colon\pi_1(S^3\setminus{T(m,n)})\to\SL(2;\C)$ associated with the meridian $\mu$, and $S_k(\xi)-\pi\sqrt{-1}u-\frac{1}{4}uv_k$ is the $\SL(2;\C)$ Chern--Simons invariant of $\rho_k$ with respect to the pair $(u,v_k)$ with $u:=\xi-2\pi\sqrt{-1}$ and $v_k:=2\dfrac{d\,S_k(\xi)}{d\,\xi}\Big|_{\xi:=2\pi\sqrt{-1}+u}-2\pi\sqrt{-1}$.
\par
Note that $\rho_k(\mu)=\begin{pmatrix}e^{u/2}&\ast\\0&e^{-u/2}\end{pmatrix}$ and that $\Delta\bigl(T(m,n);\exp(\xi)\bigr)$ can be regarded as the Reidemeister torsion of the Abelian representation $\pi_1(S^3\setminus{T(m,n)})\to\SL(2;\C)$ sending $\mu$ to $\begin{pmatrix}e^{u/2}&0\\0&e^{-u/2}\end{pmatrix}$ with $u:=\xi-2\pi\sqrt{-1}$ \cite{Milnor:ANNMA21962,Milnor:cyclic,Turaev:USPMN1986}.
\end{thm}
See \cite{Hikami/Murakami:Bonn} for more details and Theorem~\ref{thm:Hikami/Murakami} for the case of $T(2,2a+1)$.
\par
The aim of this article is to investigate the asymptotic behavior of the colored Jones polynomial for a twice-iterated torus knot.
The following is our main theorem.
\begin{thm}\label{thm:main}
Let $T(2,2a+1)^{(2,2b+1)}$ be the twice-iterated torus knot for integers $a$ and $b$ with $a>0$, $b>0$ and $2b+1-4(2a+1)>0$.
If $\xi$ is a complex number which is not purely imaginary and with $\Im{\xi}\ge0$, then the following asymptotic equivalence holds.
\begin{equation*}
\begin{split}
  &J_{N}\bigl(T(2,2a+1)^{(2,2b+1)};\xi\bigr)
  \\
  =&
  \frac{1}{\Delta(T(2,2a+1)^{2b+1});\exp\xi}
  +
  \frac{\sqrt{-\pi}}{2\sinh(\xi/2)}
  \sqrt{\frac{N}{\xi}}
  \sum_{j}
  \tau_1(\xi;j)
  \exp\left[\frac{N}{\xi}S_1(\xi;j)\right]
  \\
  &
  +
  \frac{\sqrt{-\pi}}{2\sinh(\xi/2)}
  \sqrt{\frac{N}{\xi}}
  \sum_{k}
  \tau_2(\xi;k)
  \exp\left[\frac{N}{\xi}S_2(\xi;k)\right]
  \\
  &
  +
  \frac{\pi}{2\sinh(\xi/2)}
  \frac{N}{\xi}
  \sum_{l,m}
  \tau_3(\xi;l,m)
  \exp\left[\frac{N}{\xi}S_3(\xi;l,m)\right],
\end{split}
\end{equation*}
where
\begin{align*}
  \tau_1(\xi;j)
  &:=
  (-1)^j
  \sqrt{\frac{2}{\beta}}
  \frac{\sin\left(\frac{2(2j+1)\pi}{\beta}\right)}
       {\cos\left(\frac{(2j+1)\alpha\pi}{\beta}\right)},
  \\
  S_1(\xi;j)
  &:=
  (2j+1)\xi\pi\sqrt{-1}
  +
  \frac{-\beta\xi^2}{2}
  +
  \frac{(2j+1)^2\pi^2}{2\beta},
  \\
  \tau_2(\xi;k)
  &:=
  (-1)^{k+1}
  \sqrt{\frac{2}{\alpha}}
  \frac{\sin\left(\frac{(2k+1)\pi}{\alpha}\right)}{\cosh\left(\frac{(\beta-4\alpha)\xi}{2}\right)},
  \\
  S_2(\xi;k)
  &:=
  2(2k+1)\xi\pi\sqrt{-1}
  +
  -2\alpha\xi^2
  +
  \frac{(2k+1)^2\pi^2}{2\alpha},
  \\
  \tau_3(\xi;l,m)
  &:=
  (-1)^{l+m}
  \frac{4}{\sqrt{\alpha(\beta-4\alpha)}}
  \sin\left(\frac{(2m+1)\pi}{\alpha}\right),
  \\
  S_3(\xi;l,m)
  &:=
  (2l+1)\xi\pi\sqrt{-1}
  +
  \frac{-\beta\xi^2}{2}
  \\
  &\quad+
  \frac{\pi^2}{2\alpha(\beta-4\alpha)}\bigl((2l+1)^2\alpha+(2m+1)^2\beta-4(2l+1)(2m+1)\alpha\bigr).
\end{align*}
Moreover $S_1(\xi;j)$ \rm{(}$S_2(\xi;k)$ and $S_3(\xi;l,m)$, respectively\rm{)} determines the $\SL(2;\C)$ Chern--Simons invariant of a certain irreducible representation $\rho^{\rm{AN}}_{u,j}$ \rm{(}$\rho^{\rm{NA}}_{u,k}$ and $\rho^{\rm{NN}}_{u,l,m}$, respectively\rm{)} from $\pi_1\left(S^3\setminus{T(2,2a+1)^{(2,2b+1)}}\right)$ to $\SL(2;\C)$, and $\tau_1(\xi;j)$ \rm{(}$\tau_2(\xi;k)$ and $\tau_3(\xi;l,m)$, respectively\rm{)} is related to the twisted Reidemeister torsion of $\rho^{\rm{AN}}_{u,j}$ \rm{(}$\rho^{\rm{NA}}_{u,k}$ and $\rho^{\rm{NN}}_{u,l,m}$, respectively\rm{)}.
See Subsection~\ref{subsec:topological_iterated_torus_knot} for details.
\end{thm}
This article is organized as follows.
In Section~\ref{sec:rep}, we describe non-Abelian representations of $\pi_1\left(S^3\setminus{K}\right)$ to $\SL(2;\C)$, where $K$ is the figure-eight knot $\mathcal{E}$, the torus knot of type $(2,2a+1)$, denoted $T(2,2a+1)$, and the $(2,2b+1)$-cable of $T(2,2a+1)$, denoted $T(2,2a+1)^{(2,2b+1)}$.
In Section~\ref{sec:CS}, we calculate the Chern--Simons invariant of the representations that I describe in Section~\ref{sec:rep}.
Section~\ref{sec:Reidemeister} is devoted to the calculation of the twisted Reidemeister torsion of $\mathcal{E}$, $T(2,2a+1)$, and $T(2,2a+1)^{(2,2b+1)}$.
I calculate the colored Jones polynomials of these knots and study their asymptotic behaviors in Section~\ref{sec:Jones}.
In the last section, Section~\ref{sec:topological}, I interpret the coefficients appearing in the asymptotic expansion of the colored Jones polynomial in terms of the Chern--Simons invariant and the Reidemeister torsion.
In Section
\par
When I prepared this manuscript, I often used Mathematica \cite{Mathematica}.
Some formulas are just copies from the output from Mathematica and may be incorrect because of my fault of copying.
\begin{ack}
This article is prepared for the proceedings of the conference ``The Quantum Topology and Hyperbolic Geometry'' in Nha Trang, Vietnam, 13--17 May, 2013.
I would like to thank the organizers for their hospitality.
\par
Part of this work was done when the author was visiting the Max-Planck Institute for Mathematics, Universit{\'e} Paris Diderot, and the University of Amsterdam.
The author thanks Christian Blanchet, Roland van der Veen, Jinseok Cho, and Satoshi Nawata for helpful discussions.
\end{ack}
\section{Representations of the fundamental group into $\SL(2;\C)$}
\label{sec:rep}
In this section we study representations of the fundamental group of a knot complement into $\SL(2;\C)$.
\begin{rem}
In this article I do not want to show all representations.
It may happen that I exhaust all the representations (up to conjugation), but I do not mind.
\end{rem}
\par
We start with an Abelian representation.
Let $K$ be a knot and let $\left\langle x_1,x_2,\dots,x_n\mid r_1,r_2,\dots,r_{n-1}\right\rangle$ be a Wirtinger presentation of the fundamental group $\pi_1(S^3\setminus{K})$ of a knot complement (see, for example, \cite[Chapter~11]{Lickorish:1997}).
Then the map $\rho^{\rm{A}}\colon\pi_1(S^3\setminus{K})\to\SL(2;\C)$ sending $x_i$ to $\begin{pmatrix}g&0\\0&g^{-1}\end{pmatrix}$ for any $i$ becomes a representation, which is called an Abelian representation because the image of $\rho^{\rm{A}}$ forms an Abelian subgroup of $\SL(2;\C)$.
\par
In the following subsections I will focus on non-Abelian representations.
\subsection{Figure-eight knot}\label{subsec:rep_fig8}
\begin{figure}[h!]
  \includegraphics[scale=0.3]{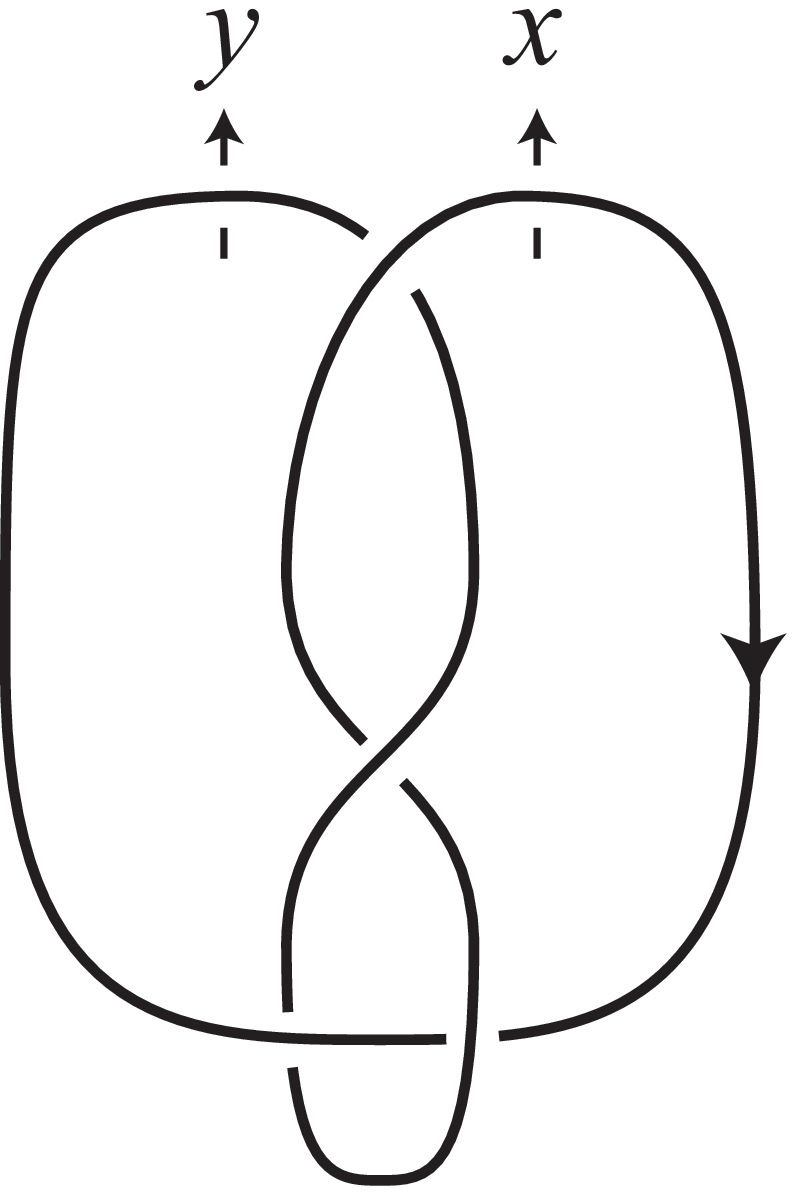}
  \caption{Generators of $\pi_1\left(S^3\setminus\FigEight\right)$}
  \label{fig:fig8}
\end{figure}

Let  $\mathcal{E}$ denote the figure-eight knot.
The fundamental group $\pi_1\left(S^3\setminus\FigEight\right)$ of the figure-eight knot complement has the following presentation (see for example \cite{Murakami:ACTMV2008}):
\begin{equation*}
  \pi_1(S^3\setminus\FigEight)
  =
  \langle x,y\mid xy^{-1}x^{-1}yx=yxy^{-1}x^{-1}y\rangle,
\end{equation*}
where $x$ and $y$ are generators depicted in Figure~\ref{fig:fig8}.
Define
\begin{align*}
  \rho_{u,\pm}(x)
  &:=
  \begin{pmatrix}e^{u/2}&1\\0&e^{-u/2}\end{pmatrix},
  \\
  \rho_{u,\pm}(y)
  &:=
  \begin{pmatrix}e^{u/2}&0\\-d_{\pm}&e^{-u/2}\end{pmatrix}
\end{align*}
with
\begin{equation*}
  d_{\pm}:=
  \cosh{u}-\frac{3}{2}\pm\frac{1}{2}\sqrt{(2\cosh{u}+1)(2\cosh{u}-3)}.
\end{equation*}
It can be easily checked (you can use Mathematica for example) that $\rho_{u,\pm}$ defines a representation of $\pi_1\left(S^3\setminus\FigEight\right)$ to $\SL(2;\C)$.
\par
We choose a longitude (a loop in the knot complement that travels along the knot with linking number zero) $\lambda$ as a loop starting at the top right.
Then it presents an element $xy^{-1}xyx^{-2}yxy^{-1}x^{-1}\in\pi_1\left(S^3\setminus\FigEight\right)$.
It can be also checked that $\rho_{u,\pm}$ sends $\lambda$ to
\begin{equation*}
  \rho_{u,\pm}(\lambda)
  =
  \begin{pmatrix}
    \ell(u)^{\pm1}&\pm2\cosh(u/2)\sqrt{(2\cosh{u}+1)(2\cosh{u}-3)} \\
    0&\ell(u)^{\mp1}
  \end{pmatrix},
\end{equation*}
where
\begin{equation}\label{eq:longitude_fig8}
  \ell(u)
  :=
  \cosh(2u)-\cosh{u}-1+\sinh{u}\sqrt{(2\cosh{u}+1)(2\cosh{u}-3)}.
\end{equation}
Note that $\rho_{0,+}$ defines the complete hyperbolic structure on $S^3\setminus\FigEight$ \cite{Riley:MATPC75,Milnor:BULAM382,Thurston:GT3M}.
\subsection{Torus knots}\label{subsec:rep_torus_knot}
\begin{figure}[h!]
  \includegraphics[scale=0.3]{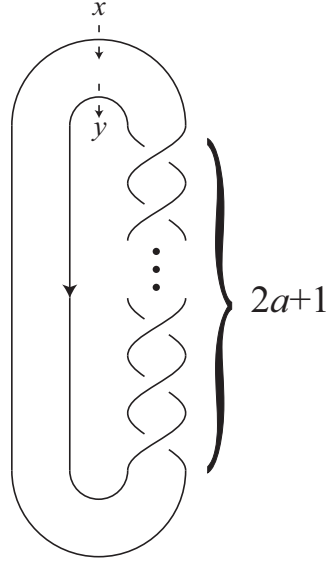}
  \caption{Generators of $\pi_1\left(S^3\setminus{T(2,2a+1)}\right)$}
  \label{fig:torus_knot}
\end{figure}

Let $T(2,2a+1)$ be the torus knot of type $(2,2a+1)$, where $a$ is a positive integer.
Then $\pi_1\left(S^3\setminus{T(2,2a+1)}\right)$ has the following presentation
\begin{equation}\label{eq:presentation_torus_knot}
  \pi_1\left(S^3\setminus{T(2,2a+1)}\right)
  =
  \langle
    x,y\mid(xy)^{a}x=y(xy)^{a}
  \rangle,
\end{equation}
where $x$ and $y$ are generators depicted in Figure~\ref{fig:torus_knot}.
For a complex number $u$ and $\omega_1$ with $\omega_1^{2a+1}=-1$, put
\begin{align*}
  \rho_{u,\omega_1}(x)
  &:=
  \begin{pmatrix}
    e^{u/2}&1 \\
    0      &e^{-u/2}
  \end{pmatrix},
  \\
  \rho_{u,\omega_1}(y)
  &:=
  \begin{pmatrix}
    e^{u/2}&0 \\
    \omega_1+\omega_1^{-1}-2\cosh{u}&e^{-u/2}
  \end{pmatrix}.
\end{align*}
Then $\rho_{u,\omega_1}$ is a representation from $\pi_1\left(S^3\setminus{T(2,2a+1)}\right)$ to $\SL(2;\C)$.
\par
If we choose the longitude $\lambda$ as a loop starting at the top right of Figure~\ref{fig:torus_knot} then it presents the element
\begin{equation*}
\begin{split}
  &(xy)^{a}x(xy)^{-a}\cdot(xy)^{a-1}x(xy)^{-a+1}\cdots(xy)x(xy)^{-1}\cdot x
  \\
  &\cdot(xy)^{a}x^{-1}(xy)^{-a+1}\cdot(xy)^{a-1}x^{-1}(xy)^{2-a}\cdots(xy)^2x^{-1}(xy)\cdot(xy)^{-1}x^{-1}
  \\
  &\cdot x^{-2a-1}
  \\
  =&
  (xy)^{a}xy^{-a}(xy)^{a}x^{-3a-1}.
\end{split}
\end{equation*}
Note that we need to add $x^{-3a-1}$ because the longitude is null-homologous, that is, its linking number with the knot is zero.
Since $(xy)^{a}x=y(xy)^a$ and $y^{-1}(xy)^{a}=(xy)^{a}x^{-1}$ from the relation \ref{eq:presentation_torus_knot}, we have
\begin{equation}\label{eq:longitude_torus_knot}
  \lambda
  =
  y(xy)^{2a}x^{-4a-1}.
\end{equation}
Its image by the representation is given by
\begin{equation}\label{eq:rep_longitude_torus_knot}
  \rho_{u,\omega_1}(\lambda)
  =
  \begin{pmatrix}
    -e^{-(2a+1)u} & \frac{\sinh\bigl((2a+1)u\bigr)}{\sinh(u/2)} \\
      0           & -e^{(2a+1)u}
  \end{pmatrix},
\end{equation}
which can be checked by Mathematica, for example.
\subsection{Twice-iterated torus knots}\label{subsec:rep_iterated_torus_knot}
\begin{figure}[h!]
  \begin{center}
    \includegraphics[scale=0.2]{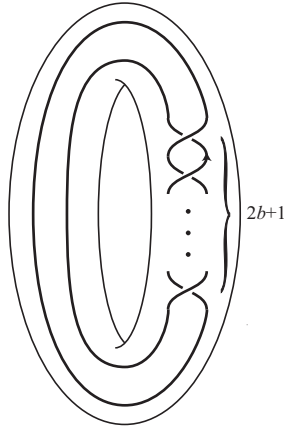}
  \end{center}
  \caption{A knot in the solid torus. There are $2b+1$ crossings.}
  \label{fig:pattern}
\end{figure}

\begin{figure}[h!]
  \includegraphics[scale=0.3]{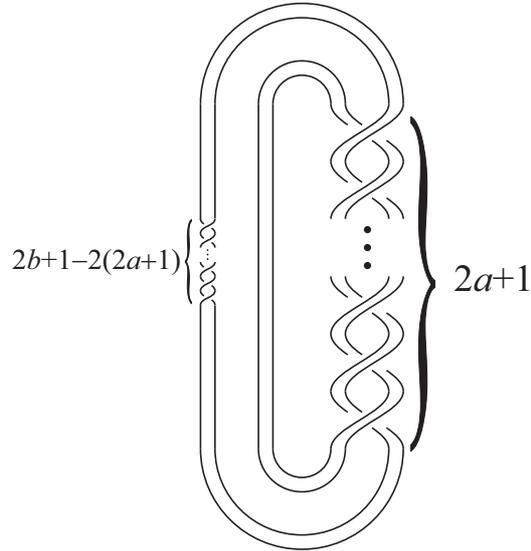}
  \caption{$T(2,2a+1)^{(2,2b+1)}$}
  \label{fig:iterated_torus_knot}
\end{figure}

Consider a solid torus $D:=D^2\times S^1$ with a knot, called a pattern knot, in it as shown in Fig\ref{fig:pattern}.
Here the knot goes twice along the solid torus and there are $2b+1$ positive crossings.
Now consider an embedding of $D$ in $S^3$ so that $D$ forms the torus knot $T(2,2a+1)$ and that the longitude of $D$ coincides with the longitude of $T(2,2a+1)$ (Figure~\ref{fig:iterated_torus_knot}).
Note that the number of positive crossings we need to add is $2b+1-2(2a+1)$, because we get $2(2a+1)$ twists from the original torus knot.
Then the knot in $D$ becomes a knot in $S^3$ called the $(2,2b+1)$-cable of $T(2,2a+1)$ denoted $T(2,2a+1)^{(2,2b+1)}$.
\par
Put $C:=S^3\setminus\Int{D}$, where $\Int{D}$ is the interior of the image of $D$ in $S^3$.
The boundary of $C$ is a torus $T_{C}$.
Let $P$ be the complement of (the interior of the regular neighborhood of) the pattern knot of $D$.
The boundary of $P$ consists of two tori $T_{P}$ and ${T'}_{P}$, where $T_{P}$ is the boundary of $D$ and ${T'}_{P}$ is the boundary of the regular neighborhood of the pattern knot.
Then $S^3\setminus\Int{N\left(T(2,2a+1)^{(2,2b+1)}\right)}=C\cup_{T_{C}=T_{P}}P$, where $N\left(T(2,2a+1)^{(2,2b+1)}\right)$ is the regular neighborhood of $T(2,2a+1)^{(2,2b+1)}$ in $S^3$.
\begin{figure}[h!]
  \begin{center}
    \includegraphics[scale=0.2]{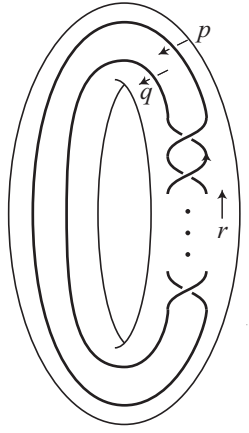}
  \end{center}
  \caption{Generators of $\pi_1(P)$}
  \label{fig:pattern_pi1}
\end{figure}
\par
Let us calculate $\pi_1(S^3\setminus{T(2,2a+1)^{(2,2b+1)}})$ by using van Kampen's theorem.
First of all from \eqref{eq:presentation_torus_knot} we have
\begin{equation}\label{eq:pi1_C}
  \pi_1(C)
  =
  \langle
    x,y\mid(xy)^ax=y(xy)^a
  \rangle.
\end{equation}
We also have
\begin{equation}\label{eq:pi1_P}
  \pi_1(P)
  =
  \langle
    p,q,r\mid(pq)r=r(pq),r(pq)^bp=q(pq)^br
  \rangle,
\end{equation}
where $p$, $q$ and $r$ are generators indicated in Figure~\ref{fig:pattern_pi1}.
Note that the longitude $\lambda_{C}$ of $T(2,2a+1)$ is given by
\begin{equation}\label{eq:longitude_C}
  \lambda_{\rm{C}}
  =
  y(xy)^{2a}x^{-4a-1}
\end{equation}
from \eqref{eq:longitude_torus_knot}.
Then by van Kampen's theorem, we have
\begin{equation*}
  \pi_1\left(S^3\setminus{T(2,2a+1)^{(2,2b+1)}}\right)
  =
  \pi_1(C)\ast_{\pi_1(T_{C})=\pi_1(T_{P})}\pi_1(P),
\end{equation*}
where we identify the meridian $x\in\pi_1(T_{C})$ with the meridian $pq\in\pi_1(T_{P})$, and the longitude $\lambda_{C}\in\pi_1(T_{C})$ with the longitude $r\in\pi_1(T_{P})$.
Therefore we have
\begin{equation}\label{eq:pi1_E1}
\begin{split}
  &\pi_1\left(S^3\setminus{T(2,2a+1)^{(2,2b+1)}}\right)
  \\
  =&
  \langle
    x,y,p,q,r
    \mid
    (xy)^ax=y(xy)^a,(pq)r=r(pq),r(pq)^bp=q(pq)^br,
  \\
  &\phantom{\langle x,y,p,q,r\mid}\quad
    x=pq,r=y(xy)^{2a}x^{-4a-1}
  \rangle
  \\
  =&
  \langle
    x,y,p,q
    \mid
    (xy)^ax=y(xy)^a,x\lambda_{C}=\lambda_{C}x,\lambda_{C}x^bp=qx^b\lambda_{C},x=pq
  \rangle,
\end{split}
\end{equation}
where we remove $r$ from the generator set introducing a word $\lambda_{C}=y(xy)^{2a}x^{-4a-1}$ in the relation set.
Since we have
\begin{equation*}
\begin{split}
  x\lambda_{C}x^{-1}\lambda^{-1}_{C}
  &=
  x\times y(xy)^{2a}x^{-4a-1}\times x^{-1}\times x^{4a+1}(xy)^{-2a}y^{-1}
  \\
  &=
  (xy)^{2a+1}x^{-1}(xy)^{-2a}y^{-1}
  \\
  &\text{(from the first relation in \eqref{eq:pi1_E1})}
  \\
  &=
  (xy)^{a+1}y^{-1}(xy)^{a}(xy)^{-2a}y^{-1}
  \\
  &=
  (xy)^{a+1}y^{-1}(xy)^{-a}y^{-1}
  \\
  &\text{(from the first relation in \eqref{eq:pi1_E1} again)}
  \\
  &=
  (xy)^{a+1}y^{-1}x^{-1}(xy)^{-a}
  \\
  &=1,
\end{split}
\end{equation*}
we do not need the relation $x\lambda_{C}=\lambda_{C}x$.
So we finally have
\begin{equation}\label{eq:pi1_E}
  \pi_1(E)
  =
  \langle
    x,y,p,q
    \mid
    (xy)^ax=y(xy)^a,\lambda_{C}x^bp=qx^b\lambda_{C},x=pq
  \rangle
\end{equation}
with $\lambda_{C}=y(xy)^{2a}x^{-4a-1}$.
\begin{rem}
We do not need $x$ in the generator set and $x=pq$ in the relation set.
However it would be convenient to include them in the following calculations.
\end{rem}
The longitude $\lambda$ of $T(2,2a+1)^{(2,2b+1)}$ is given by
\begin{equation}\label{eq:lambda}
\begin{split}
  \lambda
  &=
  r\cdot(pq)^{b}p(pq)^{-b}\cdot(pq)^{b-1}p(pq)^{-b+1}\cdots(pq)p(pq)^{-1}\cdot p
  \\
  &\quad
  \cdot r\cdot(pq)^{b}p^{-1}(pq)^{-b+1}\cdot(pq)^{b-1}p^{-1}(pq)^{2-b}\cdots(pq)^2x^{-1}(pq)\cdot(pq)^{-1}p^{-1}
  \\
  &\quad\cdot p^{-2b-1}
  \\
  &=
  r(pq)^{b}pq^{-b}r(pq)^{b}p^{-3b-1}
  \\
  &=
  \lambda_{C}x^bpq^{-b}\lambda_{C}x^bp^{-3b-1}.
\end{split}
\end{equation}
Here we read off $\lambda$ from the top right corner in Figure~\ref{fig:pattern_pi1} and multiply by $r$ any time we pass through the bottom.
\par
Let $\rho\colon\pi_1\left(S^3\setminus{T(2,2a+1)^{(2,2b+1)}}\right)\to\SL(2;\C)$ be a representation and put $\rho_{C}:=\rho\big|_{C}$ and $\rho_{P}:=\rho\big|_{P}$.
We consider the following three cases except for the case where $\Im\rho$ is Abelian.
(1) $\Im\rho_{C}$ is Abelian and $\Im\rho_{P}$ is non-Abelian,
(2) $\Im\rho_{C}$ is non-Abelian and $\Im\rho_{P}$ is Abelian, and
(3) both $\Im\rho_{C}$ and $\Im\rho_{P}$ are non-Abelian.
\subsubsection{$\Im\rho_{C}$ is Abelian and $\Im\rho_{P}$ is non-Abelian}
\label{subsubsec:AN}
If we put
\begin{equation}\label{eq:rep_AN}
\begin{split}
  \rho^{\rm{AN}}_{u,\omega_2}(p)
  &:=
  \begin{pmatrix}
    e^{u/2}&1 \\
       0   &e^{-u/2}
  \end{pmatrix},
  \\
  \rho^{\rm{AN}}_{u,\omega_2}(q)
  &:=
  \begin{pmatrix}
    e^{u/2}                         &0 \\
    \omega_2+\omega_2^{-1}-2\cosh(u)&e^{-u/2}
  \end{pmatrix},
  \\
  \rho^{\rm{AN}}_{u,\omega_2}(x)
  &=
  \rho^{\rm{AN}}_{u,\omega_2}(y)
  :=
  \rho^{\rm{AN}}_{u,\omega_2}(p)\rho^{\rm{AN}}_{u,\omega}(q)
\end{split}
\end{equation}
with $\omega_2^{2b+1}=-1$, then one can prove that $\rho^{\rm{AN}}_{u,\omega_2}$ defines a representation of $\pi_1\left(S^3\setminus{T(2,2a+1)^{(2,2b+1)}}\right)$ to $\SL(2;\C)$ from \eqref{eq:pi1_E}.
Note that
\begin{equation*}
  \rho^{\rm{AN}}_{u,\omega_2}(x)
  =
  \rho^{\rm{AN}}_{u,\omega_2}(y)
  =
  S_2^{-1}
  \begin{pmatrix}
    \omega_2&1\\
    0       &\omega_2^{-1}
  \end{pmatrix}
  S_2
\end{equation*}
with
\begin{equation*}
  S_2
  :=
  \begin{pmatrix}
    e^{u/2}                      &0 \\
    e^{u/2}\omega_2^{-1}-e^{-u/2}&1
  \end{pmatrix}.
\end{equation*}
The longitude $\lambda$ is sent to
\begin{equation}\label{eq:AN_longitude}
  \rho^{\rm{AN}}_{u,\omega_2}(\lambda)
  =
  \begin{pmatrix}
    -e^{-(2b+1)u} & \frac{\sinh((2b+1)u)}{\sinh(u/2)} \\
          0       & -e^{(2b+1)u}
  \end{pmatrix}
\end{equation}
from \eqref{eq:lambda}.
Note that $\rho^{\rm{AN}}_{u,\omega_2}\bigl|_{C}$ is Abelian and that $\lambda_{C}$ is sent to the identity matrix since $x$ and $y$ commute.
\subsubsection{$\Im\rho_{C}$ is non-Abelian and $\Im\rho_{P}$ is Abelian}
\label{subsubsec:NA}
Put
\begin{equation}\label{eq:rep_NA}
\begin{split}
  \rho^{\rm{NA}}_{u,\omega_1}(x)
  &:=
  \begin{pmatrix}e^{u}&1\\0&e^{-u}\end{pmatrix},
  \\
  \rho^{\rm{NA}}_{u,\omega_1}(y)
  &:=
  \begin{pmatrix}
     e^{u}                            &0\\
     \omega_1+\omega_1^{-1}-2\cosh(2u)&e^{-u}
  \end{pmatrix},
  \\
  \rho^{\rm{NA}}_{u,\omega_1}(p)
  &=
  \rho^{\rm{NA}}_{u,\omega_1}(q)
  :=
  \begin{pmatrix}e^{u/2}&\frac{1}{2\cosh(u/2)}\\0&e^{-u/2}\end{pmatrix},
\end{split}
\end{equation}
with $\omega_1^{2a+1}=-1$.
Then one can prove that $\rho^{\rm{NA}}_{u,\omega_1}$ defines a representation of $\pi_1\left(S^3\setminus{T(2,2a+1)^{(2,2b+1)}}\right)$ to $\SL(2;\C)$ from \eqref{eq:pi1_E}.
The longitude $\lambda_{C}$ of $T(2,2a+1)$ is sent to
\begin{equation*}
  \rho^{\rm{NA}}_{u,\omega_1}(\lambda_{C})
  =
  \begin{pmatrix}
    -e^{-2(2a+1)u}&\frac{\sinh(2(2a+1)u)}{\sinh(u)}\\
     0            &-e^{2(2a+1)u}
  \end{pmatrix}
\end{equation*}
and the longitude $\lambda$ of $T(2,2a+1)^{(2,2b+1)}$ is sent to
\begin{equation}\label{eq:NA_longitude}
  \rho^{\rm{NA}}_{u,\omega_1}(\lambda)
  =
  \rho^{\rm{NA}}_{u,\omega_1}(\lambda_{C}^2)
  =
  \begin{pmatrix}e^{-4(2a+1)u}&-\frac{\sinh(4(2a+1)u)}{\sinh(u)}\\0&e^{4(2a+1)u}\end{pmatrix}.
\end{equation}
\subsubsection{Both $\Im\rho_{C}$ and $\Im\rho_{P}$ are non-Abelian}
\label{subsubsec:NN}
If we put
\begin{equation}\label{eq:rep_NN}
\begin{split}
  \rho^{\rm{NN}}_{u,\omega_1,\omega_3}(p)
  &:=
  \begin{pmatrix}
    e^{u/2}&1 \\
       0   &e^{-u/2}
  \end{pmatrix},
  \\
  \rho^{\rm{NN}}_{u,\omega_1,\omega_3}(q)
  &:=
  \begin{pmatrix}
    e^{u/2}                         &0 \\
    \omega_3+\omega_3^{-1}-2\cosh(u)&e^{-u/2}
  \end{pmatrix},
  \\
  \rho^{\rm{NN}}_{u,\omega_1,\omega_3}(x)
  &:=
  \rho^{\rm{NN}}_{u,\omega_1,\omega_3}(pq)
  \\
  &=
  S_3^{-1}
  \begin{pmatrix}\omega_3&1\\0&\omega_3^{-1}\end{pmatrix}
  S_3,
  \\
  \rho^{\rm{NN}}_{u,\omega_1,\omega_3}(y)
  &:=
  S_3^{-1}
  \begin{pmatrix}
    \omega_3                                        & 0 \\
    \omega_1+\omega_1^{-1}-\omega_3^2-\omega_3^{-2} & \omega_3^{-1}
  \end{pmatrix}
  S_3
  \\
  &=
  \begin{pmatrix}
    \omega_3                                                                   & 0\\
    e^{-u/2}(\omega_3-\omega_3^{-1})-e^{u/2}\left(\omega_3^2+1-\omega_1-\omega_1^{-1}\right) &\omega_3^{-1}
  \end{pmatrix}
\end{split}
\end{equation}
with $\omega_1^{2a+1}=-1$, $\omega_3^{2b+1-4(2a+1)}=-1$, and
\begin{equation*}
  S_3
  :=
  \begin{pmatrix}
    e^{u/2}               &0\\
    e^{u/2}\omega_3^{-1}-e^{-u/2}&1
  \end{pmatrix},
\end{equation*}
then we have $\left(\rho^{\rm{NN}}_{u,\omega_1,\omega_3}(x)\rho^{\rm{NN}}_{u,\omega_1,\omega_3}(y)\right)^{a}\rho^{\rm{NN}}_{u,\omega_1,\omega_3}(x)=\rho^{\rm{NN}}_{u,\omega_1,\omega_3}(y)\left(\rho^{\rm{NN}}_{u,\omega_1,\omega_3}(x)\rho^{\rm{NN}}_{u,\omega_1,\omega_3}(y)\right)^a$.
The longitude $\lambda_{\rm{C}}$ of the companion knot is sent to
\begin{equation*}
\begin{split}
  &\rho^{\rm{NN}}_{u,\omega_1,\omega_3}(\lambda_{\rm{C}})
  \\
  =&
  S_3^{-1}
  \begin{pmatrix}
    -\omega_3^{-2(2a+1)}& \frac{\omega_3^{2(2a+1)}-\omega_3^{-2(2a+1)}}{\omega_3-\omega_3^{-1}} \\
     0                  & -\omega_3^{2(2a+1)}
  \end{pmatrix}
  S_3
  \\
  =&
  \frac{-\omega_3^{2(2a+1)}+\omega_3^{-2(2a+1)}}{\omega_3-\omega_3^{-1}}
  \begin{pmatrix}
    e^{-u}-
    \frac{\omega_3^{2(2a+1)-1}-\omega_3^{-2(2a+1)+1}}
         {\omega_3^{2(2a+1)}-\omega_3^{-2(2a+1)}}
    &-e^{-u/2} \\
    e^{-u/2}(e^u+e^{-u}-\omega_3-\omega_3^{-1})
    &\frac{\omega_3^{2(2a+1)+1}-\omega_3^{-2(2a+1)-1}}
          {\omega_3^{2(2a+1)}-\omega_3^{-2(2a+1)}}-e^{-u}
  \end{pmatrix}.
\end{split}
\end{equation*}
From the relation $\lambda_{C}x^bp-qx^b\lambda_{C}=O$, we have
\begin{equation*}
  \frac{\omega_3^{-b+4a+2}\left(\omega_3^{2b+1-4(2a+1)}+1\right)}{\omega_3+1}
  \begin{pmatrix}
    0                               &-1\\
    \omega_3+\omega_3^{-1}-2\cosh{u}&0
  \end{pmatrix}
  =O.
\end{equation*}
Since we choose $\omega_3$ such that $\omega_3^{2b+1-4(2a+1)}=-1$, $\lambda_{u,\omega_1,\omega_3}^{\rm{NN}}$ is well-defined.
\par
The longitude $\lambda$ is sent to
\begin{equation}\label{eq:NN_longitude}
  \rho^{\rm{NN}}_{u,\omega_1,\omega_3}(\lambda)
  =
  \begin{pmatrix}
    -e^{-(2b+1)u} & \frac{\sinh((2b+1)u)}{\sinh(u/2)}\\
         0        & -e^{(2b+1)u}
  \end{pmatrix}.
\end{equation}
\section{Chern--Simons invariant for a knot complement}
\label{sec:CS}
For the definition of the Chern--Simons invariant, refer to \cite{Chern/Simons:ANNMA21974,Kirk/Klassen:COMMP93}.
Here I only describe how to calculate it in the case of knot complements.
\subsection{How to calculate}
For a knot $K$, put $M:=S^3\setminus{K}$.
Note that $\partial{M}$ is a torus.
Let $X(M)$ be the $\SL(2;\C)$-character variety of $M$, that is, the set of characters of representations of $\pi_1(M)$ to $\SL(2;\C)$.
Then the Chern--Simons invariant $\cs_{M}$ is a map from $X(M)$ to $\Hom(\pi_1(\partial{M}),\C)\times\C^{\ast}/\approx$, where $\approx$ is defined as follows:
\begin{equation}\label{eq:CS_equivalence}
  \begin{cases}
    (s,t;z)&\approx(s+1,t;z\exp(-8\pi\sqrt{-1}t)), \\
    (s,t;z)&\approx(s,t+1;z\exp(8\pi\sqrt{-1}s)), \\
    (s,t;z)&\approx(-s,-t;z).
  \end{cases}
\end{equation}
Here we use the meridian $\mu$, a generator of $H_1(M)$ so that $\mu$ links with $K$ with linking number zero, and the preferred longitude $\lambda$, a loop on $\partial{M}$ that is parallel to $K$ such that $0=[\lambda]\in H_1(M)$, as a basis of $\pi_1(\partial{M})\cong\Z\oplus\Z$ and identify $(s,t)$ with $s\mu^{\ast}+t\lambda^{\ast}\in\Hom(\pi_1(\partial{M}),\C)$ with $\mu^{\ast}$ and $\lambda^{\ast}$ the dual elements of $\mu$ and $\lambda$, respectively.
\par
Now I will explain how to calculate $\cs_{M}$ following \cite{Kirk/Klassen:COMMP93}.
\begin{thm}[\cite{Kirk/Klassen:COMMP93}]\label{thm:Kirk/Klassen}
Let $\rho_t\colon\pi_1(M)\to\SL(2;\C)$ be a path of representations \rm{(}$0\le t\le1$\rm{)}.
Since $\mu$ and $\lambda$ commute we may assume that both $\rho_t(\mu)$ and $\rho_t(\lambda)$ are upper-tiangular.
We define $u_t$ and $v_t$ as follows:
\begin{equation*}
\begin{split}
  \rho_t(\mu)&=\begin{pmatrix}e^{u_t/2}&\ast\\0&e^{-u_t/2}\end{pmatrix},
  \\
  \rho_t(\lambda)&=\begin{pmatrix}e^{v_t/2}&\ast\\0&e^{-v_t/2}\end{pmatrix}.
\end{split}
\end{equation*}
Suppose that $\cs_M$ is given as
\begin{equation*}
  \cs_{M}([\rho_t])
  =
  \left[
  \frac{u_t}{4\pi\sqrt{-1}},\frac{v_t}{4\pi\sqrt{-1}};z_t\right].
\end{equation*}
Then we have
\begin{equation*}
  \frac{z_1}{z_0}
  =
  \exp
  \left[
    \frac{\sqrt{-1}}{2\pi}
    \int_{0}^{1}
    \left(
      u_t\frac{d\,v_t}{d\,t}-v_t\frac{d\,u_t}{d\,t}
    \right)
    dt
  \right].
\end{equation*}
\end{thm}
If a representation $\rho\colon\pi_1(M)\to\SL(2;\C)$ satisfies
\begin{equation*}
  \rho\colon\text{meridian}\to
  \begin{pmatrix}
    e^{u/2}&\ast \\
    0      &e^{-u/2}
  \end{pmatrix},
  \quad
  \text{longitude}\to
  \begin{pmatrix}
    e^{v/2}&\ast \\
    0      &e^{-v/2}
  \end{pmatrix},
\end{equation*}
then we define $\CS_{u,v}([\rho])$ so that
\vspace{-2mm}
\begin{equation*}
  \cs_{M}([\rho])
  =
  \left[
    \frac{u}{4\pi\sqrt{-1}},
    \frac{v}{4\pi\sqrt{-1}};
    \exp\left(\frac{2}{\pi\sqrt{-1}}\CS_{u,v}([\rho])\right)
  \right].
\end{equation*}
Note that $\CS_{u,v}([\rho])$ is defined modulo $\pi^2\Z$.
\subsection{$\SL(2;\C)$ Chern--Simons invariant of a hyperbolic knot}
Let $H$ be a hyperbolic knot, and $\rho_0\colon\pi_1(S^3\setminus{H})\to\SL(2;\C)$ be the representation associated with the complete hyperbolic structure.
We can deform the complete structure by a small complex parameter $u$.
Let $\rho_u$ be the representation associated with $u$.
By conjugation we assume
\begin{align*}
  \rho(\mu)&=\begin{pmatrix}e^{u/2}&\ast\\0&e^{-u/2}\end{pmatrix},
  \\
  \rho(\lambda)&=\begin{pmatrix}e^{v(u)/2}&\ast\\0&e^{-v(u)/2}\end{pmatrix},
\end{align*}
where $\mu$ is the meridian and $\lambda$ is the longitude of $\pi_1(S^3\setminus{H})$.
See for example \cite{Neumann/Zagier:TOPOL85}.
\par
Then we can define a path of representations $\rho_{tu}$ and by \cite{Kirk/Klassen:COMMP93} we can calculate the Chern--Simons invariant of $\rho_{u}$ as follows.
Let
\begin{equation*}
  \left[
    \frac{ut}{4\pi\sqrt{-1}},\frac{v(ut)}{4\pi\sqrt{-1}};z_t
  \right]
\end{equation*}
be the Chern-Simons invariant of the representation $\rho_{ut}$.
Then we have
\begin{equation*}
\begin{split}
  \frac{z_1}{z_0}
  &=
  \exp
  \left[
    -8\pi\sqrt{-1}
    \int_{0}^{1}
      \left(
        \frac{ut}{4\pi\sqrt{-1}}\frac{d}{dt}\frac{v(ut)}{4\pi\sqrt{-1}}
        -
        \frac{v(ut)}{4\pi\sqrt{-1}}\frac{d}{dt}\frac{ut}{4\pi\sqrt{-1}}
      \right)
    \,dt
  \right]
  \\
  &=
  \exp
  \left[
    \frac{\sqrt{-1}}{2\pi}
    \left(
      \Bigl[
        utv(ut)
      \Bigr]_0^1
      -
      2u
      \int_{0}^{1}v(ut)\,dt
    \right)
  \right]
  \\
  &=
  \exp
  \left[
    \frac{\sqrt{-1}}{2\pi}
    \left(
      uv(u)
      -
      2
      \int_{0}^{u}v(s)\,ds
    \right)
  \right].
\end{split}
\end{equation*}
Since $z_0=\exp\left[\frac{2}{\pi\sqrt{-1}}\CS(S^3\setminus{H})\right]$, where $\CS(S^3\setminus{H})$ is the $\SL(2;\C)$ Chern--Simons invariant associated with the Levi-Civita connection, we have
\begin{equation}\label{eq:CS}
  \CS_{u,v(u)}\left(\rho_u\right)
  =
  \CS(S^3\setminus{H})
  +
  \frac{1}{2}\int_{0}^{u}v(s)\,ds
  -\frac{1}{4}uv(u).
\end{equation}
\subsection{Chern--Simons invariant of the figure-eight knot}\label{subsec:CS_fig8}
Let $\rho_{u,\pm}$ be the representation of $\pi_1\left(S^3\setminus\FigEight\right)$ to $\SL(2;\C)$ defined in Subsection~\ref{subsec:rep_fig8}.
\par
Let $\rho_{t}:=\rho_{tu,+}$ ($0\le t\le1$) be a path of representations.
Note that $\rho_{0}=\rho_{0,+}$ and $\rho_{1}=\rho_{u,+}$.
Let
\begin{equation*}
  \left[
    \frac{ut}{4\pi\sqrt{-1}},\frac{v(t)}{4\pi\sqrt{-1}};z(t)
  \right]
\end{equation*}
be the Chern--Simons invariant for the representation $\rho_{t}$, where $v(t):=2\log\ell(ut)$.
From \eqref{eq:CS} we have
\begin{equation}\label{eq:figure8_CS}
\begin{split}
  \CS_{u,v}\left(\rho_{u,+}\right)
  &=
  \sqrt{-1}\Vol\left(S^3\setminus\FigEight\right)
  +
  \int_{0}^{u}\log\ell(s)\,ds
  -
  \frac{1}{2}u\log\ell(u).
\end{split}
\end{equation}
\subsection{Chern--Simons invariant of a torus knot}\label{subsec:CS_torus_knot}
Put $M:=S^3\setminus{T(2,2a+1)}$ and let $\rho_{u,\omega_1}\colon\pi_1(M)\to\SL(2;\C)$ be the representation defined in Subsection~\ref{subsec:rep_torus_knot}.
Note that $\omega_1^{2a+1}=-1$.
Put $\omega_1:=\exp\left(\frac{(2k+1)\pi\sqrt{-1}}{2a+1}\right)$ for $k=0,1,\dots,a-1$.
\par
Let $\alpha_{t}$ be a path of Abelian representations defined by
\begin{equation*}
  \alpha_{t}(x)=\alpha_{t}(y)
  :=
  \begin{pmatrix}
    \exp\left(\frac{(2k+1)\pi\sqrt{-1}}{2(2a+1)}t\right)&0  \\
      0                                                 &\exp\left(-\frac{(2k+1)\pi\sqrt{-1}}{2(2a+1)}t\right)
  \end{pmatrix}
\end{equation*}
and  $\beta_{t}$ be a path of non-Abelian representations defined by
\begin{equation*}
  \begin{cases}
    \beta_{t}(x)
    &:=
    \begin{pmatrix}
      e^{u_t/2}&1 \\
        0      &e^{-u_t/2}
    \end{pmatrix},
    \\
    \beta_{t}(y)
    &:=
    \begin{pmatrix}
      e^{u_t/2}                                           &0 \\
      2\cos\left(\frac{(2k+1)\pi}{2a+1}\right)-2\cosh{u_t}&e^{-u_t/2}
    \end{pmatrix},
  \end{cases}
\end{equation*}
where we put $u_t:=\dfrac{(1-t)(2k+1)\pi\sqrt{-1}}{2a+1}+tu$.
Note that
\begin{itemize}
\item
$\alpha_0$ is trivial, and so $\cs_M([\alpha_0])=1$,
\item
$\alpha_1$ and $\beta_0$ share the same trace because $\beta_0$ is upper-triangular and so $\cs_M([\alpha_1])=\cs_M([\beta_0])$, and
\item
$\beta_1=\rho_{u,\omega_1}$.
\end{itemize}
We regard $x$ as the meridian $\mu$.
From Theorem~\ref{thm:Kirk/Klassen} we can write
\begin{align*}
  \cs_M([\alpha_t])
  &:=
  \left[\frac{(2k+1)t}{4(2a+1)},0;w_t\right],
  \\
  \cs_M([\beta_t])
  &:=
  \left[\frac{u_t}{4\pi\sqrt{-1}},\frac{-2(2a+1)u_t+2\pi{l}\sqrt{-1}}{4\pi\sqrt{-1}};z_t\right]
\end{align*}
for an odd integer $l$, since
\begin{align*}
  \alpha_t(\lambda)
  &=
  \begin{pmatrix}1&0\\0&1\end{pmatrix},
  \\
  \beta_t(\lambda)
  &=
  \begin{pmatrix}e^{-(2a+1)u_t+l\pi\sqrt{-1}}&\frac{\sinh\bigl((2a+1)u_t\bigr)}{\sinh(u_t/2)}\\0&e^{(2a+1)u_t-l\pi\sqrt{-1}}\end{pmatrix}
\end{align*}
from \eqref{eq:rep_longitude_torus_knot}.
Then Kirk--Klassen's theorem (Theorem~\ref{thm:Kirk/Klassen}) shows that $\dfrac{w_1}{w_0}=1$ and
\begin{equation*}
\begin{split}
  \frac{z_1}{z_0}
  &=
  \exp
  \left(
    \frac{\sqrt{-1}}{2\pi}
    \int_{0}^{1}
    \left(
      (u_t\times\left(-2(2a+1)\frac{d\,u_t}{d\,t}\right)-(-2(2a+1)u_t+2\pi{l}\sqrt{-1})\times\frac{d\,u_t}{d\,t}
    \right)\,dt
  \right)
  \\
  &=
  \exp
  \left(l\Bigl[u_t\Bigr]_0^1\right)
  \\
  &=
  \exp\left(l\left(u-\frac{(2k+1)\pi\sqrt{-1}}{2a+1}\right)\right).
\end{split}
\end{equation*}
Since $\cs_M([\alpha_1])=\cs_M([\beta_0])$ and $w_1=w_0=1$, we have
\begin{equation*}
  \left[\frac{2k+1}{4(2a+1)},0;1\right]
  =
  \left[\frac{u_0}{4\pi\sqrt{-1}},\frac{-2(2a+1)u_0+2\pi{l}\sqrt{-1}}{4\pi\sqrt{-1}};z_0\right]
  =
  \left[\frac{2k+1}{4(2a+1)},\frac{l-2k-1}{2};z_0\right].
\end{equation*}
However, from the equivalence relation \eqref{eq:CS_equivalence}, we have
\begin{equation*}
  \left[\frac{2k+1}{4(2a+1)},0;1\right]
  \approx
  \left[\frac{2k+1}{4(2a+1)},\frac{l-2k-1}{2};\exp\left(\frac{(2k+1)(l-2k-1)\pi\sqrt{-1}}{2a+1}\right)\right].
\end{equation*}
So we have
\begin{equation*}
  z_0=\exp\left(\frac{(2k+1)(l-2k-1)\pi\sqrt{-1}}{2a+1}\right)
\end{equation*}
and
\begin{equation*}
  z_1
  =
  z_0\exp\left(l\left(u-\frac{(2k+1)\pi\sqrt{-1}}{2a+1}\right)\right)
  =
  \exp\left(lu-\frac{(2k+1)^2\pi\sqrt{-1}}{2a+1}\right).
\end{equation*}
Therefore we finally have
\begin{equation}\label{eq:torus_knot_CS}
  \CS_{u,v}([\rho_{u,\omega_1}])
  =
  \frac{1}{2}lu\pi\sqrt{-1}+\frac{(2k+1)^2\pi^2}{2(2a+1)}
\end{equation}
with $v:=-2(2a+1)u+2\pi{l}\sqrt{-1}$.
Note that this depends on the choice of an odd integer $l$.
\subsection{Chern--Simons invariant of a twice-iterated torus knot}
I will calculate the Chern--Simons invariant of $T(2,2a+1)^{(2,2b+1)}$ associated with non-Abelian representations defined in Subsection~\ref{subsec:rep_iterated_torus_knot} in a similar way to the case of $T(2,2a+1)$.
Throughout this subsection we put $M:=S^3\setminus{T(2,2a+1)^{(2,2b+1)}}$
\subsubsection{$\Im\rho_{C}$ is Abelian and $\Im\rho_{P}$ is non-Abelian}
Let $\rho^{\rm{AN}}_{u,\omega_2}$ be the representation defined in Sub-subsection~\ref{subsubsec:AN}.
Note that $\omega_2^{2b+1}=-1$.
Put $\omega_2:=\exp\left(\frac{(2j+1)\pi\sqrt{-1}}{2b+1}\right)$ for $j=0,1,\dots,b-1$.
\par
Let $\alpha_{t}$ be a path of Abelian representations defined by
\begin{align*}
  \alpha_{t}(p)&=\alpha_{t}(q)
  :=
  \begin{pmatrix}
    \exp\left(\frac{(2j+1)\pi\sqrt{-1}}{2(2b+1)}t\right)&0  \\
      0                                                 &\exp\left(-\frac{(2j+1)\pi\sqrt{-1}}{2(2b+1)}t\right)
  \end{pmatrix},
  \\
  \alpha_{t}(x)&=\alpha_{t}(y)
  :=
  \begin{pmatrix}
    \exp\left(\frac{(2j+1)\pi\sqrt{-1}}{2b+1}t\right)&0  \\
      0                                              &\exp\left(-\frac{(2j+1)\pi\sqrt{-1}}{2b+1}t\right)
  \end{pmatrix}.
\end{align*}
For $u_t:=\dfrac{(1-t)(2j+1)\pi\sqrt{-1}}{2b+1}+tu$, let $\beta_{t}$ be a path of representations $\rho^{\rm{AN}}_{u_t,\omega_2}$, that is, we define
\begin{equation*}
  \begin{cases}
    \beta_{t}(p)
    &:=
    \begin{pmatrix}
      e^{u_t/2}&1 \\
        0    &e^{-u_t/2}
    \end{pmatrix},
    \\
    \beta_{t}(q)
    &:=
    \begin{pmatrix}
      e^{u_t/2}                                           &0 \\
      2\cos\left(\frac{(2j+1)\pi}{2b+1}\right)-2\cosh{u_t}&e^{-u_t/2}
    \end{pmatrix},
    \\
    \beta_{t}(x)
    &=
    \beta_{t}(y)
    :=
    \beta_{t}(p)\beta_{t}(q).
  \end{cases}
\end{equation*}
Note that
\begin{itemize}
\item
$\alpha_0$ is trivial, and so $\cs_M([\alpha_0])=1$,
\item
$\alpha_1$ and $\beta_0$ share the same trace because $\beta_0$ is upper-triangular and so $\cs_M([\alpha_1])=\cs_M([\beta_0])$, and
\item
$\beta_1=\rho^{\rm{AN}}_{u,\omega_2}$.
\end{itemize}
We regard $p$ as the meridian $\mu$.
From Theorem~\ref{thm:Kirk/Klassen} we can write
\begin{align*}
  \cs_M([\alpha_t])
  &:=
  \left[\frac{(2j+1)t}{4(2b+1)},0;w_t\right],
  \\
  \cs_M([\beta_t])
  &:=
  \left[\frac{u_t}{4\pi\sqrt{-1}},\frac{-2(2b+1)u_t+2\pi{m}\sqrt{-1}}{4\pi\sqrt{-1}};z_t\right]
\end{align*}
for an odd integer $m$, since
\begin{align*}
  \alpha_t(\lambda)
  &=
  \begin{pmatrix}1&0\\0&1\end{pmatrix},
  \\
  \beta_t(\lambda)
  &=
  \begin{pmatrix}
    e^{-(2b+1)u_t+m\pi\sqrt{-1}}&\frac{\sinh\bigl((2b+1)u_t\bigr)}{\sinh(u_t/2)}\\
    0                           &e^{(2b+1)u_t-m\pi\sqrt{-1}}\end{pmatrix}
\end{align*}
from \eqref{eq:AN_longitude}.
Then Kirk--Klassen's theorem (Theorem~\ref{thm:Kirk/Klassen}) shows that $\dfrac{w_1}{w_0}=1$ and
\begin{equation*}
\begin{split}
  \frac{z_1}{z_0}
  &=
  \exp
  \left(
    \frac{\sqrt{-1}}{2\pi}
    \int_{0}^{1}
    \left(
      (u_t\times\left(-2(2b+1)\frac{d\,u_t}{d\,t}\right)-(-2(2b+1)u_t+2\pi{m}\sqrt{-1})\times\frac{d\,u_t}{d\,t}
    \right)\,dt
  \right)
  \\
  &=
  \exp
  \left(m\Bigl[u_t\Bigr]_0^1\right)
  \\
  &=
  \exp\left(m\left(u-\frac{(2j+1)\pi\sqrt{-1}}{2b+1}\right)\right).
\end{split}
\end{equation*}
Since $\cs_M([\alpha_1])=\cs_M([\beta_0])$ and $w_1=w_0=1$, we have
\begin{equation*}
  \left[\frac{2j+1}{4(2b+1)},0;1\right]
  =
  \left[\frac{u_0}{4\pi\sqrt{-1}},\frac{-2(2b+1)u_0+2\pi{m}\sqrt{-1}}{4\pi\sqrt{-1}};z_0\right]
  =
  \left[\frac{2j+1}{4(2b+1)},\frac{m-2j-1}{2};z_0\right].
\end{equation*}
However, from the equivalence relation \eqref{eq:CS_equivalence}, we have
\begin{equation*}
  \left[\frac{2j+1}{4(2b+1)},0;1\right]
  \approx
  \left[\frac{2j+1}{4(2b+1)},\frac{m-2j-1}{2};\exp\left(\frac{(2j+1)(m-2j-1)\pi\sqrt{-1}}{2b+1}\right)\right].
\end{equation*}
So we have
\begin{equation*}
  z_0=\exp\left(\frac{(2j+1)(m-2j-1)\pi\sqrt{-1}}{2b+1}\right)
\end{equation*}
and
\begin{equation*}
  z_1
  =
  z_0\exp\left(m\left(u-\frac{(2j+1)\pi\sqrt{-1}}{2b+1}\right)\right)
  =
  \exp\left(mu-\frac{(2j+1)^2\pi\sqrt{-1}}{2b+1}\right).
\end{equation*}
Therefore we finally have
\begin{equation}\label{eq:AN_CS}
  \CS_{u,v}([\rho^{\rm{AN}}_{u,\omega_2}])
  =
  \frac{1}{2}mu\pi\sqrt{-1}+\frac{(2j+1)^2\pi^2}{2(2b+1)}
\end{equation}
with $v:=-2(2b+1)u+2\pi{m}\sqrt{-1}$.
Note that this depends on the choice of an odd integer $m$ and that the result here can be obtained from the result for the torus knot $T(2,2a+1)$ in Subsection~\ref{subsec:CS_torus_knot} by replacing $a$ with $b$, $k$ with $j$ and $l$ with $m$.
\subsubsection{$\Im\rho_{C}$ is non-Abelian and $\Im\rho_{P}$ is Abelian}
Let $\rho^{\rm{NA}}_{u,\omega_1}$ be the representation defined in Sub-subsection~\ref{subsubsec:NA}.
Note that $\omega_1^{2b+1}=-1$.
Put $\omega_1:=\exp\left(\frac{(2k+1)\pi\sqrt{-1}}{2a+1}\right)$ for $k=0,1,\dots,a-1$.
\par
Let $\alpha_{t}$ be a path of Abelian representations defined by
\begin{align*}
  \alpha_{t}(p)&=\alpha_{t}(q)
  :=
  \begin{pmatrix}
    \exp\left(\frac{(2k+1)\pi\sqrt{-1}}{4(2a+1)}t\right)&0  \\
      0                                                 &\exp\left(-\frac{(2k+1)\pi\sqrt{-1}}{4(2a+1)}t\right)
  \end{pmatrix},
  \\
  \alpha_{t}(x)&=\alpha_{t}(y)
  :=
  \begin{pmatrix}
    \exp\left(\frac{(2k+1)\pi\sqrt{-1}}{2(2a+1)}t\right)&0  \\
      0                                                 &\exp\left(-\frac{(2k+1)\pi\sqrt{-1}}{2(2a+1)}t\right)
  \end{pmatrix}.
\end{align*}
For $u_t:=\dfrac{(1-t)(2k+1)\pi\sqrt{-1}}{2(2a+1)}+tu$, let $\beta_{t}$ be a path of representations $\rho^{\rm{NA}}_{u_t,\omega_1}$, that is, we define
\begin{equation*}
  \begin{cases}
    \beta_{t}(x)
    &:=
    \begin{pmatrix}
      e^{u_t}&1 \\
        0    &e^{-u_t}
    \end{pmatrix},
    \\
    \beta_{t}(y)
    &:=
    \begin{pmatrix}
      e^{u_t}                                              &0 \\
      2\cos\left(\frac{(2k+1)\pi}{2a+1}\right)-2\cosh(2u_t)&e^{-u_t}
    \end{pmatrix},
    \\
    \beta_{t}(p)
    &=
    \beta_{t}(q)
    :=
    \begin{pmatrix}
      e^{u_t/2}&\frac{1}{2\cosh(u_t/2)} \\
      0        &e^{-u_t/2}
    \end{pmatrix}.
  \end{cases}
\end{equation*}
Note that
\begin{itemize}
\item
$\alpha_0$ is trivial, and so $\cs_M([\alpha_0])=1$,
\item
$\alpha_1$ and $\beta_0$ share the same trace because $\beta_0$ is upper-triangular and so $\cs_M([\alpha_1])=\cs_M([\beta_0])$, and
\item
$\beta_1=\rho^{\rm{NA}}_{u,\omega_1}$.
\end{itemize}
We regard $p$ as the meridian $\mu$.
From Theorem~\ref{thm:Kirk/Klassen} we can write
\begin{align*}
  \cs_M([\alpha_t])
  &:=
  \left[\frac{(2k+1)t}{8(2a+1)},0;w_t\right],
  \\
  \cs_M([\beta_t])
  &:=
  \left[\frac{u_t}{4\pi\sqrt{-1}},\frac{-8(2a+1)u_t+4\pi{l}\sqrt{-1}}{4\pi\sqrt{-1}};z_t\right]
\end{align*}
for an integer $l$, since
\begin{align*}
  \alpha_t(\lambda)
  &=
  \begin{pmatrix}1&0\\0&1\end{pmatrix},
  \\
  \beta_t(\lambda)
  &=
  \begin{pmatrix}
    e^{-4(2a+1)u_t+2l\pi\sqrt{-1}}&-\frac{\sinh\bigl(4(2a+1)u_t\bigr)}{\sinh(u_t)}\\
    0                             &e^{4(2a+1)u_t-2l\pi\sqrt{-1}}\end{pmatrix}
\end{align*}
from \eqref{eq:NA_longitude}.
Then Kirk--Klassen's theorem (Theorem~\ref{thm:Kirk/Klassen}) shows that $\dfrac{w_1}{w_0}=1$ and
\begin{equation*}
\begin{split}
  \frac{z_1}{z_0}
  &=
  \exp
  \left(
    \frac{\sqrt{-1}}{2\pi}
    \int_{0}^{1}
    \left(
      (u_t\times\left(-8(2a+1)\frac{d\,u_t}{d\,t}\right)-(-8(2a+1)u_t+4\pi{l}\sqrt{-1})\times\frac{d\,u_t}{d\,t}
    \right)\,dt
  \right)
  \\
  &=
  \exp
  \left(2l\Bigl[u_t\Bigr]_0^1\right)
  \\
  &=
  \exp\left(2l\left(u-\frac{(2k+1)\pi\sqrt{-1}}{2(2a+1)}\right)\right).
\end{split}
\end{equation*}
Since $\cs_M([\alpha_1])=\cs_M([\beta_0])$ and $w_1=w_0=1$, we have
\begin{equation*}
  \left[\frac{2k+1}{8(2a+1)},0;1\right]
  =
  \left[\frac{u_0}{4\pi\sqrt{-1}},\frac{-8(2a+1)u_0+4\pi{l}\sqrt{-1}}{4\pi\sqrt{-1}};z_0\right]
  =
  \left[\frac{2k+1}{8(2a+1)},l-2k-1;z_0\right].
\end{equation*}
However, from the equivalence relation \eqref{eq:CS_equivalence}, we have
\begin{equation*}
  \left[\frac{2k+1}{8(2a+1)},0;1\right]
  \approx
  \left[\frac{2k+1}{8(2b+1)},l-2k-1;\exp\left(\frac{(2k+1)(l-2k-1)\pi\sqrt{-1}}{2a+1}\right)\right].
\end{equation*}
So we have
\begin{equation*}
  z_0=\exp\left(\frac{(2k+1)(l-2k-1)\pi\sqrt{-1}}{2a+1}\right)
\end{equation*}
and
\begin{equation*}
  z_1
  =
  z_0\exp\left(2l\left(u-\frac{(2k+1)\pi\sqrt{-1}}{2(2a+1)}\right)\right)
  =
  \exp\left(2lu-\frac{(2k+1)^2\pi\sqrt{-1}}{2a+1}\right).
\end{equation*}
Therefore we finally have
\begin{equation}\label{eq:NA_CS}
  \CS_{u,v}([\rho^{\rm{NA}}_{u,\omega_2}])
  =
  lu\pi\sqrt{-1}+\frac{(2k+1)^2\pi^2}{2(2a+1)}
\end{equation}
with $v:=-8(2a+1)u+4\pi{l}\sqrt{-1}$.
Note that this depends on the choice of an integer $l$.
\subsubsection{Both $\Im\rho_{C}$ and $\Im\rho_{P}$ are non-Abelian}
Let $\rho^{\rm{NN}}_{u,\omega_1,\omega_3}$ be the representation defined in Sub-subsection~\ref{subsubsec:NA}.
Note that $\omega_1^{2a+1}=\omega_3^{2b+1-4(2a+1)}-1$.
Put $\omega_1:=\exp\left(\frac{(2k+1)\pi\sqrt{-1}}{2a+1}\right)$ and $\omega_3:=\exp\left(\frac{(2h+1)}{2b+1-4(2a+1)}\pi\sqrt{-1}\right)$ for $k=0,1,\dots,a-1$ and $n=0,1,\dots,b+4a-2$.
\par
Put $u_t:=\dfrac{(1-t)(2h+1)\pi\sqrt{-1}}{2b+1-4(2a+1)}+tu$ and consider a path of representations $\beta_{t}:=\rho^{\rm{NN}}_{u_t,\omega_1,\omega_3}$.
Then $\beta_{1}=\rho^{\rm{NN}}_{u,\omega_1,\omega_3}$ and
\begin{equation*}
\begin{split}
  \beta_{0}(p)
  &=
  \begin{pmatrix}
    \omega_3^{1/2}&1 \\
       0          &\omega_3^{-1/2}
  \end{pmatrix}
  =
  S_3^{-1}
  \begin{pmatrix}
    \omega_3^{1/2}&\omega_3^{1/2} \\
    0             &\omega_3^{-1/2}
  \end{pmatrix}
  S_3,
  \\
  \beta_{0}(q)
  &=
  \begin{pmatrix}
    \omega_3^{1/2}&0 \\
    0             &\omega_3^{-1/2}
  \end{pmatrix}
  =
  S_3^{-1}
  \begin{pmatrix}
    \omega_3^{1/2}&0 \\
    0             &\omega_3^{-1/2}
  \end{pmatrix}
  S_3,
  \\
  \beta_{0}(x)
  &=
  S_3^{-1}
  \begin{pmatrix}
    \omega_3&1\\
    0       &\omega_3^{-1}
  \end{pmatrix}
  S_3,
  \\
  \beta_{0}(y)
  &=
  S_3^{-1}
  \begin{pmatrix}
    \omega_3                                       &0 \\
    \omega_1+\omega_1^{-1}-\omega_3^2-\omega_3^{-2}&\omega_3^{-1}
  \end{pmatrix}
  S_3
\end{split}
\end{equation*}
since $e^{u_0}=\omega_3$, where $S_3=\begin{pmatrix}\omega^{1/2}&0\\0&1\end{pmatrix}$.
\par
Recall that the knot complement $M=S^3\setminus{T(2,2a+1)^{(2,2b+1)}}$ can be decomposed as $M=C\cup_{T_C=T_P}P$, where $C$ is the complement of the torus knot $T(2,2a+1)$ and $P$ is the complement the pattern knot in the solid torus (see Subsection~\ref{subsec:rep_iterated_torus_knot}).
One can see that $\tr\beta_{0}\big|_C=\tr\rho^{\rm{NA}}_{u_0,\omega_1}\big|_C$ and $\tr\beta_{0}\big|_P=\tr\rho^{\rm{NA}}_{u_0,\omega_1}\big|_{P}$.
From the gluing formula of the Chern--Simons invariant \cite[Theorem~2.1]{Kirk/Klassen:COMMP93} we have
\begin{equation*}
  \cs_{M}\left(\left[\beta_{0}\right]\right)
  =
  \cs_{M}\left(\left[\rho^{\rm{NA}}_{u_0,\omega_1}\right]\right).
\end{equation*}
\par
We regard $p$ as the meridian $\mu$.
From Theorem~\ref{thm:Kirk/Klassen} we can write
\begin{equation*}
  \\
  \cs_M([\beta_t])
  :=
  \left[\frac{u_t}{4\pi\sqrt{-1}},\frac{-2(2b+1)u_t+2\pi{n}\sqrt{-1}}{4\pi\sqrt{-1}};z_t\right]
\end{equation*}
for an odd integer $n$, since
\begin{equation*}
  \beta_t(\lambda)
  =
  \begin{pmatrix}
    e^{-(2b+1)u_t+n\pi\sqrt{-1}}&\frac{\sinh\bigl((2b+1)u_t\bigr)}{\sinh(u_t/2)}\\
    0                            &e^{(2b+1)u_t-\pi\sqrt{-1}}\end{pmatrix}
\end{equation*}
from \eqref{eq:NN_longitude}.
Then Kirk--Klassen's theorem (Theorem~\ref{thm:Kirk/Klassen}) shows that
\begin{equation*}
\begin{split}
  \frac{z_1}{z_0}
  &=
  \exp
  \left(
    \frac{\sqrt{-1}}{2\pi}
    \int_{0}^{1}
    \left(
      (u_t\times\left(-2(2b+1)\frac{d\,u_t}{d\,t}\right)-(-2(2b+1)u_t+2\pi{n}\sqrt{-1})\times\frac{d\,u_t}{d\,t}
    \right)\,dt
  \right)
  \\
  &=
  \exp
  \left(n\Bigl[u_t\Bigr]_0^1\right)
  \\
  &=
  \exp\left(n\left(u-\frac{(2h+1)\pi\sqrt{-1}}{2b+1-4(2a+1)}\right)\right).
\end{split}
\end{equation*}
Since $\cs_M\left(\left[\rho^{\rm{NA}}_{u_0,\omega_1}\right]\right)=\cs_M\left(\left[\beta_0\right]\right)$, we have
\begin{equation*}
\begin{split}
  &\left[
    \frac{u_0}{4\pi\sqrt{-1}},
    \frac{-8(2a+1)u_0+4\pi{l}\sqrt{-1}}{4\pi\sqrt{-1}};
    \exp\left(2lu_0-\frac{(2k+1)^2\pi\sqrt{-1}}{2a+1}\right)
  \right]
  \\
  =&
  \left[
    \frac{u_0}{4\pi\sqrt{-1}},\frac{-2(2b+1)u_0+2\pi{n}\sqrt{-1}}{4\pi\sqrt{-1}};z_{0}
  \right].
\end{split}
\end{equation*}
from \eqref{eq:NA_CS}
However, from the equivalence relation \eqref{eq:CS_equivalence}, we have
\begin{equation*}
\begin{split}
  &\left[
    \frac{u_0}{4\pi\sqrt{-1}},
    \frac{-8(2a+1)u_0+4\pi{l}\sqrt{-1}}{4\pi\sqrt{-1}};
    \exp\left(2lu_0-\frac{(2k+1)^2\pi\sqrt{-1}}{2a+1}\right)
  \right]
  \\
  \approx&
  \left[
    \frac{u_0}{4\pi\sqrt{-1}},
    \frac{-2(2b+1)u_0+2\pi{n}\sqrt{-1}}{4\pi\sqrt{-1}};
    \exp\left(2lu_0-\frac{(2k+1)^2\pi\sqrt{-1}}{2a+1}+\frac{(n-2l-2h-1)(2h+1)\pi\sqrt{-1}}{2b+1-4(2a+1)}\right)
  \right].
\end{split}
\end{equation*}
Since $u_0=\dfrac{(2h+1)\pi\sqrt{-1}}{2b+1-4(2a+1)}$, we have
\begin{equation*}
  z_0
  =
  \exp
  \left[
    \left(
      \frac{2l(2h+1)}{2b+1-4(2a+1)}-\frac{(2k+1)^2}{2a+1}+\frac{(n-2l-2h-1)(2h+1)}{2b+1-4(2a+1)}
    \right)
    \pi\sqrt{-1}
  \right]
\end{equation*}
and
\begin{equation*}
\begin{split}
  z_1
  &=
  z_0\exp\left(n\left(u-\frac{(2h+1)\pi\sqrt{-1}}{2b+1-4(2a+1)}\right)\right)
  \\
  &=
  \exp
  \left(
    nu
    -\frac{(2k+1)^2\pi\sqrt{-1}}{2a+1}
    -\frac{(2h+1)^2\pi\sqrt{-1}}{2b+1-4(2a+1)}
  \right).
\end{split}
\end{equation*}
Therefore we finally have
\begin{equation}\label{eq:NN_CS}
  \CS_{u,v}([\rho^{\rm{NN}}_{u,\omega_1,\omega_3}])
  =
  \frac{1}{2}
  nu\pi\sqrt{-1}
  +\frac{(2k+1)^2\pi^2}{2(2a+1)}
  +\frac{(2h+1)^2\pi^2}{2(2b+1-4(2a+1))}.
\end{equation}
with $v:=-2(2b+1)u+2\pi{n}\sqrt{-1}$.
Note that this depends on the choice of an odd integer $n$.
\section{Twisted $\SL(2;\C)$ Reidemeister torsion for a knot complement}
\label{sec:Reidemeister}
In this section we study the Reidemeister torsion twisted by a representation.
It is defined as the torsion of a certain chain complex.
If the chain complex is acyclic then the torsion is well-defined without specifying a basis of the corresponding homology group.
Unfortunately in our case the homology is non-trivial, and so we need to choose a basis.
In the following subsection I will start with the definition of the torsion, and then describe how to choose a basis.
In the later subsections I will calculate the Reidemeister torsion for the figure-eight knot, torus knots, and twice-iterated torus knots.
\subsection{Definition}\label{subsec:Reidemeister_definition}
Let $\rho\colon\pi_1(S^3\setminus{K})\to\SL(2;\C)$ be a representation.
Let $\langle x_1,x_2,\dots,x_n\mid r_1,r_2,\dots,r_{n-1}\rangle$ a Wirtinger presentation of $\pi_1(S^3\setminus{K})$.
Put $\Pi:=\pi_1(S^3\setminus{K})$.
\par
For the universal cover $\widetilde{S^3\setminus{K}}$ of $S^3\setminus{K}$, the chain complex $C_{\ast}\left(\widetilde{S^3\setminus{K}};\Z\right)$ can be regarded as a $\Z[\Pi]$-module by the deck transformation, and $\sl_2(\C)$ can also be regarded as a $\Z[\Pi]$-module by $\Ad{\rho(x)}$ for $x\in\Pi$.
Here we define the adjoint action $\Ad{X}$ of $X\in\SL_2(\C)$ by $\Ad{X}(g):=X^{-1}gX$  for $g\in\sl_2(\C)$.
Then we have the following chain complex:
\begin{equation*}
  C_2\left(\widetilde{S^3\setminus{K}}\right)\otimes_{\Z[\Pi]}\sl_2(\C)
  \xrightarrow{\partial_2}
  C_1\left(\widetilde{S^3\setminus{K}}\right)\otimes_{\Z[\Pi]}\sl_2(\C)
  \xrightarrow{\partial_1}
  C_0\left(\widetilde{S^3\setminus{K}}\right)\otimes_{\Z[\Pi]}\sl_2(\C).
\end{equation*}
The associated homology group is denoted by $H_{\ast}(S^3\setminus{K};\rho)$.
\par
\begin{figure}[h]
  \includegraphics[scale=0.3]{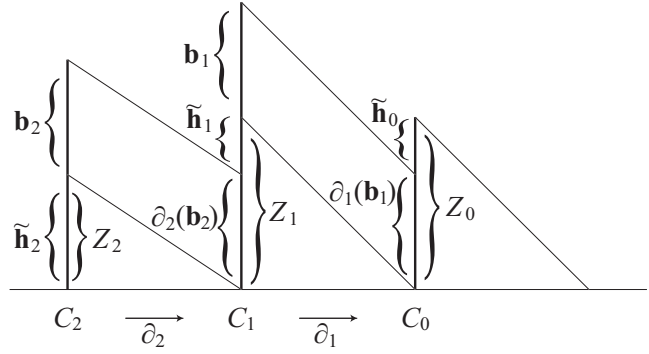}
  \caption{Chain complex and its basis}
  \label{fig:chain_complex}
\end{figure}

Let $\mathbf{c}_{i}:=\{c_{i,1},c_{i,2},\dots,c_{i,l_i}\}$ be a basis of $C_{i}:=C_{i}\left(\widetilde{S^3\setminus{K}}\right)\otimes_{\Z[\Pi]}\sl_2(\C)$, $\mathbf{b}_{i}:=\{b_{i,1},b_{i,2},\dots,b_{i,m_i}\}$ be a set of vectors such that $\{\partial_{i}(b_{i,1}),\partial_{i}(b_{i,2}),\dots,\partial_{i}(b_{i,m_{i}})\}$ forms a basis of $B_{i-1}:=\Im\partial_{i}$, $\mathbf{h}_i:=\{h_{i,1},h_{i,2},\dots,h_{i,n_i}\}$ be a basis of $H_{i}:=H_{i}(S^3\setminus{K};\rho)$, $\tilde{h}_{i,k}$ be a lift of $h_{i,k}$ in $Z_{i}:=\Ker\partial_{i}$, and $\tilde{\mathbf{h}}_{i}:=\{\tilde{h}_{i,1},\tilde{h}_{i,2},\dots,\tilde{h}_{i,n_i}\}$.
Then $\partial_{i+1}(\mathbf{b}_{i+1})\cup\tilde{\mathbf{h}}_{i}\cup\mathbf{b}_{i}$ forms a basis of $C_{i}$ (Figure~\ref{fig:chain_complex}).
For two bases $\mathbf{u}$ and $\mathbf{v}$ of a vector space $W$, let $[\mathbf{u}\mid\mathbf{v}]$ be the determinant of the basis change matrix from $\mathbf{u}$ to $\mathbf{v}$.
Put
\begin{equation}\label{eq:Reidemeister_def}
\begin{split}
  \operatorname{Tor}(C_{\ast},\mathbf{b}_{\ast},\mathbf{h}_{\ast})
  &:=
  \prod_{i=0}^{2}
  \left[
    \partial_{i+1}(\mathbf{b}_{i+1})\cup\tilde{\mathbf{h}}_{i}\cup\mathbf{b}_{i}\Bigm|\mathbf{c}_{i}
  \right]^{(-1)^{i+1}}
  \\
  &=
  \frac{
  \left[
    \partial_{2}(\mathbf{b}_{2})\cup\tilde{\mathbf{h}}_{1}\cup\mathbf{b}_{1}\Bigm|\mathbf{c}_{1}
  \right]}
  {
  \left[
    \partial_{1}(\mathbf{b}_{1})\cup\tilde{\mathbf{h}}_{0}\Bigm|\mathbf{c}_{0}
  \right]
  \left[
    \tilde{\mathbf{h}}_{2}\cup\mathbf{b}_{2}\Bigm|\mathbf{c}_{2}
  \right]}
  \in\C^{\ast}.
\end{split}
\end{equation}
Note that this does not depend on the choices of $\mathbf{b}_{i}$ nor the choices of lifts of $\mathbf{h}_{i}$ (see for example \cite{Turaev:2001}).
This does not depend on the choice of a basis of $\sl_2(\C)$, either, since the Euler characteristic of a knot complement is zero.
\par
For actual computation, we regard $S^3\setminus{K}$ as a CW-complex $\mathcal{K}$ with one $0$-cell $p$, $n$ $1$-cells $x_1,x_2,\dots,x_{n}$, and $(n-1)$ $2$-cells $r_1,r_2,\dots,r_{n-1}$.
For the following argument, I refer the reader to \cite[Chapter~11]{Lickorish:1997}.
\par
As a $\Z[\Pi]$-module $C_0\left(\tilde{\mathcal{K}}\right)$ is generated by a point $\tilde{p}$, $C_1\left(\tilde{\mathcal{K}}\right)$ is generated by $\tilde{x}_1,\tilde{x}_2,\dots,\tilde{x}_n$, where $\tilde{x}_i$ is a lift of $x_i$ attached to $\tilde{p}$, and $C_2\left(\tilde{\mathcal{K}}\right)$ is generated by $\tilde{r}_1,\tilde{r}_2,\dots,\tilde{r}_{n-1}$, where $\tilde{r}_j$ is a lift of $r_j$ whose boundary starts at $\tilde{p}$.
\par
If we put $E:=\begin{pmatrix}0&1\\0&0\end{pmatrix}$, $H:=\begin{pmatrix}1&0\\0&-1\end{pmatrix}$ and $F:=\begin{pmatrix}0&0\\1&0\end{pmatrix}$, then $\{E,H,F\}$ forms a basis of $\sl_2(\C)$.
\par
Therefore $C_0\left(\tilde{\mathcal{K}}\right)\otimes_{\Z[\Pi]}\sl_2(\C)$ is generated by $\tilde{p}\otimes{E}$, $\tilde{p}\otimes{H}$, $\tilde{p}\otimes{F}$, $C_1\left(\tilde{\mathcal{K}}\right)\otimes_{\Z[\Pi]}\sl_2(\C)$ is generated by $\tilde{x}_1\otimes{E}$, $\tilde{x}_1\otimes{H}$, $\tilde{x}_1\otimes{F}$, $\tilde{x}_2\otimes{E}$, $\tilde{x}_2\otimes{H}$, $\tilde{x}_2\otimes{F}$, $\dots$, $\tilde{x}_n\otimes{E}$, $\tilde{x}_n\otimes{H}$, $\tilde{x}_n\otimes{F}$, and  $C_2\left(\tilde{\mathcal{K}}\right)\otimes_{\Z[\Pi]}\sl_2(\C)$ is generated by $\tilde{r}_1\otimes{E}$, $\tilde{r}_1\otimes{H}$, $\tilde{r}_1\otimes{F}$, $\tilde{r}_2\otimes{E}$, $\tilde{r}_2\otimes{H}$, $\tilde{r}_2\otimes{F}$, $\dots$, $\tilde{r}_{n-1}\otimes{E}$, $\tilde{r}_{n-1}\otimes{H}$, $\tilde{r}_{n-1}\otimes{F}$.
So we can choose the ordered bases
\begin{equation*}
  \left\{
    \tilde{r}_1\otimes{E}, \tilde{r}_1\otimes{H}, \tilde{r}_1\otimes{F},
    \tilde{r}_2\otimes{E}, \tilde{r}_2\otimes{H}, \tilde{r}_2\otimes{F},
    \dots,
    \tilde{r}_{n-1}\otimes{E}, \tilde{r}_{n-1}\otimes{H}, \tilde{r}_{n-1}\otimes{F}
  \right\}
\end{equation*}
for $C_2\left(\tilde{\mathcal{K}}\right)\otimes_{\Z[\Pi]}\sl_2(\C)$,
\begin{equation*}
  \left\{
    \tilde{x}_1\otimes{E}, \tilde{x}_1\otimes{H}, \tilde{x}_1\otimes{F},
    \tilde{x}_2\otimes{E}, \tilde{x}_2\otimes{H}, \tilde{x}_2\otimes{F},
    \dots,
    \tilde{x}_n\otimes{E}, \tilde{x}_n\otimes{H}, \tilde{x}_n\otimes{F}
  \right\}
\end{equation*}
for $C_1\left(\tilde{\mathcal{K}}\right)\otimes_{\Z[\Pi]}\sl_2(\C)$, and
\begin{equation*}
  \left\{
    \tilde{p}\otimes{E}, \tilde{p}\otimes{H}, \tilde{p}\otimes{F}
  \right\}
\end{equation*}
for $C_0\left(\tilde{\mathcal{K}}\right)\otimes_{\Z[\Pi]}\sl_2(\C)$,
\par
With respect to these bases the differentials $\partial_2$ and $\partial_1$ are given by the Fox free differential calculus.
Let $\frac{\partial\,r_i}{\partial\,x_j}$ be the Fox derivative \cite{Fox:free_differential_calculus_I}, which is defined as follows:
\begin{itemize}
\item
for words $u$ and $v$ in the $x_j$, $\frac{\partial\,(uv)}{\partial\,x_j}=\frac{\partial\,u}{\partial\,x_j}+u\frac{\partial\,v}{\partial\,x_j}$,
\item
for the empty word $1$, $\frac{\partial\,1}{\partial\,x_j}=0$,
\item
$\frac{\partial\,x_i}{\partial\,x_j}=\delta^i_j$, where $\delta^i_j$ is Kronecker's delta.
\end{itemize}
Note that since $0=\frac{\partial\,(x\cdot x^{-1})}{\partial\,x}=1+x\cdot\frac{\partial\,x^{-1}}{\partial\,x}$, we have $\frac{\partial\,x^{-1}}{\partial\,x}=-x^{-1}$.
\par
The differential $\partial_2$ is given by the following $3(n-1)\times3n$ matrix:
\begin{equation*}
  \partial_2
  =
  \begin{pmatrix}
    \Ad{\rho}\left(\frac{\partial\,r_1}{\partial\,x_1}\right)&\cdots
    &\Ad{\rho}\left(\frac{\partial\,r_{n-1}}{\partial\,x_1}\right)
    \\
    \vdots&\ddots&\vdots
    \\
    \Ad{\rho}\left(\frac{\partial\,r_1}{\partial\,x_n}\right)&\cdots
    &\Ad{\rho}\left(\frac{\partial\,r_{n-1}}{\partial\,x_n}\right)
  \end{pmatrix},
\end{equation*}
where we abuse the notation and write $\Ad{\rho\left(\frac{\partial\,r_j}{\partial\,x_i}\right)}$ as $\Ad{\rho}\left(\frac{\partial\,r_j}{\partial\,x_i}\right)$, and $\Ad{\rho(x_i)}$ is given by the $3\times3$ matrix
\begin{equation*}
  \begin{pmatrix}
    \Ad{\rho(x_i)}(E)_{E}&\Ad{\rho(x_i)}(H)_{E}&\Ad{\rho(x_i)}(F)_{E} \\
    \Ad{\rho(x_i)}(E)_{H}&\Ad{\rho(x_i)}(H)_{H}&\Ad{\rho(x_i)}(F)_{H} \\
    \Ad{\rho(x_i)}(E)_{F}&\Ad{\rho(x_i)}(H)_{F}&\Ad{\rho(x_i)}(F)_{F}
  \end{pmatrix}.
\end{equation*}
Here we use the subscripts so that
\begin{align*}
  \Ad{\rho(x_i)}(E)
  &=
  \Ad{\rho(x_i)}(E)_{E}E+\Ad{\rho(x_i)}(E)_{H}H+\Ad{\rho(x_i)}(E)_{F}F,
  \\
  \Ad{\rho(x_i)}(H)
  &=
  \Ad{\rho(x_i)}(H)_{E}E+\Ad{\rho(x_i)}(H)_{H}H+\Ad{\rho(x_i)}(H)_{F}F,
  \\
  \Ad{\rho(x_i)}(F)
  &=
  \Ad{\rho(x_i)}(F)_{E}E+\Ad{\rho(x_i)}(F)_{H}H+\Ad{\rho(x_i)}(F)_{F}F.
\end{align*}
The differential $\partial_1$ is given by the following $3\times3n$ matrix:
\begin{equation*}
  \begin{pmatrix}
    \Ad{\rho}\left(x_1-1\right)&\cdots&\Ad{\rho}\left(x_n-1\right)
  \end{pmatrix}.
\end{equation*}
Note that $\partial_1\circ\partial_2$, which is a $3\times3(n-1)$ matrix, vanishes, since its $(1,k)$-entry as a block matrix with $3\times3$ blocks is
\begin{equation*}
  \Ad{\rho}(x_1-1)\Ad{\rho}\left(\frac{\partial\,r_k}{\partial\,x_1}\right)
  +\dots+
  \Ad{\rho}(x_n-1)\Ad{\rho}\left(\frac{\partial\,r_k}{\partial\,x_n}\right),
\end{equation*}
which vanishes from the fundamental formula of the free differential calculus \cite[(2.3)]{Fox:free_differential_calculus_I}.
\par
Now we need to fix a basis $\mathbf{h}_{\ast}$ for $H_{\ast}(S^3\setminus{K};\rho)$.
To do this we need several definitions.
\begin{defn}
An irreducible $\SL(2;\C)$ representation $\rho$ is called regular if $\dim H_1(S^3\setminus{K};\rho)=1$.
\end{defn}
\begin{rem}
If a representation $\rho$ is irreducible, then $H_0(S^3\setminus{K};\rho)=0$.
See the calculation in Subsection~\ref{subsec:Reidemeister_fig8_scratch} for a proof.
\par
Therefore if an irreducible representation $\rho$ is regular, then $\dim H_{2}(S^3\setminus{K};\rho)=\dim H_{1}(S^3\setminus{K};\rho)=1$, and that $\dim H_{i}(S^3\setminus{K};\rho)=0$ for $i\ne1,2$ since the Euler characteristic of $S^3\setminus{K}$ is zero.
\end{rem}
So for a regular, irreducible representation we need to choose one homology class for each of $H_{1}$ and $H_{2}$.
To choose a homology class for $H_{1}$ we pick up a simple closed curve on $\partial\left(S^3\setminus{K}\right)$.
The following definition is due to Porti \cite[D{\'e}finition~3.21]{Porti:MAMCAU1997} (see also \cite{Dubois/Huynh/Yamaguchi:JKNOT2009})
\begin{defn}
Let $\gamma$ be a simple closed curve on $\partial(S^3\setminus{K})$.
An irreducible representation is called $\gamma$-regular if
\begin{itemize}
\item
The inclusion $\gamma\hookrightarrow S^3\setminus{K}$ induces a surjective map $H_{1}(\gamma;\rho)\twoheadrightarrow H_{1}(S^3\setminus{K};\rho)$, and
\item
if $\tr(\rho(\pi_1(\partial(S^3\setminus{K}))))\subset\{2,-2\}$, then $\rho(\gamma)\ne\pm\begin{pmatrix}1&0\\0&1\end{pmatrix}$.
\end{itemize}
\end{defn}
For a regular representation $\rho$ we fix a basis $P^{\rho}$ of $H_{0}(\partial(S^3\setminus{K});\rho)$ such that $P^{\rho}$ is invariant under the adjoint action of $\rho(x)$ for any $x\in\pi_1(\partial(S^3\setminus{K}))$ (see the calculation in Subsection~\ref{subsec:Reidemeister_fig8_scratch}).
\par
Then we choose $i_{\ast}\left([\lambda]\otimes P^{\rho}\right)$ as a basis of $H_{1}(S^3\setminus{K};\rho)$, called the reference generator, and $i_{\ast}\left([\partial(S^3\setminus{K})]\otimes P^{\rho}\right)$ as a basis of $H_{1}(S^3\setminus{K};\rho)$, also called the reference generator, where $[\lambda]\in H_1\bigl(\partial(S^3\setminus{K})\bigr)$ is the homology class of the curve $\gamma$, $[\partial(S^3\setminus{K})]\in H_2\bigl(\partial(S^3\setminus{K})\bigr)$ is the fundamental class, and $i\colon\partial(S^3\setminus{K})\to S^3\setminus{K}$ is the inclusion map.
\par
The twisted Reidemeister torsion $\mathbb{T}^{K}_{\gamma}(\rho)$ of $\rho$ associated with $\gamma$ is defined as $\operatorname{Tor}(C_{\ast},\mathbf{b}_{\ast},\mathbf{h}_{\ast})$ defined in \eqref{eq:Reidemeister_def} for these bases.
\subsection{Calculation of the Reidemeister torsion of the figure-eight knot from scratch}
\label{subsec:Reidemeister_fig8_scratch}
Let $\FigEight$ be the figure-eight knot.
In this subsection I calculate $\mathbb{T}^{\FigEight}_{\mu}(\rho_{u,\pm})$.
\par
From Subsection~\ref{subsec:rep_fig8}, we have $\pi_1(S^3\setminus{\FigEight})=\langle x,y\mid xy^{-1}x^{-1}yx=yxy^{-1}x^{-1}y\rangle$.
In this case there is only one relation $r:=xy^{-1}x^{-1}yxy^{-1}xyx^{-1}y^{-1}$ and we have
\begin{align*}
  \frac{\partial\,r}{\partial\,x}
  &=
  1-xy^{-1}x^{-1}+xy^{-1}x^{-1}y+xy^{-1}x^{-1}yxy^{-1}-xy^{-1}x^{-1}yxy^{-1}xyx^{-1},
  \\
  \frac{\partial\,r}{\partial\,y}
  &=
  -xy^{-1}+xy^{-1}x^{-1}-xy^{-1}x^{-1}yxy^{-1}+xy^{-1}x^{-1}yxy^{-1}x
  \\
  &\quad
  -xy^{-1}x^{-1}yxy^{-1}xyx^{-1}y^{-1}.
\end{align*}
\par
Let $\rho_{u,\pm}$ be the representation defined in Subsection~\ref{subsec:CS_fig8}.
We first calculate the homology groups $H_{\ast}(S^3\setminus{\FigEight};\rho_{u,\pm})$.
\par
The adjoint actions of $\rho_{u,\pm}(x)$ are given as follows:
\begin{align*}
  \Ad{\rho_{u,\pm}(x)}(E)
  &:=
  \rho^{-1}_{u,\pm(x)}\cdot E\cdot\rho_{u,\pm}
  =
  \begin{pmatrix}0&e^{-u}\\0&0\end{pmatrix},
  \\
  \Ad{\rho_{u,\pm}(x)}(H)
  &=
  \begin{pmatrix}1&2e^{-u/2}\\0&-1\end{pmatrix},
  \\
  \Ad{\rho_{u,\pm}(x)}(F)
  &=
  \begin{pmatrix}-e^{u/2}&-1\\e^{u}&e^{u/2}\end{pmatrix}.
\end{align*}
So with respect to the basis $\{E,H,F\}$, $\Ad{\rho_{u,\pm}(x)}$ is given by the $3\times3$ matrix
\begin{equation*}
  X:=
  \begin{pmatrix}
    e^{-u}&2e^{-u/2}&-1       \\
    0     &1        &-e^{u/2} \\
    0     &0        &e^{u}
  \end{pmatrix}.
\end{equation*}
Similarly, $\Ad{\rho_{u,\pm}(y)}$ is given by
\begin{equation*}
  Y:=
  \begin{pmatrix}
    e^{-u}          &0              &0\\
    -e^{-u/2}d_{\pm}&1              &0\\
    -d^{2}_{\pm}    &2e^{u/2}d_{\pm}&e^{u}
  \end{pmatrix}
\end{equation*}
with respect to the same basis.
\par
Now the differential $\partial_2$ is given by the $6\times3$ matrix
\begin{equation*}
  \partial_2
  =
  \begin{pmatrix}
    \Ad{\rho_{u,\pm}}\left(\frac{\partial\,r}{\partial\,x}\right) \\[2mm]
    \Ad{\rho_{u,\pm}}\left(\frac{\partial\,r}{\partial\,y}\right) 
  \end{pmatrix},
\end{equation*}
where
\begin{equation*}
\begin{split}
  &\Ad{\rho_{u,\pm}}\left(\frac{\partial\,r}{\partial\,x}\right)
  \\
  =&
  I_3-X^{-1}Y^{-1}X+YX^{-1}Y^{-1}X+Y^{-1}XYX^{-1}Y^{-1}X
  \\
  &
  -X^{-1}YXY^{-1}XYX^{-1}Y^{-1}X,
\end{split}
\end{equation*}
with $I_3$ the $3\times3$ identity matrix, and
\begin{equation*}
\begin{split}
  &\Ad{\rho_{u,\pm}}\left(\frac{\partial\,r}{\partial\,y}\right)
  \\
  =&
  -Y^{-1}X+X^{-1}Y^{-1}X-Y^{-1}XYX^{-1}Y^{-1}X+XY^{-1}XYX^{-1}Y^{-1}X
  \\
  &-Y^{-1}X^{-1}YXY^{-1}XYX^{-1}Y^{-1}X.
\end{split}
\end{equation*}
Note that we need to reverse the order of the multiplication.
\par
By Mathematica we calculate
\begin{equation*}
  \partial_2
  =
  \begin{pmatrix}
    D_{11}&D_{12}&D_{13} \\
    D_{21}&D_{22}&D_{23} \\
    D_{31}&D_{32}&D_{33} \\
    D_{41}&D_{42}&D_{43} \\
    D_{51}&D_{52}&D_{53} \\
    D_{61}&D_{62}&D_{63}
  \end{pmatrix},
\end{equation*}
where
\begin{align*}
  D_{11}
  &:=
  e^{-2u}
  \left(
    (2\cosh(u)-3)
    (e^{3u}+2e^{2u}+1)
    -
    e^{-u}d_{\pm}
    (3e^{3u}-2e^{2u}+4e^{u}-1)
  \right),
  \\
  D_{12}
  &:=
  2e^{-5u/2}
  \\&\quad
  \times
  \left(
  e^{5u}-2e^{4u}-2e^{3u}+4e^{2u}-4e^{u}+1
  -
  d_{\pm}
  (e^{4u}+e^{3u}-2e^{2u}+3e^{u}-1)
  \right),
  \\
  D_{13}
  &:=
  -e^{-u}
  (2\cosh(u)-1)
  \left(
    e^{3u}-e^{2u}-2e^{u}+1-d_{\pm}(e^{2u}+e^{u}-1)
  \right),
  \\
  D_{21}
  &:=
  e^{-5u/2}d_{\pm}
  \left(
    e^{u}(-3e^{u}+2)
    +
    d_{\pm}
    (e^{3u}-e^{2u}+2e^{u}-1)
  \right),
  \\
  D_{22}
  &:=
  -e^{-3u}
  \left(
    -5e^{3u}+14e^{2u}-10e^{u}+2
    +
    2d_{\pm}
    (e^{4u}-e^{3u}+4e^{2u}-4e^{u}+1)
  \right),
  \\
  D_{23}
  &:=
  -e^{-5u/2}
  (e^{u}-1)
  \left(
    (e^{u}-1)
    (e^{2u}+2e^{u}-1)
    -
    d_{\pm}
    (2e^{u}-1)
  \right),
  \\
  D_{31}
  &:=
  e^{-u}
  d_{\pm}^2
  \left(
    -e^{-u}(e^{u}+1)^2(2\cosh(u)-3)
    +
    2d_{\pm}
  \right),
  \\
  D_{32}
  &:=
  2e^{-3u/2}d_{\pm}
  (2\cosh(u)-3)
  \left(
    e^{2u}+e^{u}-1-d_{\pm}
  \right),
  \\
  D_{33}
  &:=
  e^{-2u}(e^{u}-1)(2\cosh(u)-3)
  \left(
    e^{2u}+e^{u}-1-d_{\pm}
  \right),
  \\
  D_{41}
  &:=
  2(2\cosh(u)-3)
  \left(
    -1
    +
    d_{\pm}\cosh(u)
  \right),
  \\
  D_{42}
  &:=
  -
  2e^{-u/2}
  (2\cosh(u)-3)
  (e^{2u}+e^{u}-1-d_{\pm}),
  \\
  D_{43}
  &:=
  (e^{u}-1)
  \left(
    (e^{u}+1)(2\cosh(u)-3)
    +
    e^{-u}d_{\pm}
    (-e^{2u}-e^{u}+1)
  \right),
  \\
  D_{51}
  &:=
  e^{-3u/2}d_{\pm}
  (2\cosh(u)-3)
  (e^{2u}-e^{u}-1+e^{2u}d_{\pm}),
  \\
  D_{52}
  &:=
  -
  \left(
    8\cosh^2(u)-16\cosh(u)+7
    +2d_{\pm}(4\cosh^2(u)-6\cosh(u)+1)
  \right),
  \\
  D_{53}
  &:=
  -e^{-3u/2}(e^{u}-1)
  \left(
    (e^{u}-1)^2
    +
    d_{\pm}(e^{2u}-e^{u}+1)
  \right),
  \\
  D_{61}
  &:=
  -(e^{u}-1)(2\cosh(u)-3)
  \left(
    (2\cosh(u)-3)
    +
    e^{-2u}d_{\pm}
    (e^{3u}-2e^{2u}-1)
  \right),
  \\
  D_{62}
  &:=
  2e^{-u/2}d_{\pm}(e^{u}-1)
  \left(
    2\cosh(u)-2
    +
    d_{\pm}(2\cosh(u)-1)
  \right),
  \\
  D_{63}
  &:=
  (2\cosh(u)-1)
  \left(
    2\cosh(u)-3
    +
    2d_{\pm}(\cosh(u)-1)
  \right).
\end{align*}
(I hope that I could copy them well from the output of Mathematica.)
So we have
\begin{align*}
  \Ker\partial_2
  &=
  \left\langle
    \begin{pmatrix}
      e^{u/2}(d_{\pm}-e^u+1)\\
      -d_{\pm} \\
      d_{\pm}(e^{u/2}-e^{-u/2})
    \end{pmatrix}
  \right\rangle,
  \\
  \Im\partial_2
  &=
  \left\langle
    \partial_2\begin{pmatrix}1\\0\\0\end{pmatrix},
    \partial_2\begin{pmatrix}0\\1\\0\end{pmatrix}
  \right\rangle.
\end{align*}
The differential $\partial_1$ is given by the $3\times6$ matrix
\begin{equation*}
\begin{split}
  \partial_1
  &=
  \begin{pmatrix}
    \Ad{\rho_{u,\pm}}(x-1)&\Ad{\rho_{r,\pm}}(y-1)
  \end{pmatrix}
  \\
  &=
  \begin{pmatrix}
    (X-I_3)&(Y-I_3)
  \end{pmatrix}
  \\
  &=
  \begin{pmatrix}
    e^{-u}-1&2e^{-u/2}&-1      &e^{-u}-1        &0              &0\\
    0       &0        &-e^{u/2}&-e^{-u/2}d_{\pm}&0              &0\\
    0       &0        &e^{u}-1 &-d^{2}_{\pm}    &2e^{u/2}d_{\pm}&e^{u}-1
  \end{pmatrix}.
\end{split}
\end{equation*}
So we have
\begin{align*}
  \Ker\partial_1
  &=
  \Im\partial_2
  +
  \left\langle
    \begin{pmatrix}
      2e^{u/2}\\e^{u}-1\\0\\0\\0\\0
    \end{pmatrix}
  \right\rangle,
  \\
  \Im\partial_1
  &=
  \left\langle
    \partial_1\begin{pmatrix}0\\1\\0\\0\\0\\0\end{pmatrix},
    \partial_1\begin{pmatrix}0\\0\\1\\0\\0\\0\end{pmatrix},
    \partial_1\begin{pmatrix}0\\0\\0\\0\\1\\0\end{pmatrix}
  \right\rangle.
\end{align*}
So $H_2(S^3\setminus{\FigEight};\rho_{u,\pm})=H_1(S^3\setminus{\FigEight};\rho_{u,\pm})\cong\C$.
We can also see that $H_2(S^3\setminus{\FigEight};\rho_{u,\pm})$ is generated by (the homology class of)
$\begin{pmatrix}
    e^{u/2}(d_{\pm}-e^u+1)\\
    -d_{\pm} \\
    d_{\pm}(e^{u/2}-e^{-u/2})
  \end{pmatrix}$
and (the homology class of) $H_1(S^3\setminus{\FigEight};\rho_{u,\pm})$ is generated by $\begin{pmatrix}2e^{u/2}\\e^{u}-1\\0\\0\\0\\0\end{pmatrix}$.
\par
Next we calculate $H_{\ast}(\partial(S^3\setminus{\FigEight});\rho_{u})$.
If we choose $x$ as the meridian, then longitude $\lambda$ is given by $xy^{-1}xyx^{-2}yxy^{-1}x^{-1}$.
So $\Ad{\rho_{u,\pm}(\lambda)}$ acts on $\sl_2(\C)$ as the matrix
\begin{multline*}
  \Lambda
  \\
  {\tiny
  :=
  \hspace{-1mm}
  \begin{pmatrix}
    \ell^{\mp2}&\pm4\ell^{\mp}\cosh(u/2)\sqrt{(2\cosh{u}+1)(2\cosh{u}-3)}
               &-4\cosh^2(u/2)(2\cosh{u}+1)(2\cosh{u}-3)\\
    0          &1&\mp2\ell^{\pm}\cosh(u/2)\sqrt{(2\cosh{u}+1)(2\cosh{u}-3)}   \\
    0          &0                                  &\ell^{\pm2}
  \end{pmatrix}.}
\end{multline*}
Let us consider the twisted chain complex $C'_{\ast}:=C_{\ast}(\partial(S^3\setminus{\FigEight}))\otimes_{\Z\left[\pi_1(\partial(S^3\setminus{\FigEight}))\right]}\sl_2(\C)$.
Since
\begin{equation*}
  \pi_1(\partial(S^3\setminus{\FigEight}))
  =
  \langle
    x,\lambda\mid x\lambda=\lambda x
  \rangle,
\end{equation*}
we can assume that $\dim C'_{2}=3$ (generated by $\{\tilde{r}'\otimes{E},\tilde{r}'\otimes{H},\tilde{r}'\otimes{F}\}$, where $\tilde{r}'$ is a lift of the relation $r':=x\lambda x^{-1}\lambda^{-1}$), $\dim C'_{1}=6$ (generated by $\{\tilde{x}\otimes{E},\tilde{x}\otimes{H},\tilde{x}\otimes{F},\tilde{\lambda}\otimes{E},\tilde{\lambda}\otimes{H},\tilde{\lambda}\otimes{F}\}$, where $\tilde{x}$ and $\tilde{\lambda}$ are lifts of $x$ and $\lambda$, respectively) and $\dim C'_{0}=3$ (generated by $\{\tilde{p}'\otimes{E},\tilde{p}'\otimes{H},\tilde{p}'\otimes{F}\}$, where $\tilde{p}'$ is a lift of a point $p'$ where $\tilde{r}'$, $\tilde{x}$ and $\tilde{\lambda}$ are attached), and that the differentials are given by
\begin{equation*}
\begin{split}
  \partial'_2
  &=
  \begin{pmatrix}
    \Ad{\rho_{u,\pm}}
    \left(\frac{\partial\,(\lambda x\lambda^{-1}x^{-1})}{\partial\,x}\right)
    \\
    \Ad{\rho_{u,\pm}}
    \left(\frac{\partial\,(\lambda x\lambda^{-1} x^{-1})}{\partial\,\lambda}\right)
  \end{pmatrix}
  \\
  &=
  \begin{pmatrix}
    \Ad{\rho_{u,\pm}}(\lambda-\lambda x\lambda^{-1}x^{-1})
    \\
    \Ad{\rho_{u,\pm}}(1-\lambda x\lambda^{-1})
  \end{pmatrix}
  \\
  &=
  \begin{pmatrix}
    \Lambda-I_3 \\
    I_3-X
  \end{pmatrix}
\end{split}
\end{equation*}
and
\begin{equation*}
\begin{split}
  \partial'_1
  &=
  \begin{pmatrix}
    \Ad{\rho_{u,\pm}}(x-1)&\Ad{\rho_{u,\pm}}(\lambda-1)
  \end{pmatrix}
  \\
  &=
  \begin{pmatrix}
    X-I_3&\Lambda-I_3
  \end{pmatrix}.
\end{split}
\end{equation*}
Therefore the kernel of $\partial'_2$ is generated by a non-zero vector $P^{\rho_{u,\pm}}$ such that $\Lambda P^{\rho_{u,\pm}}=P^{\rho_{u,\pm}}$ and $X P^{\rho_{u,\pm}}=P^{\rho_{u,\pm}}$.
We can put
\begin{equation*}
  P^{\rho_{u,\pm}}
  :=
  \begin{pmatrix}
    2e^{u/2}\\
    e^{u}-1 \\
    0
  \end{pmatrix}.
\end{equation*}
We also have
\begin{equation*}
  \Im\partial'_2
  =
  \left\langle
    \partial'_2\begin{pmatrix}1\\0\\0\end{pmatrix},
    \partial'_2\begin{pmatrix}0\\0\\1\end{pmatrix}
  \right\rangle
\end{equation*}
and
\begin{align*}
  \Ker\partial'_1
  &=
  \Im\partial'_2+
  \left\langle
  \begin{pmatrix}
    0\\0\\0\\2e^{u/2}\\e^{u}-1\\0
  \end{pmatrix},
  \begin{pmatrix}
    2e^{u/2}\\e^{u}-1\\0\\0\\0\\0
  \end{pmatrix}
  \right\rangle,
  \\
  \Im\partial'_1
  &=
  \left\langle
    \begin{pmatrix}
      e^{u/2}\\0\\0
    \end{pmatrix},
    \begin{pmatrix}
      -1\\e^{-u/2}\\e^{-u}-1
    \end{pmatrix}
  \right\rangle.
\end{align*}
Therefore we see that $H_{*}(\partial(S^3\setminus{\FigEight});\rho_{u})\cong H_{\ast}(\partial(S^3\setminus{\FigEight});\C)$, and we can choose
\begin{equation*}
  \tilde{h}_1
  :=
  i_{\ast}([\mu]\otimes P^{\rho_{u,\pm}})
  =
  \begin{pmatrix}
    2e^{u/2}\\e^{u}-1\\0\\0\\0\\0
  \end{pmatrix}
\end{equation*}
as the reference generator of $H_1(S^3\setminus{\FigEight};\rho_{u,\pm})$, where $i\colon\partial(S^3\setminus{\FigEight})\to S^3\setminus{\FigEight}$ is the inclusion map.
The reference generator of $H_2(S^3\setminus{\FigEight};\rho_{u,\pm})$ is given by $i_{\ast}([\partial(S^3\setminus{\FigEight})]\otimes P^{\rho_{u,\pm}})$, where $[\partial(S^3\setminus{\FigEight})]\in H_2(\partial(S^3\setminus{\FigEight});\Z)$ is the fundamental class.
Since the fundamental class $[\partial(S^3\setminus{\FigEight})]$ is represented by the $2$-cell whose boundary is attached to
\begin{equation*}
\begin{split}
  \lambda x\lambda^{-1}x^{-1}
  &=
  xy^{-1}xyx^{-2}yxy^{-1}x^{-1}\cdot
  x\cdot
  xyx^{-1}y^{-1}x^{2}y^{-1}x^{-1}yx^{-1}\cdot
  x^{-1}
  \\
  &=
  \operatorname{ad}_{xy^{-1}xyx^{-1}}
  \left(
    \operatorname{ad}_{yx^{-1}}(r){r}^{-1}
  \right),
\end{split}
\end{equation*}
we have
\begin{equation*}
\begin{split}
  \tilde{h}_2
  &:=
  i_{\ast}([H_2(S^3\setminus{\FigEight})]\otimes P)
  \\
  &=
  \Ad{\rho_{u,\pm}(xy^{-1}xyx^{-1})}(\Ad{\rho_{u,\pm}(yx^{-1})}r-r)
  \\
  &=
  (X^{-1}Y-I_3)X^{-1}YXY^{-1}X\begin{pmatrix}2e^{u/2}\\e^{u}-1\\0\end{pmatrix}
  \\
  &=
  \begin{pmatrix}
    2e^{u/2}(d_{\pm}-e^{u}+1)\\
    -2d_{\pm}\\
    2e^{-u/2}(e^{u}-1)d_{\pm}
  \end{pmatrix}.
\end{split}
\end{equation*}
Here we put $\operatorname{ad}_{z}(w):=zwz^{-1}$ for $z,w\in\pi_1(S^3\setminus{\FigEight})$.
\par
Now the twisted Reidemeister torsion $\mathbb{T}^{\FigEight}_{\mu}(\rho_{\mu,\pm})$ is given as
\begin{equation}\label{eq:fig8_Reidemeister}
\begin{split}
  &\frac{
  \left.
  \left[
    \left\{
      \partial_2\begin{pmatrix}1\\0\\0\end{pmatrix},
      \partial_2\begin{pmatrix}0\\1\\0\end{pmatrix}
    \right\}
    \cup
    \left\{\tilde{h}_1\right\}
    \cup
    \left\{
      \begin{pmatrix}0\\1\\0\\0\\0\\0\end{pmatrix},
      \begin{pmatrix}0\\0\\1\\0\\0\\0\end{pmatrix},
      \begin{pmatrix}0\\0\\0\\0\\1\\0\end{pmatrix}
    \right\}
    \right|
    I_6
  \right]}
  {
  \left[
    \left.
    \left\{
      \partial_1\begin{pmatrix}0\\1\\0\\0\\0\\0\end{pmatrix},
      \partial_1\begin{pmatrix}0\\0\\1\\0\\0\\0\end{pmatrix},
      \partial_1\begin{pmatrix}0\\0\\0\\0\\1\\0\end{pmatrix}
    \right\}
    \right|
    I_3
  \right]
  \times
  \left[
    \left.
    \left\{\tilde{h}_2\right\}
    \cup
    \left\{
      \begin{pmatrix}1\\0\\0\end{pmatrix},
      \begin{pmatrix}0\\1\\0\end{pmatrix}
    \right\}
    \right|
    I_3
  \right]}
  \\
  &=
  \mp\frac{\sqrt{(2\cosh{u}+1)(2\cosh{u}-3)}}{2}.
\end{split}
\end{equation}
\subsection{How to calculate the Reidemeister torsion from the twisted Alexander polynomial}
\label{subsec:Reidemeister_Alexander}
In general, as you see in the previous subsection, the calculation of the twisted Reidemeister is very hard.
In this subsection I will describe an easier way.
\par
Here I will explain how to calculate the twisted Reidemeister torsion associated with the {\em longitude} $\lambda$ when our representation is $\lambda$-regular.
\par
Let $\langle x_1,\dots,x_n\mid r_1,\dots,r_{n-1}\rangle$ be a Wirtinger presentation of $\pi_1\left(S^3\setminus{K}\right)$.
For a representation $\rho\colon\pi_1\left(S^3\setminus{K}\right)\to\SL(2;\C)$, put $\Phi:=\Ad{\rho}\otimes\alpha$, where $\alpha\colon\pi_1(S^3\setminus{K})\to\Z\cong H_1\left(S^3\setminus{K}\right)$ is the Abelianization sending $x_i\mapsto{t}$ for any $i$, where $t$ is a generator of $\Z$ and we denote $\Phi(x_i)$ by $t\Ad{\rho(x_i)}$, noting that $x_i$ is sent to the generator of $\Z$ by $\alpha$.
\par
Now we follow the technique used in \cite{Dubois/Huynh/Yamaguchi:JKNOT2009}.
We use the following theorem.
\begin{thm}[{\cite[Theorem~3.1.2]{Yamaguchi:ANNIF2008}}]
Suppose that a representation $\rho$ is $\lambda$-regular for the longitude $\lambda$.
Then the twisted Reidemeister torsion $\mathbb{T}^{K}_{\lambda}(\rho)$ of $\rho$ associated with $\lambda$ is given as
\begin{equation*}
  \mathbb{T}^{K}_{\lambda}(\rho)
  =
  -\lim_{t\to1}\frac{\mathcal{T}^{K}(\Phi;t)}{t-1},
\end{equation*}
where $\mathcal{T}^{K}(\Phi;t)$ is the twisted Reidemeister torsion of $\Phi$ with parameter $t$.
\end{thm}
\begin{rem}
Note that from \cite[Theorem~A]{Kitano:PACJM1996} $\mathcal{T}^{K}(\rho)$ coincides with the twisted Alexander polynomial of $K$ associated with $\Phi$.
\end{rem}
\begin{rem}
From \cite[Proposition~3.1.1]{Yamaguchi:ANNIF2008}, if $\rho$ is $\lambda$-regular, then the twisted chain complex associated with $\Phi$ is acyclic.
Therefore its Reidemeister torsion is well-defined without worrying about basis for the homology group.
\end{rem}
For actual computation we use the Fox free differential calculus again \cite{Fox:free_differential_calculus_I}.
Consider the  $3(n-1)\times3(n-1)$ matrix $\Phi\left(\frac{\partial r_i}{\partial x_j}\right)$ ($i\in\{1,\dots,n-1\}$,$j\in\{1,\dots,n\}\setminus\{{l}\}$ for some $l$ with entry $\C[t,t^{-1}]$, where $t$ is a generator of $\Z$ and $\frac{\partial r_i}{\partial x_j}$ is described in \ref{subsec:Reidemeister_definition}.
Combining a result by T.~Kitano and Y.~Yamaguchi, we have the following theorem.
\begin{thm}[\cite{Kitano:PACJM1996,Kirk/Livingston:TOPOL1999}]\label{thm:Reidemeister_Alexander}
Let $K$ be a knot and $\rho\colon\pi_1(S^3\setminus{K})\to\SL(2;\C)$ a $\lambda$-regular representation, where $\lambda$ is the preferred longitude.
Then we have
\begin{equation*}
  \mathbb{T}_{\lambda}^{K}(\rho)
  =
  \pm
  \lim_{t\to1}
  \frac{\det\Phi\left(\frac{\partial r_i}{\partial x_j}\right)}{(t-1)\det\Phi(x_{{l}}-1)}.
\end{equation*}
\end{thm}
\begin{rem}
We can determine the sign in the formula above.
See \cite{Dubois/Huynh/Yamaguchi:JKNOT2009} for details.
\end{rem}
To obtain the twisted Reidemeister torsion associated with the meridian $\mu$, we need the following result of Porti \cite[Th\'eor\`eme~4.1]{Porti:MAMCAU1997}.
\begin{thm}[{\cite[Th\'eor\`eme~4.1]{Porti:MAMCAU1997}}]\label{thm:Porti}
Let $\rho$ be a representation that sends the meridian $\mu$ to $\begin{pmatrix}e^{u/2}&\ast\\0&e^{-u/2}\end{pmatrix}$ and the longitude $\lambda$ to $\begin{pmatrix}e^{v(u)/2}&\ast\\0&e^{-v(u)/2}\end{pmatrix}$ with $u$ a complex parameter.
Then we have
\begin{equation*}
  \mathbb{T}_{\mu}^{K}(\rho)
  =
  \pm
  \frac{\mathbb{T}_{\lambda}^{K}(\rho)}{d\,v(u)/d\,u}.
\end{equation*}
\end{thm}
\subsection{Twisted Reidemeister torsion of the figure-eight knot again}
Here we calculate the twisted Reidemeister torsion of the figure-eight knot again by using Theorem~\ref{thm:Reidemeister_Alexander}
\par
Put $\Phi_{u,\pm}:=\Ad{\rho_{u,\pm}}\otimes\alpha$.
Then we have $\Phi_{u,\pm}(x)=tX$ and $\Phi_{u,\pm}(y)=tY$ since $\alpha(x)=\alpha(y)=t$, where $X$ and $Y$ are given in Subsection~\ref{subsec:Reidemeister_fig8_scratch}.
So we have
\begin{equation*}
\begin{split}
  &\Phi_{u,\pm}\left(\frac{\partial\,r}{\partial\,x}\right)
  \\
  =&
  I_3-t^{-1}X^{-1}Y^{-1}X+YX^{-1}Y^{-1}X+Y^{-1}XYX^{-1}Y^{-1}X
  \\
  &
  -tX^{-1}YXY^{-1}XYX^{-1}Y^{-1}X
\end{split}
\end{equation*}
and Mathematica tells us
\begin{multline*}
  \det\Phi_{u,\pm}\left(\frac{\partial\,r}{\partial\,x}\right)
  \\
  =
  -t^{-3}e^{-2u}(t-1)^2(t-e^u)(te^u-1)\bigl(e^u+t(-2+(t-1)e^u-2e^{2u})\bigr).
\end{multline*}
Since $\det\Phi_{u,\pm}(y-1)=(te^{-u}-1)(t-1)(te^{u}-1)$, we have
\begin{equation*}
  \mathbb{T}^{\FigEight}_{\lambda}\left(\rho_{u,\pm}\right)
  =
  \lim_{t\to1}\frac{\det\Phi_{u,\pm}\left(\frac{\partial\,r}{\partial\,x}\right)}{(t-1)\det\Phi_{u,\pm}(x-1)}
  =
  4\cosh{u}-1
\end{equation*}
from Theorem~\ref{thm:Reidemeister_Alexander}.
\par
Now we apply Theorem~\ref{thm:Porti}.
From \eqref{eq:longitude_fig8}, we have
\begin{equation*}
\begin{split}
  \frac{d\,v(u)}{d\,u}
  &=
  \pm2
  \frac{d}{d\,u}
  \log\left(\cosh(2u)-\cosh{u}-1-\sinh{u}\sqrt{(2\cosh{u}+1)(2\cosh{u}-3)}\right)
  \\
  &=
  \pm\frac{2(1-4\cosh{u})}{\sqrt{(2\cosh{u}+1)(2\cosh{u}-3)}}.
\end{split}
\end{equation*}
Therefore we finally have
\begin{equation*}
  \mathbb{T}^{\FigEight}_{\mu}\left(\rho_{u,\pm}\right)
  =
  \frac{\mathbb{T}^{\FigEight}_{\lambda}\left(\rho_{u,\pm}\right)}{d\,v(u)/d\,u}
  =
  \frac{\sqrt{(2\cosh{u}+1)(2\cosh{u}-3)}}{2}
\end{equation*}
up to a sign, which coincides with our previous calculation \eqref{eq:fig8_Reidemeister}.
\subsection{Twisted Reidemeister torsion of a torus knot}\label{subsec:Reidemeister_torus_knot}
Now we calculate the twisted Reidemeister torsion of a torus knot by using the twisted Alexander polynomial described in Subsection~\ref{subsec:Reidemeister_Alexander}.
Here we assume that our representation is both $\lambda$-regular and $\mu$-regular.
It is known that any irreducible representation of $\pi_1\left(S^3\setminus{T(a,b)}\right)$ is $\lambda$-regular and $\mu$-regular (see \cite[Example~1]{Dubois:CANMB2006}).
So the representation $\rho_{u,\omega_1}$ given in Subsection~\ref{subsec:rep_torus_knot} is $\lambda$-regular and $\mu$-regular unless $u=\exp\left(\frac{(2k+1)\pi\sqrt{-1}}{2a+1}\right)$.
\par
Putting $r:=(xy)^ax(xy)^{-a}y^{-1}$, we have
\begin{equation*}
  \frac{\partial\,r}{\partial\,x}
  =
  \sum_{i=0}^{a-1}(xy)^{i}
  +
  (xy)^{a}
  \left(
    1-x(xy)^{-a}\left(\sum_{i=0}^{a-1}(xy)^{i}\right)
  \right).
\end{equation*}
For $z\in\pi_1\left(S^3\setminus{T(a,b)}\right)$, put $\Phi_{u,\omega_1}(z):=\alpha(z)\Ad{\rho_{u,\omega_1}(z)}$.
We also put $X:=\Phi_{u,\omega_1}(x)$ and $Y:=\Phi_{u,\omega_1}(y)$.
Then we have
\begin{align*}
  X
  &=
  t
  \begin{pmatrix}
    e^{-u}&2e^{-u/2}&-1 \\
    0     &1        &-e^{u/2}\\
    0     &0        &e^{u}
  \end{pmatrix},
  \\
  Y
  &=
  t
  \begin{pmatrix}
    e^{-u}&0&0 \\
    e^{-u/2}\left(\omega_1+\omega_1^{-1}-2\cosh{u}\right)&1&0\\
    -\left(\omega_1+\omega_1^{-1}-2\cosh{u}\right)^2&
    -2e^{u/2}\left(\omega_1+\omega_1^{-1}-2\cosh{u}\right)&e^{u}
  \end{pmatrix},
  \intertext{and}
  \Phi_{u,\omega_1}\left(\frac{\partial\,r}{\partial\,x}\right)
  &=
  \sum_{i=0}^{a-1}(YX)^{i}
  +
  \left(
    I_3-\left(\sum_{i=0}^{a-1}(YX)^{i}\right)(YX)^{-a}X
  \right)
  (YX)^{a}.
\end{align*}
By Mathematica we have
\begin{equation*}
  \det\Phi_{u,\omega_1}\left(\frac{\partial\,r}{\partial\,x}\right)
  =
  \frac{(t^{2a+1}-1)^2(t^{2a+1}+1)(te^{u}-1)(te^{-u}-1)}{(t+1)(t^2-\omega_1^2)(t^2-\omega_1^{-2})}.
\end{equation*}
Since
\begin{equation*}
  \det\Phi_{u,\omega_1}(y-1)
  =
  (t-1)(te^u-1)(te^{-u}-1),
\end{equation*}
we have
\begin{equation*}
  \mathbb{T}^{T(2,2a+1)}_{\lambda}\left(\rho_{u,\omega_1}\right)
  =
  \pm
  \lim_{t\to1}
  \frac{\det\Phi_{u,\omega_1}\left(\frac{\partial\,r}{\partial\,x}\right)}{(t-1)\det\Phi(y-1)}
  =
  \pm
  \left(
    \frac{2a+1}{\omega_1-\omega_1^{-1}}
  \right)^2
\end{equation*}
Since $d\,v(u)/d\,u=-2(2a+1)$ from \eqref{eq:rep_longitude_torus_knot}, we have
\begin{equation}\label{eq:torus_knot_Reidemeister}
  \mathbb{T}^{T(a,b)}_{\mu}\left(\rho_{u,\omega_1}\right)
  =
  \frac{\mathbb{T}^{T(a,b)}_{\lambda}\left(\rho_{u,\pm}\right)}{d\,v(u)/d\,u}
  =
  \frac{2a+1}{2\left(\omega_1-\omega_1^{-1}\right)^2}
\end{equation}
up to a sign.
\subsection{Twisted Reidemeister torsion of a twice-iterated torus knot}
In this subsection I calculate the twisted Reidemeister torsion assuming that representations are both $\mu$-regular and $\lambda$-regular.
Unfortunately I do not know whether this is true or not.
So a safer way is to say that I will calculate the twisted Alexander polynomial.
\par
Throughout this subsection we put
\begin{align*}
  r_1&:=(xy)^{a}x(xy)^{-a}y^{-1},
  \\
  r_2&:=pqx^{-1}
  \\
  \intertext{and}
  r_3&:=\lambda_{C}x^bp\lambda^{-1}_{C}x^{-b}q^{-1}
\end{align*}
with $\lambda_{C}:=y(xy)^{2a}x^{-4a-1}$ so that $\pi_1(E)=\langle x,y,p,q\mid r_1,r_2,r_3\rangle$ (see Subsection~\ref{subsec:rep_iterated_torus_knot}).
\par
From Theorem~\ref{thm:Reidemeister_Alexander}, $\mathbb{T}^{T(2,2a+1)^{(2,2b+1)}}_{\lambda}(\rho)$ is given by the following:
\begin{equation*}
  \pm
  \lim_{t\to1}
  \frac
  {
  \det
  \begin{pmatrix}
    \Phi\left(\frac{\partial\,r_1}{\partial\,x}\right)&
    \Phi\left(\frac{\partial\,r_2}{\partial\,x}\right)&
    \Phi\left(\frac{\partial\,r_3}{\partial\,x}\right)
    \\
    \Phi\left(\frac{\partial\,r_1}{\partial\,p}\right)&
    \Phi\left(\frac{\partial\,r_2}{\partial\,p}\right)&
    \Phi\left(\frac{\partial\,r_3}{\partial\,p}\right)
    \\
    \Phi\left(\frac{\partial\,r_1}{\partial\,q}\right)&
    \Phi\left(\frac{\partial\,r_2}{\partial\,q}\right)&
    \Phi\left(\frac{\partial\,r_3}{\partial\,q}\right)
  \end{pmatrix}
  }
  {(t-1)\det\Phi(y-1)},
\end{equation*}
where $\Phi:=\Ad{\rho}\otimes\alpha$ with $\rho=\rho^{\rm{AN}}_{u,\omega_2}$, $\rho^{\rm{NA}}_{u,\omega_1}$ or $\rho^{\rm{NN}}_{u,\omega_1,\omega_3}$ given in Subsection~\ref{subsec:rep_iterated_torus_knot}.
Since $r_1$ does not contain $p$ or $q$, we have $\frac{\partial\,r_1}{\partial\,p}=\frac{\partial\,r_1}{\partial\,q}=0$ and so we have
\begin{equation}\label{eq:Reidemeister_det}
  \mathbb{T}^{T(2,2a+1)^{(2,2b+1)}}_{\lambda}(\rho)
  =
  \pm
  \lim_{t\to1}
  \frac{
  \det\Phi\left(\frac{\partial\,r_1}{\partial\,x}\right)
  \det
  \begin{pmatrix}
    \Phi\left(\frac{\partial\,r_2}{\partial\,p}\right)&
    \Phi\left(\frac{\partial\,r_3}{\partial\,p}\right)
    \\
    \Phi\left(\frac{\partial\,r_2}{\partial\,q}\right)&
    \Phi\left(\frac{\partial\,r_3}{\partial\,q}\right)
  \end{pmatrix}
  }
  {(t-1)\det\Phi(y-1)}.
\end{equation}
\par
From the Fox free differential calculus, we have
\begin{align}
  \frac{\partial r_1}{\partial x}
  &=
  \sum_{i=1}^{a-1}(xy)^i+(xy)^a\left(1-x(xy)^{-a}\left(\sum_{i=0}^{a-1}(xy)^i\right)\right),
  \label{eq:Fox_iterated_torus_knot_r1x}
  \\
  \frac{\partial\,r_2}{\partial\,p}
  &=
  1,
  \label{eq:Fox_iterated_torus_knot_r2p}
  \\
  \frac{\partial\,r_2}{\partial\,q}
  &=p,
  \label{eq:Fox_iterated_torus_knot_r2q}
  \\
  \frac{\partial\,r_3}{\partial\,p}
  &=
  \lambda_{C}x^b,
  \label{eq:Fox_iterated_torus_knot_r3p}
  \\
  \frac{\partial\,r_3}{\partial\,q}
  &=
  -\lambda_{C}x^bp\lambda^{-1}_{C}x^{-b}q^{-1}
  =-r_3.
  \label{eq:Fox_iterated_torus_knot_r3q}
\end{align}
Note that since $\rho(r_3)=I_3$ and $\alpha(r_3)=0$, $\Phi(r_3)=I_3$ for any choice of $\rho$.
\subsubsection{$\Im\rho_{C}$ is Abelian and $\Im\rho_{P}$ is non-Abelian.}
Let $\rho^{\rm{AN}}_{u,\omega_2}$ be the representation given in \eqref{eq:rep_AN}.
For $z\in\pi_1\left(S^3\setminus{T(2,2a+1)^{(2b+1)}}\right)$, put $\Phi(z):=\alpha(z)\Ad{\rho^{\rm{AN}}_{u,\omega_2}(z)}$.
We also put $P:=\Phi(p)$, $Q:=\Phi(q)$, $X:=\Phi(x)$, and $Y:=\Phi(y)$.
\par
Then we have
\begin{align*}
  P
  &=
  t
  \begin{pmatrix}
    e^{-u}&2e^{-u/2}&-1\\
    0     & 1       &-e^{u/2} \\
    0     & 0       &e^{u}
  \end{pmatrix},
  \\
  Q
  &=
  t
  \begin{pmatrix}
    e^{-u}                                               &0&0 \\
    e^{-u/2}\left(\omega_2+\omega_2^{-1}-2\cosh{u}\right)&1&0 \\
    -\left(\omega_2+\omega_2^{-1}-2\cosh{u}\right)^2&
    -2e^{u/2}\left(\omega_2+\omega_2^{-1}-2\cosh{u}\right)&e^{u}
  \end{pmatrix},
  \\
  X
  &=
  Y
  =
  QP
  =
  t^2
  T_2
  \begin{pmatrix}
    \omega_2^{-2}&2\omega_2^{-1}&-1\\
    0            &1             &-\omega_2\\
    0            &0             &\omega_2^2
  \end{pmatrix}
  T_2^{-1},
\end{align*}
where
\begin{equation*}
  T_2
  :=
  \begin{pmatrix}
    e^{-u/2}                        &0                              &0 \\
    \omega_2^{-1}-e^{-u}            &1                              &0 \\
    -e^{u/2}(\omega_2^{-1}-e^{-u})^2&2e^{-u/2}-2\omega_2^{-1}e^{u/2}&e^{u/2}
  \end{pmatrix}.
\end{equation*}
\par
Now we calculate the determinants in \eqref{eq:Reidemeister_det}.
\par
First we calculate $\det\Phi\left(\frac{\partial\,r_1}{\partial\,x}\right)$.
Since $X$ and $Y$ commute, we have
\begin{equation*}
\begin{split}
  \Phi\left(\frac{\partial r_1}{\partial x}\right)
  &=
  \sum_{i=1}^{a-1}X^{2i}+\left(1-\left(\sum_{i=0}^{a-1}X^{2i}\right)X^{1-2a}\right)X^{2a}
  \\
  &=
  \left(X^{2a+1}+I_3\right)\left(X+I_3\right)^{-1}
\end{split}
\end{equation*}
from \eqref{eq:Fox_iterated_torus_knot_r1x}.
Since the eigenvalues of $X$ are $t^2$, $t^2\omega_2^{2}$, and $t^2\omega_2^{-2}$, we have
\begin{equation*}
  \det\Phi\left(\frac{\partial r_1}{\partial x}\right)
  =
  \frac{\left(t^{2(2a+1)}+1\right)
        \left(t^{2(2a+1)}\omega_2^{2(2a+1)}+1\right)
        \left(t^{2(2a+1)}\omega_2^{-2(2a+1)}+1\right)}
       {\left(t^2+1\right)\left(t^2\omega_2^2+1\right)\left(t^2\omega_2^{-2}+1\right)}
\end{equation*}
and so
\begin{equation*}
  \det\Phi\left(\frac{\partial r_1}{\partial x}\right)\Bigm|_{t=1}
  =
  \left(\frac{\omega_2^{2a+1}+\omega_2^{-(2a+1)}}{\omega_2+\omega_2^{-1}}\right)^2
\end{equation*}
Note that this coincides with $\Delta\left(T(2,2a+1);\omega_2^2\right)^2$.
From \eqref{eq:Fox_iterated_torus_knot_r2p}--\eqref{eq:Fox_iterated_torus_knot_r3q}, we have
\begin{equation*}
\begin{split}
  \det
  \begin{pmatrix}
    \Phi\left(\frac{\partial\,r_2}{\partial\,p}\right)&
    \Phi\left(\frac{\partial\,r_3}{\partial\,p}\right)
    \\
    \Phi\left(\frac{\partial\,r_2}{\partial\,q}\right)&
    \Phi\left(\frac{\partial\,r_3}{\partial\,q}\right)
  \end{pmatrix}
  &=
  \begin{vmatrix}
    I_3&X^b \\
    P &-I_3
  \end{vmatrix}
  \\
  &=
  \begin{vmatrix}
    I_3&X^b \\
    O_3&-I_3-PX^b
  \end{vmatrix}
  \\
  &=
  -\det(I_3+PX^b)
  \\
  &=
  -\left(t^{2b+1}+1\right)\left(t^{2b+1}-1\right)^2
\end{split}
\end{equation*}
since $\Phi(\lambda)=\Phi(r_3)=I_3$.
\par
Since we obtain
\begin{equation*}
  \det\begin{pmatrix}\Phi(y-1)\end{pmatrix}
  =
  \left(t^2-1\right)\left(t^2\omega_2^2-1\right)\left(t^2\omega_2^{-2}-1\right),
\end{equation*}
we finally have
\begin{equation*}
\begin{split}
  \mathbb{T}^{T(2,2a+1)^{(2,2b+1)}}_{\lambda}(\rho^{\rm{AN}}_{u,\omega_2})
  =&
  \pm
  \left(\frac{\omega_2^{2a+1}+\omega_2^{-(2a+1)}}{\omega_2+\omega_2^{-1}}\right)^2
  \lim_{t\to1}
  \frac{-\left(t^{2b+1}+1\right)\left(t^{2b+1}-1\right)^2}
       {\left(t^2-1\right)\left(t^2\omega_2^2-1\right)\left(t^2\omega_2^{-2}-1\right)}
  \\
  =&
  \pm
  \left(
    (2b+1)\times
    \frac{\omega_2^{2a+1}+\omega_2^{-(2a+1)}}{\omega_2^{2}-\omega_2^{-2}}
  \right)^2
\end{split}
\end{equation*}
from \eqref{eq:Reidemeister_det}.
Note that this also follows from \cite[Theorem~3.7]{Kirk/Livingston:TOPOL1999}.
\par
Since $d\,v(u)/d\,u=-2(2b+1)$ from \eqref{eq:AN_longitude}, we have
\begin{equation}\label{eq:AN_Reidemeister}
  \mathbb{T}^{T(2,2a+1)^{(2,2b+1)}}_{\mu}(\rho^{\rm{AN}}_{u,\omega_2})
  =
  \frac{(2b+1)}{2}
  \times
  \left(
    \frac{\omega_2^{2a+1}+\omega_2^{-(2a+1)}}{\omega_2^{2}-\omega_2^{-2}}
  \right)^2
\end{equation}
up to a sign from Theorem~\ref{thm:Porti}.
Note that here we assume that $\rho^{\rm{AN}}_{u,\omega_2}$ is $\mu$-regular and $\lambda$-regular.
\subsubsection{$\Im\rho_{C}$ is non-Abelian and $\Im\rho_{P}$ is Abelian.}
Let $\rho^{\rm{AN}}_{u,\omega_1}$ be the representation given in \eqref{eq:rep_NA}.
For $z\in\pi_1\left(S^3\setminus{T(2,2a+1)^{(2b+1)}}\right)$, put $\Phi(z):=\alpha(z)\Ad{\rho^{\rm{NA}}_{u,\omega_1}(z)}$.
We also put $P:=\Phi(p)$, $Q:=\Phi(q)$, $X:=\Phi(x)$, $Y:=\Phi(y)$, and $\Lambda_C:=\Phi(\lambda_C)$.
\par
We have
\begin{align*}
  X
  &=
  t^2
  \begin{pmatrix}
    e^{-2u}&2e^{-u}&-1\\
    0      & 1     &-e^{u} \\
    0      & 0     &e^{2u}
  \end{pmatrix},
  \\
  Y
  &=
  t^2
  \begin{pmatrix}
    e^{-2u}                                                               &0&0\\
    e^{-u}\left(\omega_1+\omega_1^{-1}-2\cosh{2u}\right)&1&0\\
    -\left(\omega_1+\omega_1^{-1}-2\cosh{2u}\right)^2
    &-2e^{u}\left(\omega_1+\omega_1^{-1}-2\cosh{2u}\right)&e^{2u}
  \end{pmatrix},
  \\
  P
  &=
  Q=
  t
  \begin{pmatrix}
    e^{-u}&\frac{2}{1+e^{u}}&\frac{-1}{(e^{u/2}+e^{-u/2})^2}\\
    0     & 1               &\frac{-e^{u}}{1+e^{u}} \\
    0     & 0                &e^{u}
  \end{pmatrix},
  \\
  \Lambda_{C}
  &=
  \begin{pmatrix}
    e^{4(2a+1)u}&2e^{2(2a+1)u}\frac{\sinh\bigl(2(2a+1)u\bigr)}{\sinh(u)}&-\frac{\sinh^2(2(2a+1)u)}{\sinh^2(u)}
    \\
    0&1&-e^{-2(2a+1)u}\frac{\sinh\bigl(2(2a+1)u\bigr)}{\sinh(u)}
    \\
    0&0&e^{-4(2a+1)u}
  \end{pmatrix}.
\end{align*}
\par
Now we calculate the determinants.
As in the case of the torus knot (Subsection~\ref{subsec:Reidemeister_torus_knot}), we have
\begin{equation*}
  \det\left(\Phi\left(\frac{\partial\,r_1}{\partial\,x}\right)\right)
  =
  \frac{(t^{2(2a+1)}-1)^2(t^{2(2a+1)}+1)(t^2e^{2u}-1)(t^2e^{-2u}-1)}
       {(t^2+1)(t^4-\omega_1^2)(t^4-\omega_1^{-2})}
\end{equation*}
and
\begin{equation*}
  \det\Phi(y-1)
  =
  (t^2-1)(t^2e^{2u}-1)(t^2e^{-2u}-1).
\end{equation*}
\par
We also have
\begin{equation*}
\begin{split}
  \det
  \begin{pmatrix}
    \Phi\left(\frac{\partial\,r_2}{\partial\,p}\right)&
    \Phi\left(\frac{\partial\,r_3}{\partial\,p}\right)
    \\
    \Phi\left(\frac{\partial\,r_2}{\partial\,q}\right)&
    \Phi\left(\frac{\partial\,r_3}{\partial\,q}\right)
  \end{pmatrix}
  &=
  \begin{vmatrix}
    I_3&\Lambda_CX^b \\
    P  &-I_3
  \end{vmatrix}
  \\
  &=
  \begin{vmatrix}
    I_3&\Lambda_CX^b \\
    O_3&-I_3-P\Lambda_CX^b
  \end{vmatrix}
  \\
  &=
  -\det(I_3+P\Lambda_CX^b).
\end{split}
\end{equation*}
Since $P$, $\Lambda_C$, and $X$ are all upper-triangle matrices, we have
\begin{multline*}
 \det(I_3+P\Lambda_CX^b)
 \\
 =
 (1+t^{2b+1})(1+t^{2b+1}e^{(2b+1-4(2a+1))u})(1+t^{2b+1}e^{-(2b+1-4(2a+1))u}).
\end{multline*}
\par
So we finally have
\begin{equation*}
  \mathbb{T}^{T(2,2a+1)^{(2,2b+1)}}_{\lambda}(\rho^{\rm{NA}}_{u,\omega_1})
  =
  \pm
  \left(
    (2a+1)
    \times
    \frac{\cosh\left(\frac{(2b+1-4(2a+1))u}{2}\right)}{\omega_1-\omega_1^{-1}}
  \right).
\end{equation*}
Since $d\,v(u)/d\,u=-8(2a+1)$ from \eqref{eq:NA_longitude}, we have
\begin{equation}\label{eq:NA_Reidemeister}
  \mathbb{T}^{T(2,2a+1)^{(2,2b+1)}}_{\mu}(\rho^{\rm{NA}}_{u,\omega_1})
  =
  \frac{(2a+1)}{8}
  \times
  \left(
    \frac{\cosh\left(\frac{(2b+1-4(2a+1))u}{2}\right)}{\omega_1-\omega_1^{-1}}
  \right)^2
\end{equation}
up to a sign, assuming that $\rho^{\rm{NA}}_{u,\omega_1}$ is $\mu$-regular and $\lambda$-regular.
\subsubsection{Both $\Im\rho_{C}$ and $\Im\rho_{P}$ are non-Abelian.}
Let $\rho^{\rm{NN}}_{u,\omega_1,\omega_3}$ be the representation given in \eqref{eq:rep_NN}.
For $z\in\pi_1\left(S^3\setminus{T(2,2a+1)^{(2b+1)}}\right)$, put $\Phi(z):=\alpha(z)\Ad{\rho^{\rm{NN}}_{u,\omega_1,\omega_3}(z)}$.
We also put $P:=\Phi(p)$, $Q:=\Phi(q)$, $X:=\Phi(x)$, $Y:=\Phi(y)$, and $\Lambda_C:=\Phi(\lambda_C)$.
\par
We have
\begin{align*}
  &X
  \\
  =&
  t^2
  \begin{pmatrix}
    e^{-2u}&2e^{-3u/2}&-e^{-u}
    \\
    -e^{-5u/2}(e^u-g)(e^u-g^{-1})&-2e^{-2u}-1+2e^{-u}(g+g^{-1})&-e^{-u/2}(g+g^{-1}-e^{-u})
    \\
    -e^{-u}(g+g^{-1}-e^u-e^{-u})^2&2e^{-5u/2}(e^u-g)(e^u-g^{-1})(e^u(g+g^{-1})-1)&(g+g^{-1}-e^{-u})^2
  \end{pmatrix}
  \\
  =&
  t^2
  T_3
  \begin{pmatrix}
    \omega_3^{-2}&2\omega_3^{-1}&-1 \\
    0            &1             &-\omega_3 \\
    0            &0             &\omega_3^{2}
  \end{pmatrix}
  T_3^{-1},
  \\
  &Y
  \\
  =&
  {\footnotesize
  t^2
  \begin{pmatrix}
    \omega_3^{-2}&0&0
    \\
    e^{-u/2}\omega_3^{-1}\left(\omega_3-\omega_3^{-1}-e^{u}\left(\omega_3^2+1-\omega_1-\omega_1^{-1}\right)\right)
    &1&0
    \\
    -e^{-u}\left(\omega_3^{-1}-\omega_3+e^{u}\left(\omega_3^2+1-\omega_1-\omega_1^{-1}\right)\right)^2
    &2e^{-u/2}\left(1-\omega_3^2+e^{u}\omega_3\left(\omega_3^2+1-\omega_1-\omega_1^{-1}\right)\right)&\omega_3^{2}
  \end{pmatrix}},
  \\
  &P
  \\
  =&
  t
  \begin{pmatrix}
    e^{-u}&2e^{-u/2}&-1 \\
     0 & 1          &-e^{u/2}\\
     0 & 0          &e^{u}
  \end{pmatrix},
  \\
  &Q
  \\
  =&
  5
  \begin{pmatrix}
    e^{-u}&0&0 \\
    e^{u/2}(\omega_3+\omega_3^{-1}-2\cosh{u})&1&0 \\
    -(\omega_3+\omega_3^{-1}-2\cosh{u})^2&-2e^{u/2}(\omega_3+\omega_3-2\cosh{u})&e^{u}
  \end{pmatrix},
  \\
  &\Lambda_C
  \\
  =&
  T_3
  \begin{pmatrix}
    \omega_3^{4(2a+1)}&-\frac{2(\omega_3^{4(2a+1)}-1)}{\omega_3-\omega_3^{-1}}
    &-\frac{\omega_3^{-4(2a+1)}(\omega_3^{4(2a+1)}-1)^2}{(\omega_3-\omega_3^{-1})^2}
    \\
    0&1&-\frac{\omega_3^{-4(2a+1)}-1}{\omega_3-\omega_3^{-1}}
    \\
    0&0&\omega_3^{-4(2a+1)}
  \end{pmatrix}
  T_3^{-1}
\end{align*}
where
\begin{equation*}
  T_3:=
  \begin{pmatrix}
    e^{-u/2}&0&0 \\
    \omega_3^{-1}-e^{-u}&1&0 \\
    -e^{u/2}(e^{-u}-\omega_3^{-1})^2&2\left(e^{-u/2}-\omega_3^{-1}e^{u/2}\right)&e^{u/2}
  \end{pmatrix}.
\end{equation*}
As in the previous case, Mathematica tells us that
\begin{equation*}
\begin{split}
  \det
  \begin{pmatrix}
    \Phi\left(\frac{\partial\,r_2}{\partial\,p}\right)&
    \Phi\left(\frac{\partial\,r_3}{\partial\,p}\right)
    \\
    \Phi\left(\frac{\partial\,r_2}{\partial\,q}\right)&
    \Phi\left(\frac{\partial\,r_3}{\partial\,q}\right)
  \end{pmatrix}
  &=
  -\det\left(I_3+P\Lambda_CX^b\right)
  =
  -(t^{2b+1}-1)^2(t^{2b+1}+1),
  \\
  \det\Phi\left(\frac{\partial\,r_1}{\partial\,x}\right)
  &=
  \frac{(t^{2(2a+1)}-1)^2(t^{2(2a+1)}+1)(t^2\omega_3^{2}-1)(t^2\omega_3^{-2}-1)}
       {(t^2+1)(t^4-\omega_1^2)(t^4-\omega_1^{-2})},
  \\
  \intertext{and}
  \det\Phi(y-1)
  &=
  (t^2-1)(t^2\omega_3^2-1)(t^2\omega_3^{-2}-1).
\end{split}
\end{equation*}
So we see that the twisted Alexander polynomial $\mathcal{T}^{T(2,2a+1)^{(2,2b+1)}}(\Phi^{\rm{NN}}_{u,\omega_1,\omega_3})$ is given as follows.
\begin{equation*}
\begin{split}
  &\mathcal{T}^{T(2,2a+1)^{(2,2b+1)}}(\Phi^{\rm{NN}}_{u,\omega_1,\omega_3})
  \\
  =&
  -
  \frac{(t^{2b+1}-1)^2(t^{2b+1}+1)(t^{2(2a+1)}-1)^2(t^{2(2a+1)}+1)}
       {(t^2+1)(t^4-\omega_1^2)(t^4-\omega_1^{-2})(t^2-1)}.
\end{split}
\end{equation*}
Therefore we conclude that the twisted Alexander polynomial is divisible by $(t-1)$ three times, not once as in the case where the representation is regular.
We also have
\begin{equation*}
  \lim_{t\to\infty}
  \frac{\mathcal{T}^{T(2,2a+1)^{(2,2b+1)}}(\Phi^{\rm{NN}}_{u,\omega_1,\omega_3})}{(t-1)^3}
  =
  \left(
    \frac{2(2a+1)(2b+1)}{(\omega_1-\omega_1^{-1})}
  \right)^2.
\end{equation*}
Since $d\,v(u)/d\,u=-2(2b+1)$, we have
\begin{equation}\label{eq:NN_Reidemeister}
  \lim_{t\to\infty}
  \frac{\mathcal{T}^{T(2,2a+1)^{(2,2b+1)}}(\Phi^{\rm{NN}}_{u,\omega_1,\omega_3})}{(t-1)^3\frac{d\,v(u)}{d\,u}}
  =
  -2(2b+1)
  \left(
    \frac{(2a+1)}{(\omega_1-\omega_1^{-1})}
  \right)^2.
\end{equation}
I do not know whether this is the twisted Reidemeister torsion or not.
\section{The colored Jones polynomial and its asymptotic behavior}
\label{sec:Jones}
For an unoriented link diagram $|D|$, let $\langle{D}\rangle$ be the Kauffman bracket defined by the following two axioms.
\begin{align*}
  \left\langle
    \raisebox{-4mm}{\includegraphics[scale=0.2]{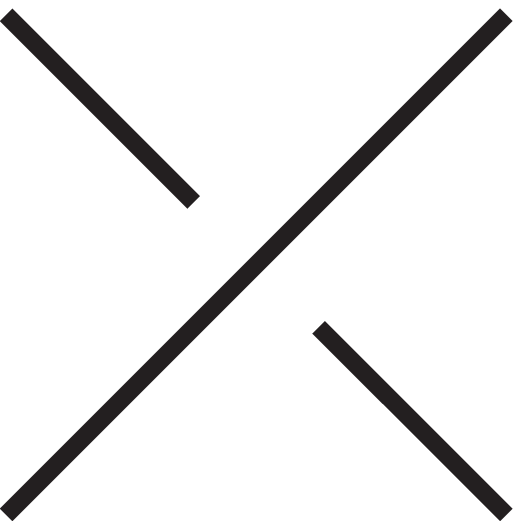}}
  \right\rangle
  &=
  A
  \left\langle
    \raisebox{-4mm}{\includegraphics[scale=0.2]{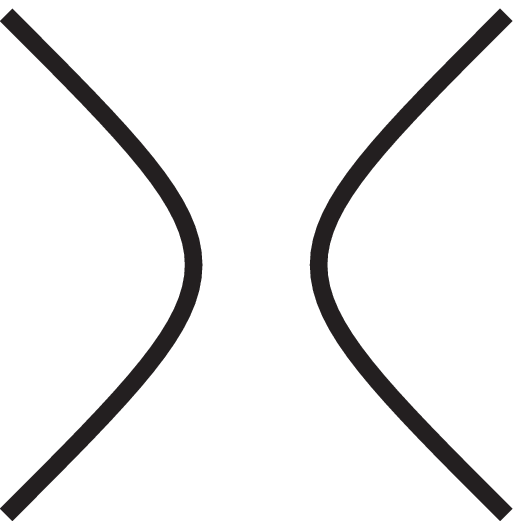}}
  \right\rangle
  +
  A^{-1}
  \left\langle
    \raisebox{-4mm}{\includegraphics[scale=0.2]{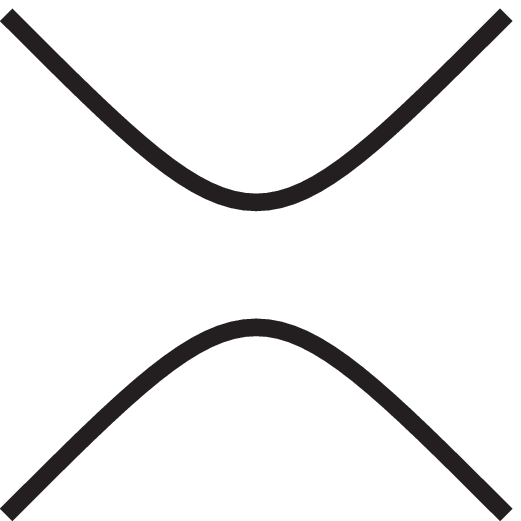}}
  \right\rangle,
  \\
  \left\langle
    U^{c}
  \right\rangle
  &=
  (-A^2-A^{-2})^{c},
\end{align*}
where $U^c$ is the trivial $c$ component link diagram \cite{Kauffman:TOPOL87}.
The Jones polynomial $J_2(K;t)$ of a knot $K$ with a diagram $D$ is defined as
\begin{equation*}
  \frac{\left(-A^3\right)^{-w(D)}\left\langle|D|\right\rangle}{-A^2-A^{-2}}
  \Bigg|_{t:=A^{4}},
\end{equation*}
where $|D|$ is the unoriented diagram obtained from $D$ by forgetting the orientation, and $w(D)$ is the writhe of $D$ that is the sum of the signs of $D$.
Note that $J_2(K;t)$ is a ``quantum'' normalized version of the original Jones polynomial $V_K(t)$ \cite{Jones:BULAM385}.
More precisely it satisfies the following two axioms.
\begin{gather*}
  tJ_2
  \left(
    \raisebox{-4mm}{\includegraphics[scale=0.2]{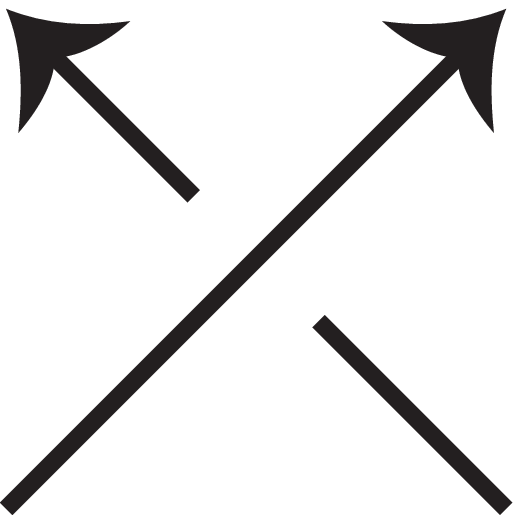}}
  \right)
  -t^{-1}J_2
  \left(
    \raisebox{-4mm}{\includegraphics[scale=0.2]{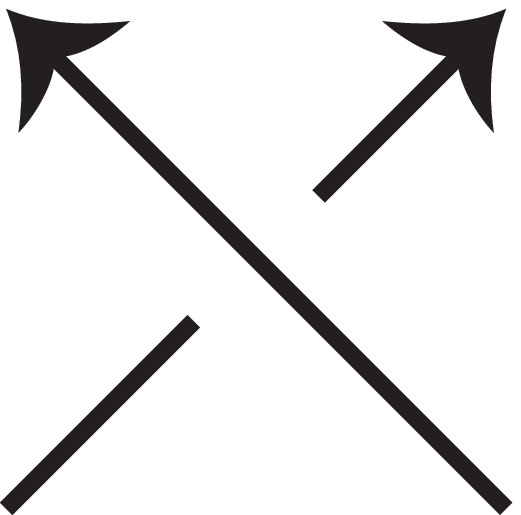}}
  \right)
  =
  \left(t^{1/2}-t^{-1/2}\right)J_2
  \left(
    \raisebox{-4mm}{\includegraphics[scale=0.2]{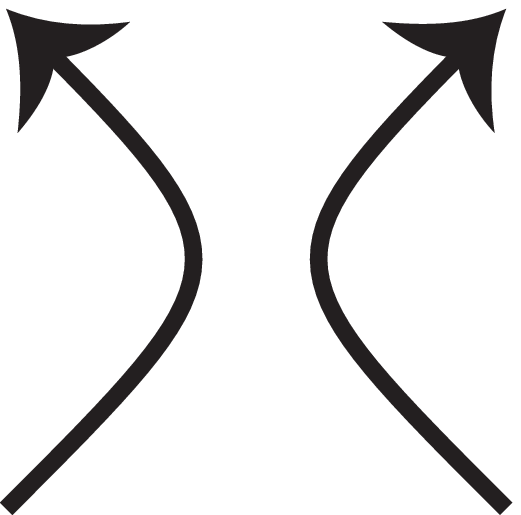}}
  \right),
  \\
  J_2(O)
  =1,
\end{gather*}
where $O$ is the unknot.
\par
The $N$-colored Jones polynomial is defined as
\begin{equation*}
  \frac{\left((-1)^{N-1}A^{N^2-1}\right)^{-w(D)}
    \left\langle
      \raisebox{-9mm}{\includegraphics[scale=0.2]{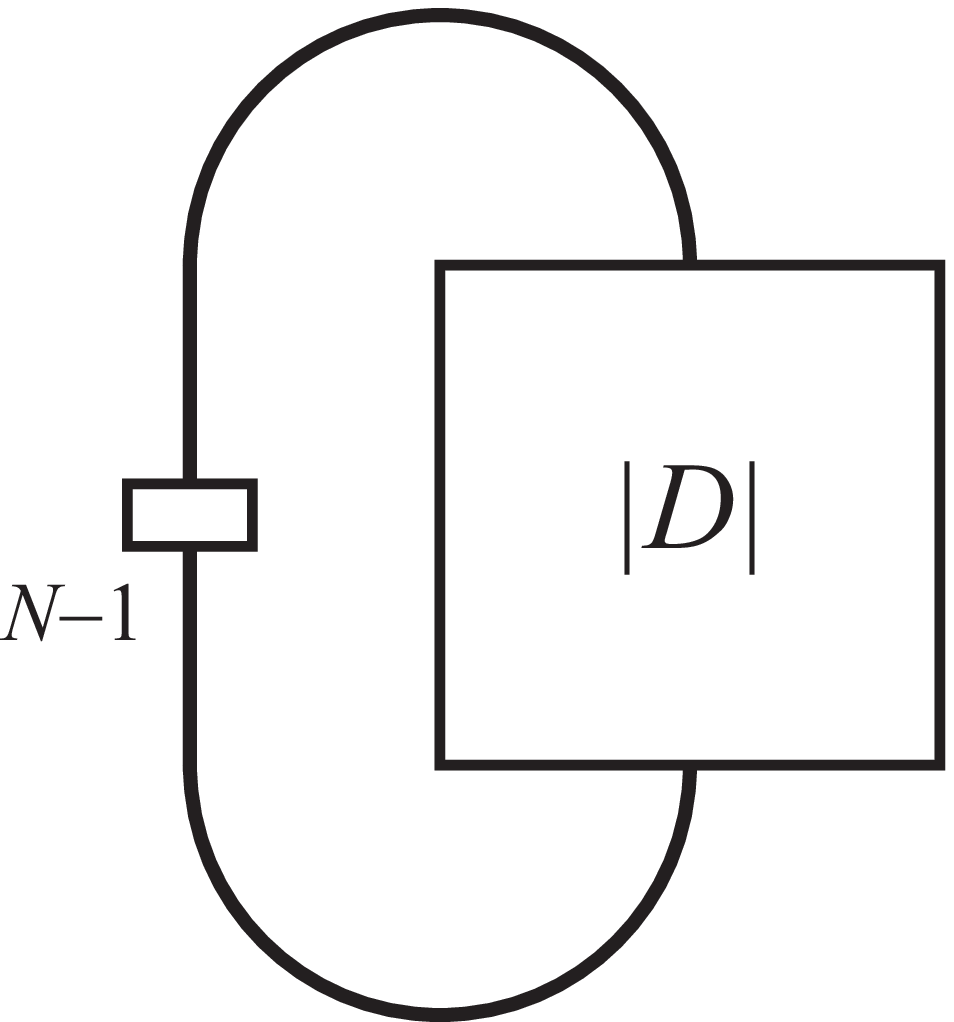}}
    \right\rangle
  }{(-1)^{N-1}\frac{A^{2N}-A^{-2N}}{A^2-A^{-2}}}
  \left.
    \vphantom{\raisebox{-9mm}{\includegraphics[scale=0.2]{idempotent.eps}}}
  \right|_{t:=A^{4}},
\end{equation*}
where a small box with $N-1$ beside it denotes the Jones--Wenzl idempotent \cite{Wenzl:CRMAR87}.
See, for example, \cite{Masbaum/Vogel:PACJM1994} or \cite[Chapter~14]{Lickorish:1997}) for more details.
\subsection{Torus knot}
In this subsection I calculate the colored Jones polynomials of $T(2,2a+1)$ and $T(2,2a+1)^{(2,2b+1)}$ by using linear skein theory.
\subsubsection{The colored Jones polynomial}
Let $T(2,2a+1)$ be the torus knot of type $(2,2a+1)$ as depicted in Figure~\ref{fig:torus_knot}, where $a$ is a positive integer.
\par
We first calculate the Kauffman bracket of the diagram replacing the knot diagram in Figure~\ref{fig:torus_knot} with the Jones--Wenzl idempotent.
By linear skein theory (see \cite[Chapter~14]{Lickorish:1997} for example) we have
\begin{equation}\label{eq:torus_knot_idempotent}
\begin{split}
  &\left\langle
    \quad
    \raisebox{-22mm}{\includegraphics[scale=0.2]{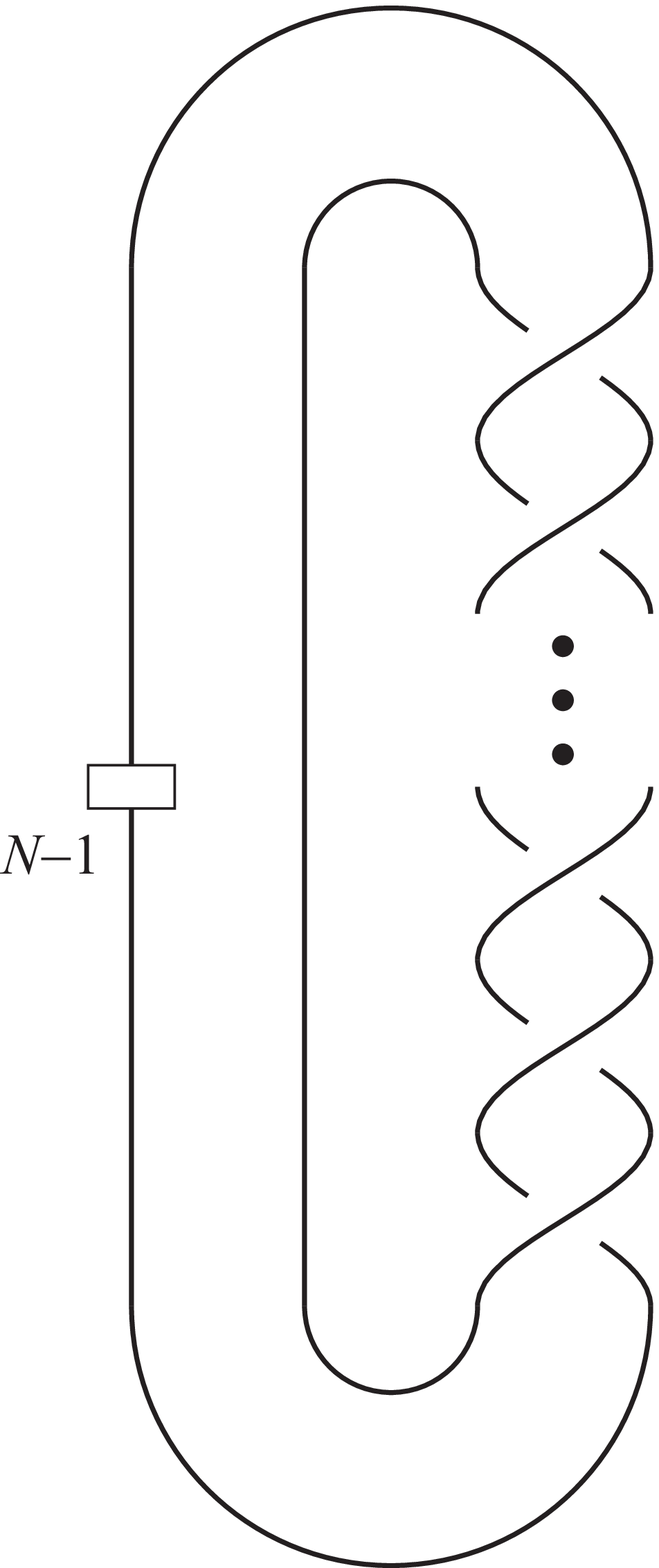}}
    \quad
  \right\rangle
  \\
  =&
  \sum_{c=0}^{N-1}\frac{\Delta_{2c}}{\theta(N-1,N-1,2c)}
  \left\langle
    \quad
    \raisebox{-22mm}{\includegraphics[scale=0.2]{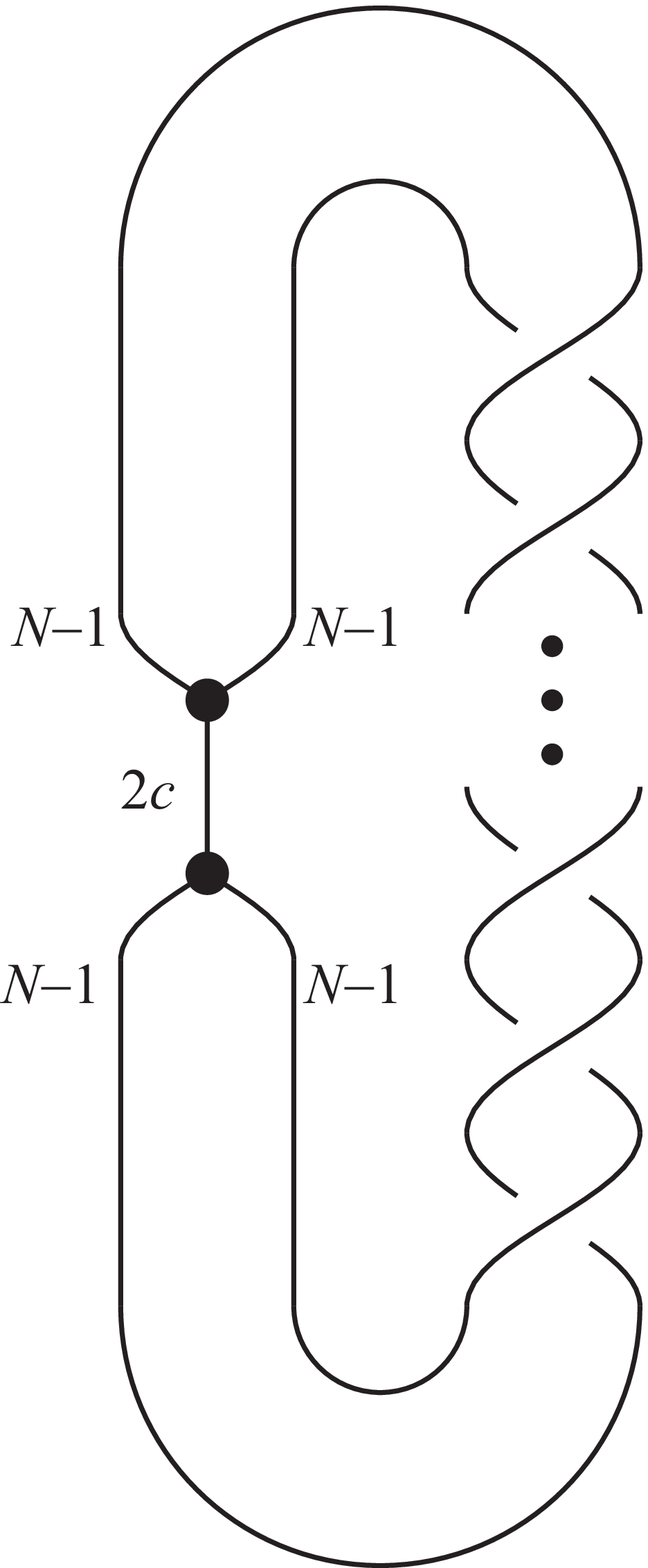}}
    \quad
  \right\rangle
  \\
  =&
  \sum_{c=0}^{N-1}\frac{\Delta_{2c}}{\theta(N-1,N-1,2c)}
  \left((-1)^{c-N+1}A^{-2(N-1)+2c+2c^2-(N-1)^2}\right)^{2a+1}
  \\
  &\quad\quad\times
  \left\langle
    \quad
    \raisebox{-12mm}{\includegraphics[scale=0.2]{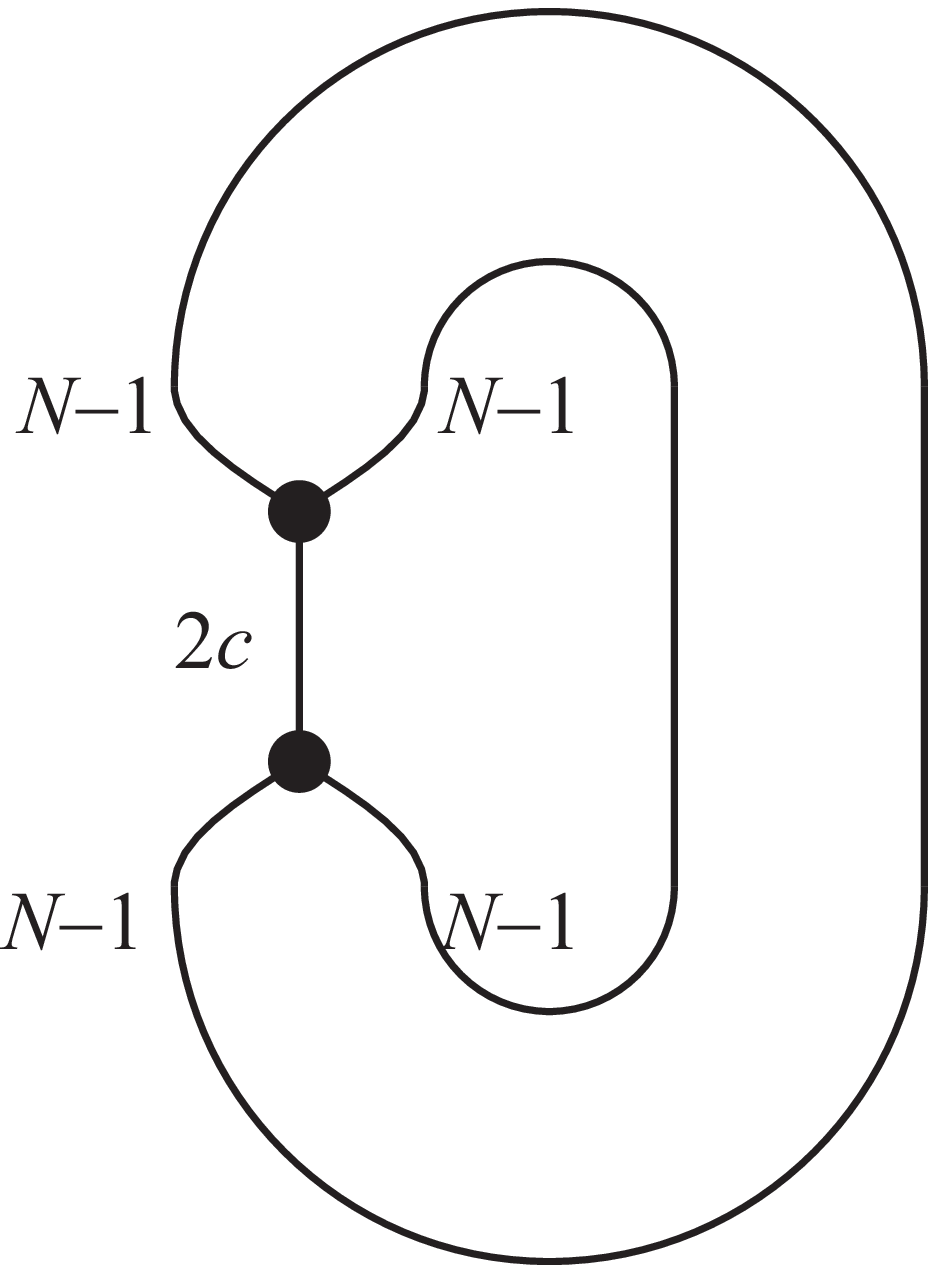}}
    \quad
  \right\rangle
  \\
  =&
  \sum_{c=0}^{N-1}
  (-1)^{c-N+1}A^{(2a+1)(2c^2+2c-N^2+1)}
  \frac{A^{2(2c+1)}-A^{-2(2c+1)}}{A^2-A^{-2}}.
\end{split}
\end{equation}
The colored Jones polynomial is given from this by multiplying $\left((-1)^{N-1}A^{N^2-1}\right)^{-(2a+1)}$, dividing by $(-1)^{N-1}\frac{A^{2N}-A^{-2N}}{A^2-A^{-2}}$ and replacing $A$ with $t^{1/4}$.
Therefore we have
\begin{equation}\label{eq:colored_Jones_torus_knot}
\begin{split}
  &J_{N}(T(2,2a+1);t)
  \\
  =&
  \frac{(-1)^{N-1}t^{-(2a+1)(N^2-1)/4}}{t^{N/2}-t^{-N/2}}
  \\
  &\times
  \sum_{c=0}^{N-1}
  (-1)^{c}t^{(2a+1)(2c^2+2c-N^2+1)/4}
  \left(t^{(2c+1)/2}-t^{-(2c+1)/2}\right)
  \\
  =&
  \frac{(-1)^{N-1}t^{-(2a+1)(N^2-1)/2}}{t^{N/2}-t^{-N/2}}
  \sum_{c=0}^{N-1}
  (-1)^{c}t^{(2a+1)(c^2+c)/2}
  \left(t^{(2c+1)/2}-t^{-(2c+1)/2}\right).
\end{split}
\end{equation}
\par
For a formula for a general torus knot, see \cite{Rosso/Jones:JKNOT93,Morton:MATPC95}.
\subsubsection{Asymptotic behavior of the colored Jones polynomial}
In this sub-subsection we study the asymptotic behavior of the colored Jones polynomial $J_{N}\bigl(T(2,2a+1);\exp(\xi/N)\bigr)$ for large $N$.
We assume that $\xi$ is not purely imaginary and $\Im\xi\ge0$.
To do that we will use a special case of the the saddle point method (see for example \cite[Theorems~7.2.9]{Marsden/Hoffman:Complex_Analysis}).
\begin{thm}\label{thm:saddle}
Let $C_{\theta}=\{t\exp(\theta\sqrt{-1})\}$ be a line in the complex plane passing through the origin.
For a constant $H$ with $\Re\bigl(H^{-1}\exp(2\theta\sqrt{-1})\bigr)>0$, we have
\begin{equation*}
  \int_{\gamma}g(\zeta)\exp\left[\frac{-N\zeta^2}{H}\right]\,d\zeta
  \underset{N\to\infty}{\sim}
  \sqrt{\frac{\pi H}{N}}g(0),
\end{equation*}
where $f(N)\underset{N\to\infty}{\sim}g(N)$ means that $f(N)$ and $g(N)$ are asymptotically equivalent, that is, $f(N)=g(N)\bigl(1+o(1)\bigr)$ for $N\to\infty$.
\end{thm}
Note that the assumption $\Re\bigl(H^{-1}\exp(2\theta\sqrt{-1})\bigr)>0$ is to make the integral converge.
\par
Replacing $t$ with $\exp(\xi/N)$ in \eqref{eq:colored_Jones_torus_knot} we have
\begin{equation*}
\begin{split}
  &J_{N}(T(2,2a+1);\exp(\xi/N))
  \\
  =&
  \frac{(-1)^{N-1}\exp\left[\frac{-(2a+1)(N^2-1)\xi}{2N}\right]}{2\sinh(\xi/2)}
  \\
  &\times
  \left(
    \sum_{c=0}^{N-1}(-1)^{c}\exp\left[\frac{\bigl((2a+1)(c^2+c)+2c+1\bigr)\xi}{2N}\right]
  \right.
  \\
  &\quad\quad-
  \left.
    \sum_{c=0}^{N-1}(-1)^{c}\exp\left[\frac{\bigl((2a+1)(c^2+c)-2c-1\bigr)\xi}{2N}\right]
  \right).
\end{split}
\end{equation*}
Note that since $\xi$ is not purely imaginary, $\sinh(\xi/2)$ does not vanish.
Put
\begin{equation*}
  \Sigma_{\pm}
  :=
  \sum_{c=0}^{N-1}(-1)^{c}\exp\left[\frac{\bigl((2a+1)(c^2+c)\pm(2c+1)\bigr)\xi}{2N}\right]
\end{equation*}
so that
\begin{equation*}
  J_{N}(T(2,2a+1);\exp(\xi/N)
  =
  \frac{(-1)^{N-1}\exp\left[\frac{-(2a+1)(N^2-1)\xi}{2N}\right]}{2\sinh(\xi/2)}
  (\Sigma_{+}-\Sigma_{-}).
\end{equation*}
We have
\begin{equation*}
\begin{split}
  \Sigma_{\pm}
  =&
  \exp\left[\frac{-\xi}{N}\left(\frac{2a+1}{8}+\frac{1}{2(2a+1)}\right)\right]
  \\
  &\times
  \sum_{c=0}^{N-1}
  (-1)^{c}
  \exp\left[\frac{(2a+1)\xi}{2N}\left(c+\frac{2a+1\pm2}{2(2a+1)}\right)^2\right].
\end{split}
\end{equation*}
Now we use the following formula:
\begin{equation*}
  \sqrt{\frac{\alpha}{\pi}}
  \int_{C_{\theta}}\exp(-\alpha x^2+px)\,dx
  =
  \exp\left(\frac{p^2}{4\alpha}\right),
\end{equation*}
where $C_{\theta}$ is the line $\{t\exp(\theta\sqrt{-1})\mid t\in\R\}$.
We choose $\theta$ so that $\Re(\alpha\exp(2\theta\sqrt{-1}))>0$ to make the integral converge.
\par
If we choose $\varphi$ so that $\Re(\xi^{-1}\exp(2\varphi\sqrt{-1}))>0$ we have
\begin{equation*}
\begin{split}
  &\Sigma_{\pm}
  \\
  =&
  \sqrt{\frac{N}{2(2a+1)\xi\pi}}
  \exp\left[\frac{-\xi}{N}\left(\frac{2a+1}{8}+\frac{1}{2(2a+1)}\right)\right]
  \\
  &\times
  \sum_{c=0}^{N-1}
  (-1)^{c}
  \int_{C_{\varphi}}
  \exp\left[\frac{-N}{2(2a+1)\xi}x^2+\left(c+\frac{2a+1\pm2}{2(2a+1)}\right)x\right]\,dx
  \\
  =&
  \sqrt{\frac{N}{2(2a+1)\xi\pi}}
  \exp\left[\frac{-\xi}{N}\left(\frac{2a+1}{8}+\frac{1}{2(2a+1)}\right)\right]
  \\
  &\times
  \int_{C_{\varphi}}
  \exp\left[\frac{-N}{2(2a+1)\xi}x^2+\frac{x}{2}\pm\frac{x}{2a+1}\right]
  \left(
    \sum_{c=0}^{N-1}(-1)^{c}\exp(cx)
  \right)
  \,dx
  \\
  =&
  \sqrt{\frac{N}{2(2a+1)\xi\pi}}
  \exp\left[\frac{-\xi}{N}\left(\frac{2a+1}{8}+\frac{1}{2(2a+1)}\right)\right]
  \\
  &\times
  \int_{C_{\varphi}}
  \exp\left[\frac{-N}{2(2a+1)\xi}x^2\right]
  \exp\left(\frac{\pm x}{2a+1}\right)
  \left(\frac{1-(-1)^N\exp(Nx)}{\exp(x/2)+\exp(-x/2)}\right)
  \,dx.
\end{split}
\end{equation*}
Therefore we have
\begin{equation*}
\begin{split}
  &\Sigma_{+}-\Sigma_{-}
  \\
  =&
  \sqrt{\frac{N}{2(2a+1)\xi\pi}}
  \exp\left[\frac{-\xi}{N}\left(\frac{2a+1}{8}+\frac{1}{2(2a+1)}\right)\right]
  \\
  &\times
  \left(
    \int_{C_{\varphi}}
    \frac{\sinh\left(\frac{x}{2a+1}\right)}{\cosh\left(\frac{x}{2}\right)}
    \exp\left[\frac{-N}{2(2a+1)\xi}x^2\right]
    \,dx
  \right.
  \\
  &\quad
  -(-1)^N
  \left.
    \int_{C_{\varphi}}
    \frac{\sinh\left(\frac{x}{2a+1}\right)}{\cosh\left(\frac{x}{2}\right)}
    \exp\left[\frac{-N}{2(2a+1)\xi}x^2\right]
   \exp(Nx)
    \,dx
  \right).
\end{split}
\end{equation*}
By the saddle point method (Theorem~\ref{thm:saddle}), the first integral is asymptotically equivalent to $0$.
We calculate the second integral.
We have
\begin{equation*}
\begin{split}
  &\int_{C_{\varphi}}
  \frac{\sinh\left(\frac{x}{2a+1}\right)}{\cosh\left(\frac{x}{2}\right)}
  \exp\left[\frac{-N}{2(2a+1)\xi}x^2\right]
  \exp(Nx)
  \,dx
  \\
  =&
  \exp\left[\frac{(2a+1)\xi N}{2}\right]
  \int_{C_{\varphi}}
  \frac{\sinh\left(\frac{x}{2a+1}\right)}{\cosh\left(\frac{x}{2}\right)}
  \exp\left[\frac{-N}{2(2a+1)\xi}\bigl(x-(2a+1)\xi\bigr)^2\right]
  \,dx.
\end{split}
\end{equation*}
If we let ${C}_{\varphi}+(2a+1)\xi$ be the line $\{t\exp(\varphi\sqrt{-1})+(2a+1)\xi\mid t\in\R\}$, then by the residue theorem we have
\begin{equation*}
\begin{split}
  &\int_{C_{\varphi}}
  \frac{\sinh\left(\frac{x}{2a+1}\right)}{\cosh\left(\frac{x}{2}\right)}
  \exp\left[\frac{-N}{2(2a+1)\xi}\bigl(x-(2a+1)\xi\bigr)^2\right]
  \,dx
  \\
  =&
  \int_{C_{\varphi}+(2a+1)\xi}
  \frac{\sinh\left(\frac{x}{2a+1}\right)}{\cosh\left(\frac{x}{2}\right)}
  \exp\left[\frac{-N}{2(2a+1)\xi}\bigl(x-(2a+1)\xi\bigr)^2\right]
  \,dx
  \\
  &+
  2\pi\sqrt{-1}
  \\
  &\quad\times
  \sum_{k}\Res
  \left(
    \frac{\sinh\left(\frac{x}{2a+1}\right)}{\cosh\left(\frac{x}{2}\right)}
    \exp\left[\frac{-N}{2(2a+1)\xi}\bigl(x-(2a+1)\xi\bigr)^2\right];
    x=(2k+1)\pi\sqrt{-1}
  \right)
  \\
  =&
  \int_{C_{\varphi}+(2a+1)\xi}
  \frac{\sinh\left(\frac{x}{2a+1}\right)}{\cosh\left(\frac{x}{2}\right)}
  \exp\left[\frac{-N}{2(2a+1)\xi}\bigl(x-(2a+1)\xi\bigr)^2\right]
  \,dx
  \\
  &+
  2\pi\sqrt{-1}
  \sum_{k}(-1)^{k+1}\sqrt{-1}
  \sinh\left(\frac{(2k+1)\pi\sqrt{-1}}{2a+1}\right)
  \\
  &\quad\times
  \exp\left[\frac{-N}{2(2a+1)\xi}\bigl((2k+1)\pi\sqrt{-1}-(2a+1)\xi\bigr)^2\right],
\end{split}
\end{equation*}
where $\Res(f(x);x=x_0)$ is the residue of $f(x)$ at $x_0$ and $k$ runs over all integers such that $(2k+1)\pi\sqrt{-1}$ is between $C_{\varphi}$ and $C_{\varphi}+(2a+1)\xi$.
If $C_{\varphi}+(2a+1)\xi$ passes through a pole of $\cosh(x/2)$ we avoid it by changing $C_{\varphi}+(2a+1)\xi$ slightly.
Note that we assume that $\xi$ is not on the imaginary axis.
Note also that since $\Im\xi\ge0$, $C_{\varphi}+(2a+1)\xi$ is above $C$.
\par
Putting $y:=x-(2a+1)\xi$ the integral becomes
\begin{equation*}
  \int_{C_{\varphi}}
  \frac{\sinh\left(\frac{y+(2a+1)\xi}{2a+1}\right)}{\cosh\left(\frac{y+(2a+1)\xi}{2}\right)}
  \exp\left[\frac{-N}{2(2a+1)\xi}y^2\right]
  \,dy.
\end{equation*}
By the saddle point method (Theorem~\ref{thm:saddle}) this is asymptotically equivalent to
\begin{equation*}
\begin{split}
  \sqrt{\frac{2(2a+1)\pi\xi}{N}}
  \frac{\sinh\left(\xi\right)}{\cosh\left(\frac{(2a+1)\xi}{2}\right)}.
\end{split}
\end{equation*}
Therefore we finally have the following asymptotic equivalence, which is a special case of Theorem~\ref{thm:Hikami/Murakami_intro}.
\begin{equation*}
\begin{split}
  &J_{N}\bigl(T(2,2a+1);\exp(\xi/N)\bigr)
  \\
  \underset{N\to\infty}{\sim}&
  \frac{(-1)^{N-1}\exp\left[\frac{-(2a+1)(N^2-1)\xi}{2N}\right]}{2\sinh(\xi/2)}
  \sqrt{\frac{N}{2(2a+1)\xi\pi}}
  \exp\left[\frac{-\xi}{N}\left(\frac{2a+1}{8}+\frac{1}{2(2a+1)}\right)\right]
  \\
  &\times
  (-1)^{N+1}
  \exp\left[\frac{(2a+1)\xi N}{2}\right]
  \\
  &\times
  \left(
    \sqrt{\frac{2(2a+1)\pi\xi}{N}}
    \frac{\sinh\left(\xi\right)}{\cosh\left(\frac{(2a+1)\xi}{2}\right)}
  \right.
  \\
  &+
  \left.
    2\pi\sqrt{-1}
    \sum_{k}(-1)^{k+1}\sqrt{-1}
    \sinh\left(\frac{(2k+1)\pi\sqrt{-1}}{2a+1}\right)
    \exp\left[\frac{-N}{2(2a+1)\xi}\bigl((2k+1)\pi\sqrt{-1}-(2a+1)\xi\bigr)^2\right]
  \right)
  \\
  =&
  \frac{\sinh\left(\xi\right)}{2\sinh(\xi/2)\cosh\left(\frac{(2a+1)\xi}{2}\right)}
  \\
  &+
  \frac{1}{2\sinh(\xi/2)}
  \sqrt{\frac{N}{2(2a+1)\xi\pi}}
  \\
  &\times
  2\pi\sqrt{-1}
  \sum_{k}(-1)^{k}
  \sin\left(\frac{(2k+1)\pi}{2a+1}\right)
  \exp\left[\frac{-N}{2(2a+1)\xi}\bigl((2k+1)\pi\sqrt{-1}-(2a+1)\xi\bigr)^2\right]
  \\
  =&
  \frac{1}{\Delta\left(T(2,2a+1);\exp(\xi)\right)}
  +
  \frac{\sqrt{-\pi}}{2\sinh(\xi/2)}
  \sum_{k}
  \exp\left[\frac{N}{\xi}S_{k}(\xi)\right]
  \sqrt{\frac{N}{\xi}}
  \tau_{k},
\end{split}
\end{equation*}
where $\Delta(K;t)$ is the Alexander polynomial of a knot $K$ and
\begin{align}
  S_{k}(\xi)
  &:=
  \frac{-\bigl((2k+1)\pi\sqrt{-1}-(2a+1)\xi\bigr)^2}{2(2a+1)},
  \label{eq:asymptotic_torus_CS}
  \\
  \tau_{k}
  &:=
  (-1)^{k}
  \frac{4\sin\left(\frac{(2k+1)\pi}{2a+1}\right)}{\sqrt{2(2a+1)}}.
  \label{eq:asymptotic_torus_Reidemeister}
\end{align}
Note that $\Delta(T(2,2a+1);t)=\dfrac{\left(t^{2a+1}-t^{-(2a+1)}\right)\left(t^{1/2}-t^{-1/2}\right)}{\left(t^{(2a+1)/2}-t^{-(2a+1)/2}\right)\left(t-t^{-1}\right)}$.
\par
For the range of the index $k$ and a full asymptotic expansion see \cite{Hikami/Murakami:Bonn}.
\subsection{Twice-iterated torus knots}
Let $T(2,2a+1)^{(2,2b+1)}$ be the $(2,2b+1)$-cable of the torus knot of type $(2,2a+1)$ as depicted in Figure~\ref{fig:iterated_torus_knot}, where $a$ is a positive integer and $b$ is an integer such that $2b+1-4(2a+1)>0$.
\subsubsection{The colored Jones polynomial}
We first calculate the Kauffman bracket of the diagram replacing the knot diagram in Figure~\ref{fig:iterated_torus_knot} with the Jones--Wenzl idempotent.
We put
\begin{align*}
  \alpha&:=2a+1, \\
  \beta&:=2b+1.
\end{align*}
By linear skein theory, we have
\begin{align*}
  &
  \left\langle
    \quad
    \raisebox{-22mm}{\includegraphics[scale=0.2]{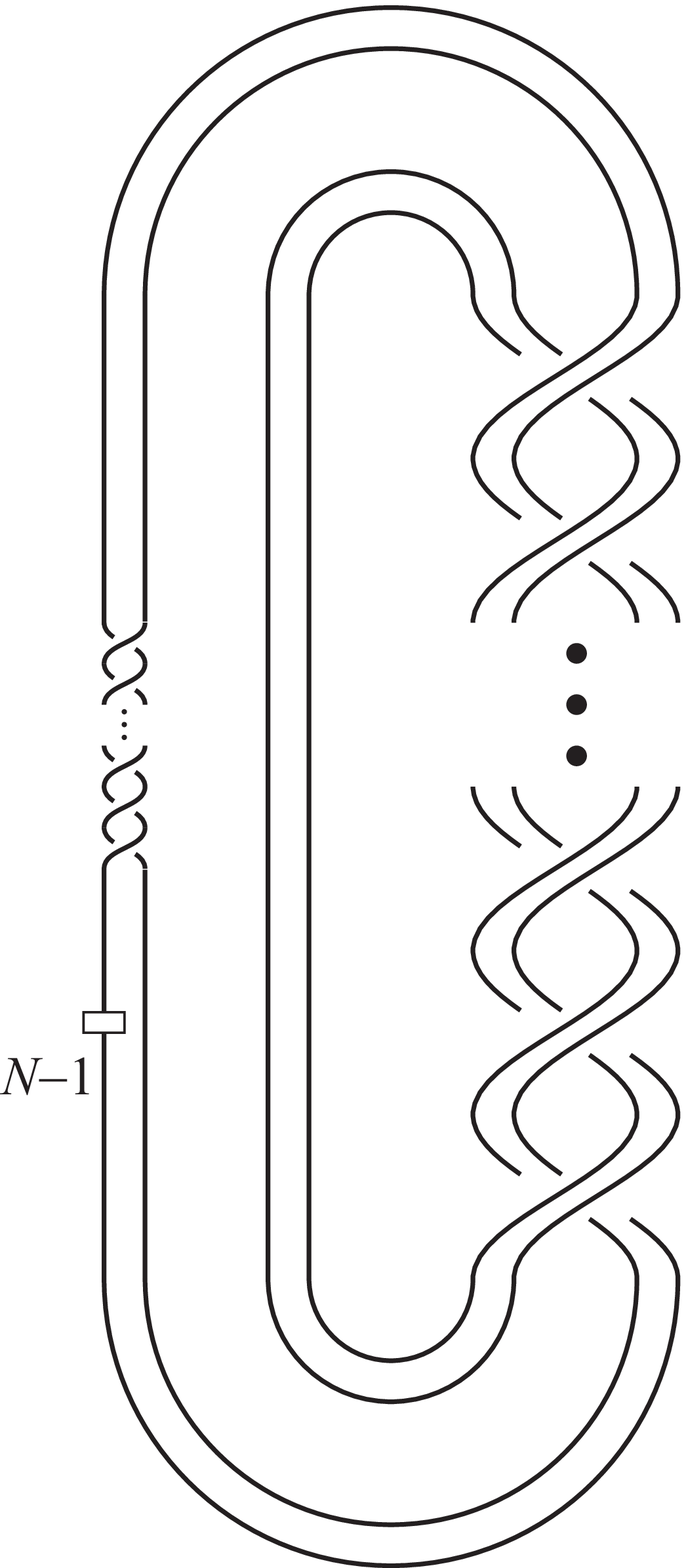}}
    \quad
  \right\rangle
  \\
  =&
  \sum_{d=0}^{N-1}\frac{\Delta_{2d}}{\theta(N-1,N-1,2d)}
  \left\langle
    \quad
    \raisebox{-22mm}{\includegraphics[scale=0.2]{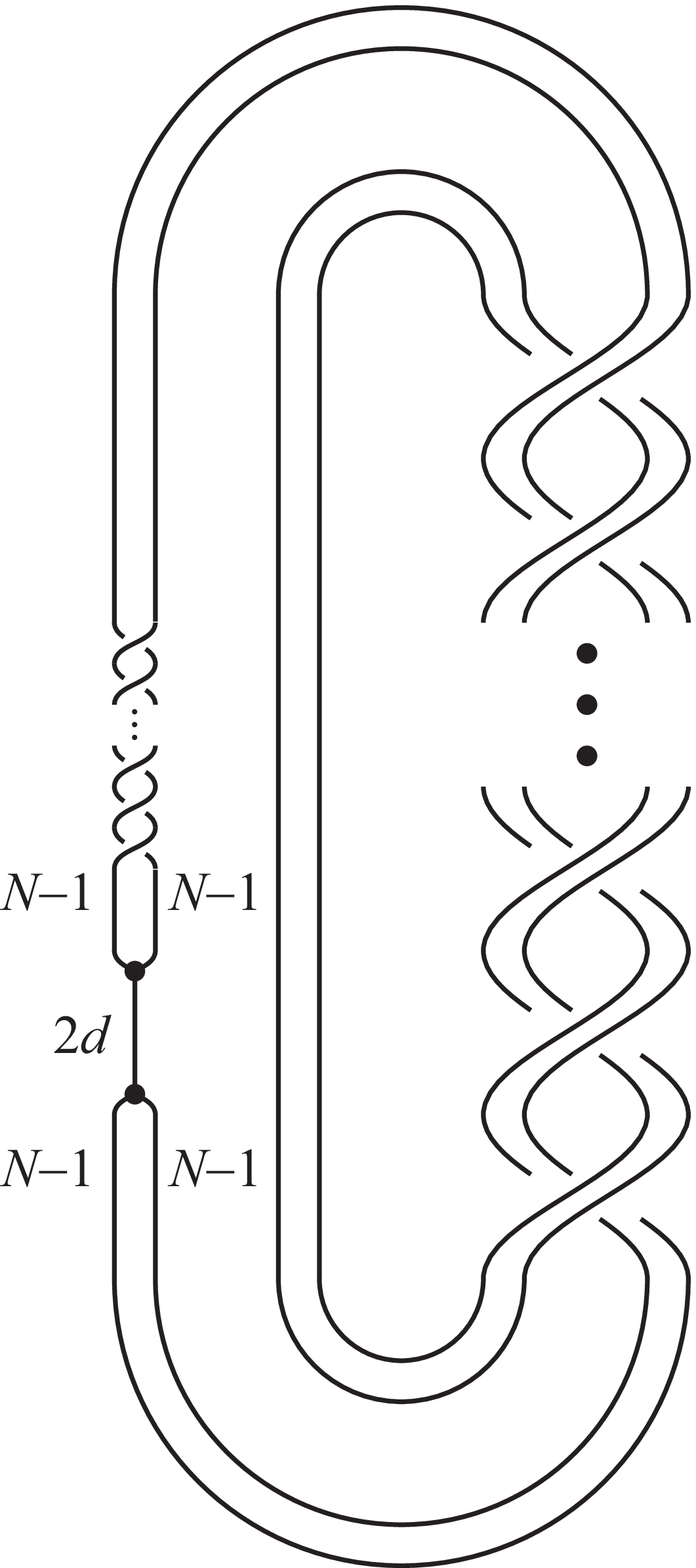}}
    \quad
  \right\rangle
  \\
  =&
  \sum_{d=0}^{N-1}\frac{\Delta_{2d}}{\theta(N-1,N-1,2d)}
  \left((-1)^{d-N+1}A^{-2(N-1)+2d+2d^2-(N-1)^2}\right)^{\beta-2\alpha}
  \left\langle
    \quad
    \raisebox{-22mm}{\includegraphics[scale=0.2]{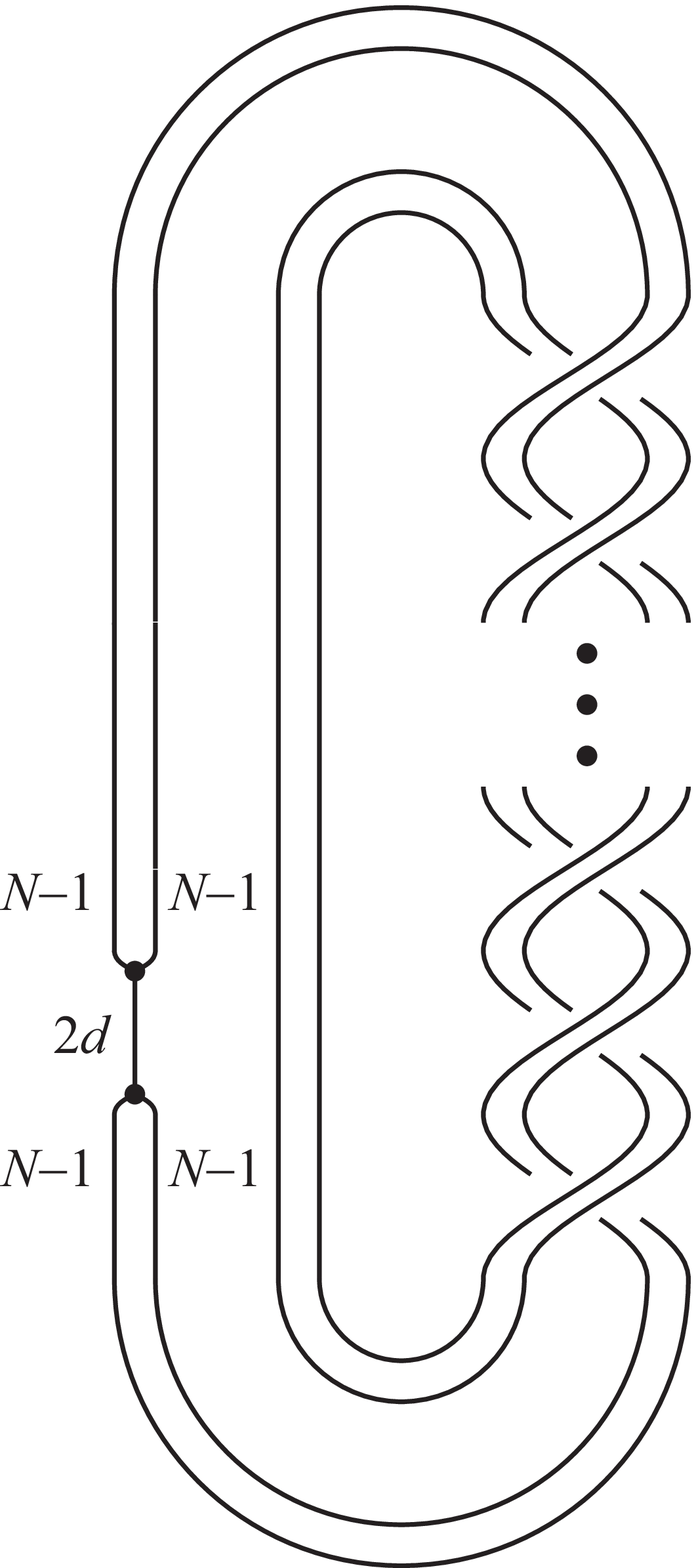}}
    \quad
  \right\rangle
  \\
  =&
  \sum_{d=0}^{N-1}
  (-1)^{d-N+1}A^{(\beta-2\alpha)(2d^2+2d-N^2+1)}
  \frac{\Delta_{2d}}{\theta(N-1,N-1,2d)}
  \left\langle
    \quad
    \raisebox{-22mm}{\includegraphics[scale=0.2]{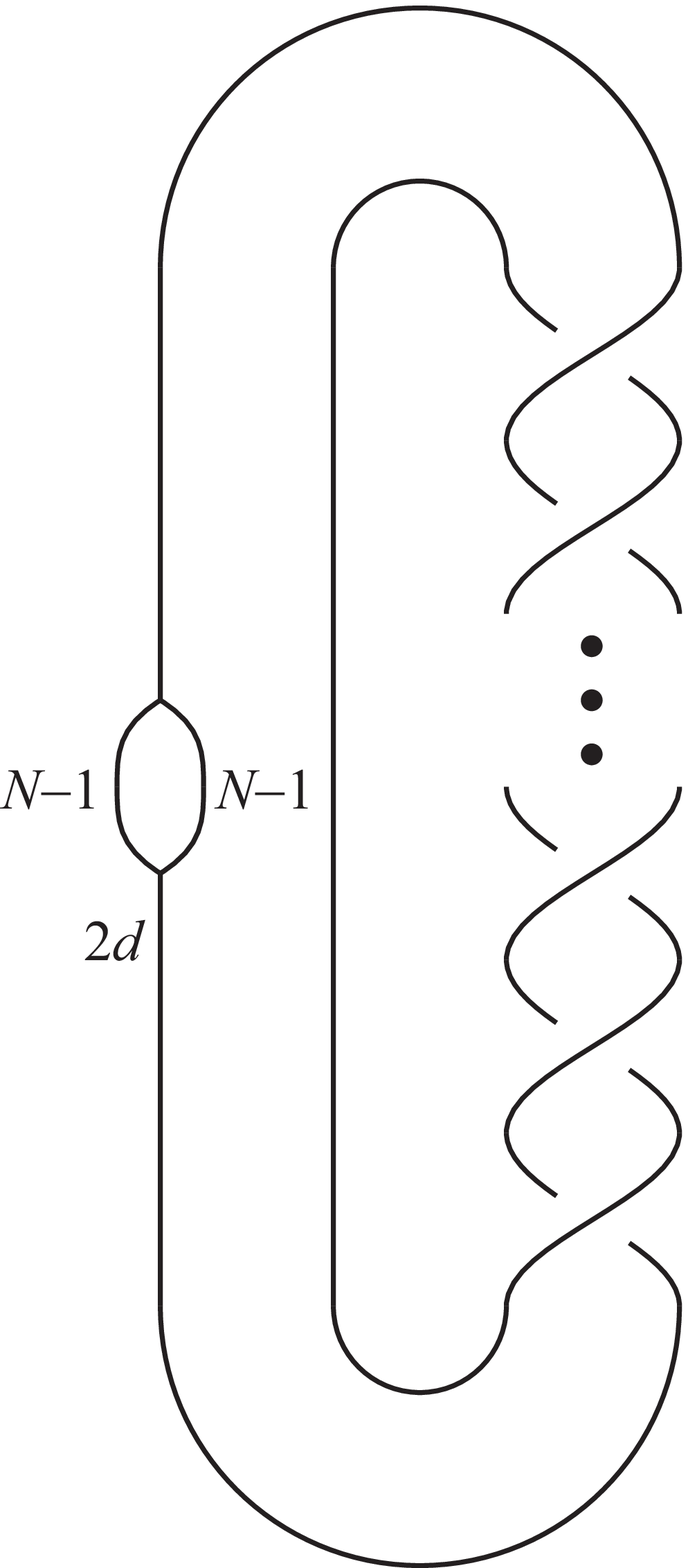}}
    \quad
  \right\rangle
  \\
  =&
  \sum_{d=0}^{N-1}
  (-1)^{d-N+1}A^{(\beta-2\alpha)(2d^2+2d-N^2+1)}
  \left\langle
    \quad
    \raisebox{-22mm}{\includegraphics[scale=0.2]{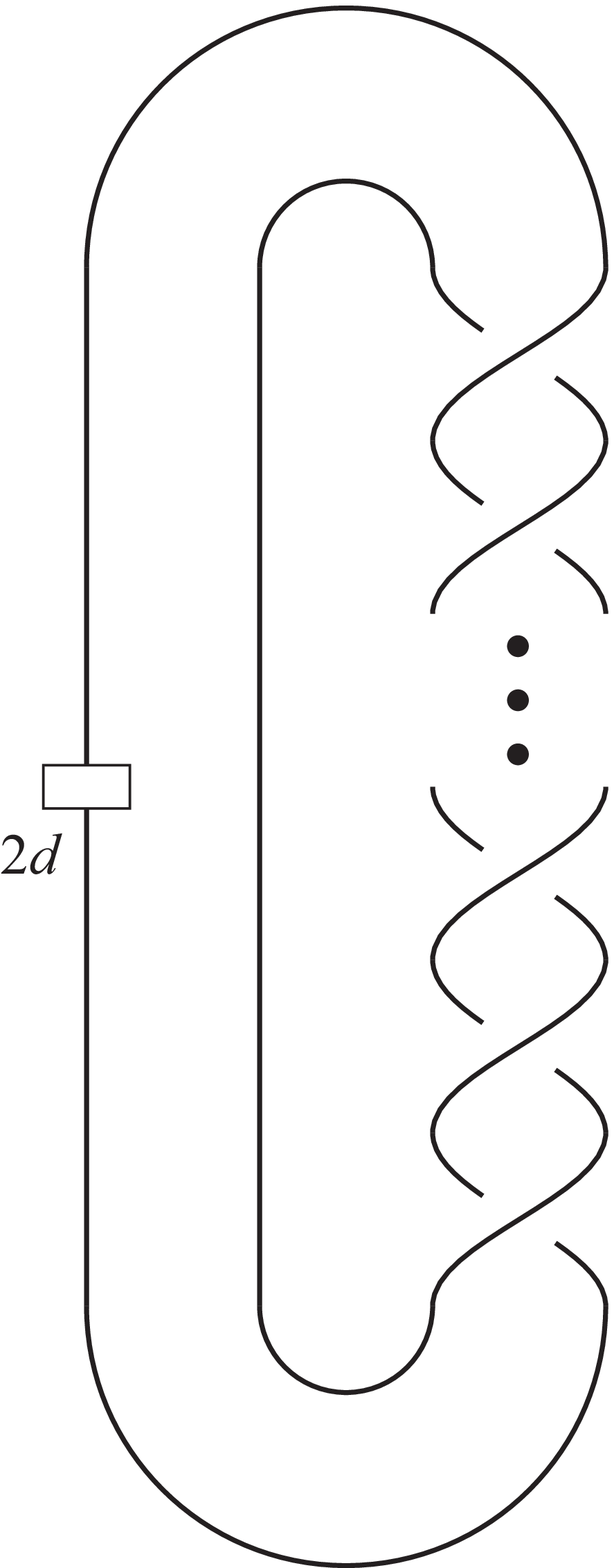}}
    \quad
  \right\rangle
  \\
  =&
  \sum_{d=0}^{N-1}
  (-1)^{d-N+1}A^{(\beta-2\alpha)(2d^2+2d-N^2+1)}
  \\
  &\quad\quad\times
  \left(
    \sum_{c=0}^{2d}
    (-1)^{c-2d}A^{\alpha(2c^2+2c-(2d+1)^2+1)}
    \frac{A^{2(2c+1)}-A^{-2(2c+1)}}{A^2-A^{-2}}
  \right)
  \\
  =&
  \frac{(-1)^{N-1}A^{-(\beta-2\alpha)(N^2-1)}}{A^2-A^{-2}}
  \\
  &\quad\quad\times
  \sum_{d=0}^{N-1}\sum_{c=0}^{2d}
  (-1)^{d+c}A^{2\beta(d^2+d)-8\alpha)(d^2+d)+2\alpha(c^2+c)}
  \left(A^{2(2c+1)}-A^{-2(2c+1)}\right).
\end{align*}
Here the last equality follows from \eqref{eq:torus_knot_idempotent}.
The colored Jones polynomial is given by multiplying $\left((-1)^{N-1}A^{N^2-1}\right)^{-(2b+1)-2(2a+1)}=(-1)^{N-1}A^{-\bigl(2b+1+2(2a+1)\bigr)(N^2-1)}$, dividing by $(-1)^{N-1}\frac{A^{2N}-A^{-2N}}{A^2-A^{-2}}$ and replacing $A$ with $t^{1/4}$.
Therefore we have
\begin{equation*}
\begin{split}
  &J_{N}\left(T(2,\alpha)^{(2,\beta)};t\right)
  \\
  =&
  \frac{(-1)^{N-1}t^{-\beta(N^2-1)/2}}{(t^{N/2}-t^{-N/2})}
  \\
  &\quad
  \times
  \sum_{d=0}^{N-1}\sum_{c=0}^{2d}
  (-1)^{d+c}
  t^{\beta(d^2+d)/2-2\alpha(d^2+d)+\alpha(c^2+c)/2}
  \left(t^{(2c+1)/2}-t^{-(2c+1)/2}\right).
\end{split}
\end{equation*}
\par
See \cite{van_der_Veen:2008} for formulas for general iterated torus knots.
\subsubsection{Asymptotic behavior of the colored Jones polynomial}
Putting $t:=\exp(\xi/N)$ we have
\begin{equation*}
  J_{N}\left(T(2,\alpha)^{(2,\beta)};\exp(\xi/N)\right)
  =
  \frac{(-1)^{N-1}\exp
  \left[
    \frac{-\beta(N^2-1)\xi}{2N}
  \right]}{2\sinh(\xi/2)}
  \left(
    \tilde{\Sigma}_{+}-\tilde{\Sigma}_{-}
  \right),
\end{equation*}
where
\begin{equation*}
\begin{split}
  \tilde{\Sigma}_{\pm}
  &:=
  \sum_{d=0}^{N-1}\sum_{c=0}^{2d}
  (-1)^{d+c}
  \exp\left[\frac{\bigl(\beta(d^2+d)-4\alpha(d^2+d)+\alpha(c^2+c)\pm(2c+1)\bigr)\xi}{2N}\right]
  \\
  &=
  \exp
  \left[
    -\frac{\xi}{2N}
    \left(
      \frac{\beta-4\alpha)}{4}
      +
      \frac{(\alpha\pm2)^2}{4\alpha}
      \mp1
    \right)
  \right]
  \\
  &\quad
  \times
  \sum_{d=0}^{N-1}\sum_{c=0}^{2d}
  (-1)^{d+c}
  \exp
  \left[
    \frac{\xi}{2N}
    \left(
      \alpha\left(c+\frac{\alpha\pm2}{2\alpha}\right)^2
      +
      \bigl(\beta-4\alpha\bigr)\left(d+\frac{1}{2}\right)^2
    \right)
  \right].
\end{split}
\end{equation*}
As in the case of the torus knot, we assume that $\xi$ is not purely imaginary and $\Im\xi\ge0$.
We use the following formula:
\begin{equation*}
  \iint_{C_{x,\varphi}\times C_{y,\psi}}
    e^{-gx^2-hy^2+px+qy}
  \,dx\,dy
  =
  \frac{\pi}{\sqrt{gh}}
  e^{\frac{p^2}{4g}+\frac{q^2}{4h}},
\end{equation*}
where $x\in C_{x,\varphi}$ and $y\in C_{y,\psi}$ with $C_{x,\varphi}$ the line $\left\{te^{\varphi\sqrt{-1}}\mid t\in\R\right\}$ with $\Re\left(ge^{2\varphi\sqrt{-1}}\right)>0$ and $C_{y,\psi}$ the line $\left\{te^{\psi\sqrt{-1}}\mid t\in\R\right\}$ with $\Re\left(he^{2\psi\sqrt{-1}}\right)>0$.
Note that the double integral converges absolutely.
Putting
\begin{align*}
  g&:=\frac{N}{2(2a+1)\xi}, \\
  h&:=\frac{N}{2\bigl(2b+1-4(2a+1)\bigr)\xi},\\
  p&:=c+\frac{2a+1\pm2}{2(2a+1)},\\
  q&:=d+\frac{1}{2}
\end{align*}
we have
\begin{equation*}
\begin{split}
  &\tilde{\Sigma}_{\pm}
  \\
  =&
  \frac{N}{2\xi\pi\sqrt{\alpha(\beta-4\alpha)}}
  \exp
  \left[
    -\frac{\xi}{2N}
    \left(
      \frac{\beta-4\alpha}{4}
      +
      \frac{\alpha}{4}+\frac{1}{\alpha}
    \right)
  \right]
  \\
  &
  \times
  \sum_{d=0}^{N-1}\sum_{c=0}^{2d}
  (-1)^{d+c}
  \\
  &\quad
  \iint_{C_{x,\varphi}\times C_{y,\psi}}
  \exp
  \left[
    -\frac{Nx^2}{2\alpha\xi}
    -\frac{Ny^2}{2(\beta-4\alpha)\xi}
    +\left(c+\frac{\alpha\pm2}{2\alpha}\right)x
    +\left(d+\frac{1}{2}\right)y
  \right]
  \,dx\,dy.
\end{split}
\end{equation*}
Note that we need to choose $\varphi$ so that $\Re\left(e^{2\varphi\sqrt{-1}}\xi^{-1}\right)>0$ and $\psi$ so that $\Re\left((\beta-4\alpha)e^{2\psi\sqrt{-1}}\xi^{-1}\right)>0$.
The summation becomes
\begin{equation*}
\begin{split}
  &\sum_{d=0}^{N-1}
  (-1)^{d}
  \\
  &\times
  \iint_{C_{x,\varphi}\times C_{y,\psi}}
  \exp
  \left[
    -\frac{Nx^2}{2\alpha\xi}
    -\frac{Ny^2}{2(\beta-4\alpha)\xi}
    +\left(\frac{\alpha\pm2}{2\alpha}\right)x
    +\left(d+\frac{1}{2}\right)y
  \right]
  \left(
    \sum_{c=0}^{2d}(-e^x)^c
  \right)
  \,dx\,dy
  \\
  =
  &\sum_{d=0}^{N-1}
  (-1)^{d}
  \\
  &\times
  \iint_{C_{x,\varphi}\times C_{y,\psi}}
  \exp
  \left[
    -\frac{Nx^2}{2\alpha\xi}
    -\frac{Ny^2}{2(\beta-4\alpha)\xi}
    +\frac{x}{2}
    +\frac{y}{2}
  \right]
  \exp\left[\frac{\pm x}{\alpha}\right]
  \frac{e^{dy}\left(1+e^{(2d+1)x}\right)}{1+e^x}
  \,dx\,dy.
\end{split}
\end{equation*}
Therefore we have
\begin{equation*}
\begin{split}
  &\tilde{\Sigma}_{+}-\tilde{\Sigma}_{-}
  \\
  =&
  \frac{N}{\xi\pi\sqrt{\alpha(\beta-4\alpha)}}
  \exp
  \left[
    -\frac{\xi}{2N}
    \left(
      \frac{\beta-4\alpha}{4}
      +
      \frac{\alpha}{4}+\frac{1}{\alpha}
    \right)
  \right]
  \\
  &
  \times
  \iint_{C_{x,\varphi}\times C_{y,\psi}}
  \exp
  \left[
    -\frac{Nx^2}{2\alpha\xi}
    -\frac{Ny^2}{2(\beta-4\alpha)\xi}
    +\frac{x}{2}
    +\frac{y}{2}
  \right]
  \\
  &\quad\times
  \frac{\sinh\left(\frac{x}{\alpha}\right)}{1+e^x}
  \left(
    \sum_{d=0}^{N-1}
    (-1)^{d}e^{dy}\left(1+e^{(2d+1)x}\right)
  \right)
  \,dx\,dy
  \\
  =&
  \frac{N}{\xi\pi\sqrt{\alpha(\beta-4\alpha)}}
  \exp
  \left[
    -\frac{\xi}{2N}
    \left(
      \frac{\beta-4\alpha}{4}
      +
      \frac{\alpha}{4}+\frac{1}{\alpha}
    \right)
  \right]
  \\
  &
  \times
  \iint_{C_{x,\varphi}\times C_{y,\psi}}
  \exp
  \left[
    -\frac{Nx^2}{2\alpha\xi}
    -\frac{Ny^2}{2(\beta-4\alpha)\xi}
    +\frac{x}{2}
    +\frac{y}{2}
  \right]
  \\
  &\quad\times
  \frac{\sinh\left(\frac{x}{\alpha}\right)}{1+e^x}
  \left(
    \frac{1-(-1)^Ne^{Ny}}{1+e^y}
    +
    \frac{e^x\left(1-(-1)^Ne^{N(2x+y)}\right)}{1+e^{2x+y}}
  \right)
  \,dx\,dy
  \\
  =&
  \frac{N}{4\xi\pi\sqrt{\alpha(\beta-4\alpha)}}
  \exp
  \left[
    -\frac{\xi}{2N}
    \left(
      \frac{\beta-4\alpha}{4}
      +
      \frac{\alpha}{4}+\frac{1}{\alpha}
    \right)
  \right]
  \left(
    I_1-(-1)^{N}I_2-(-1)^{N}I_3
  \right)
\end{split}
\end{equation*}
where
\begin{align*}
  I_1
  &:=
  \iint_{C_{x,\varphi}\times C_{y,\psi}}
  \exp
  \left[
    -\frac{Nx^2}{2\alpha\xi}
    -\frac{Ny^2}{2(\beta-4\alpha)\xi}
  \right]
  \\
  &\quad\times
  \frac{\sinh\left(\frac{x}{\alpha}\right)}{\cosh\left(\frac{x}{2}\right)}
  \left(
    \frac{1}{\cosh\left(\frac{y}{2}\right)}
    +
    \frac{1}{\cosh\left(\frac{2x+y}{2}\right)}
  \right)
  \,dx\,dy,
  \\
  I_2
  &:=
  \iint_{C_{x,\varphi}\times C_{y,\psi}}
  \exp
  \left[
    -\frac{Nx^2}{2\alpha\xi}
    -\frac{Ny^2}{2(\beta-4\alpha)\xi}
    +Ny
  \right]
  \frac{\sinh\left(\frac{x}{\alpha}\right)}
       {\cosh\left(\frac{x}{2}\right)\cosh\left(\frac{y}{2}\right)}
  \,dx\,dy,
  \\
  I_3
  &:=
  \iint_{C_{x,\varphi}\times C_{y,\psi}}
  \exp
  \left[
    -\frac{Nx^2}{2\alpha\xi}
    -\frac{Ny^2}{2(\beta-4\alpha)\xi}
    +2Nx
    +Ny
  \right]
  \\
  &\quad\times
  \frac{\sinh\left(\frac{x}{\alpha}\right)}
       {\cosh\left(\frac{x}{2}\right)\cosh\left(\frac{2x+y}{2}\right)}
  \,dx\,dy.
\end{align*}
We apply the saddle point method (Theorem~\ref{thm:saddle}) to obtain the asymptotic behaviors of $I_1$, $I_2$ and $I_3$.
\par
Since we assume that $\beta-4\alpha>0$, we can put $\psi=\varphi$.
We also assume that $\varphi\ne\pm\pi/2$ so that the denominators of the integrands in $I_1$, $I_2$ and $I_3$ do not vanish.
\begin{rem}
The condition $\beta-4\alpha>0$ is also a condition that the iterated torus knot $T(2,\alpha)^{(2,\beta)}$ becomes the link of a singularity \cite[Appendix to Chapter~I]{Eisenbud/Neumann:1985}.
I wish to thank Roland van der Veen for pointing this out.
\end{rem}
\par
We first calculate $I_3$.
Put
\begin{equation*}
  J(x)
  :=
  \int_{C_{\varphi}}
  \frac{\exp\left[N\left(-\frac{y^2}{2(\beta-4\alpha)\xi}+y\right)\right]}
       {\cosh\left(\frac{2x+y}{2}\right)}
  \,dy
\end{equation*}
so that
\begin{equation*}
  I_3
  =
  \int_{C_{\varphi}}
  \exp
  \left[N\left(-\frac{x^2}{2\alpha\xi}+2x\right)\right]
  \frac{\sinh\left(\frac{x}{\alpha}\right)}{\cosh\left(\frac{x}{2}\right)}
  J(x)
  \,dx.
\end{equation*}
\begin{equation*}
\end{equation*}
By the residue theorem, we have
\begin{equation*}
\begin{split}
  J(x)
  &=
  \int_{C_{\varphi}+(\beta-4\alpha)\xi}
  \frac{\exp\left[N\left(-\frac{y^2}{2(\beta-4\alpha)\xi}+y\right)\right]}
       {\cosh\left(\frac{2x+y}{2}\right)}
  \,dy
  \\
  &\quad
  +
  2\pi\sqrt{-1}
  \sum_{j}\Res
  \left(
    \frac{\exp\left[N\left(-\frac{y^2}{2(\beta-4\alpha)\xi}+y\right)\right]}
         {\cosh\left(\frac{2x+y}{2}\right)}
    ;y=-2x+(2j+1)\pi\sqrt{-1}
  \right)
  \\
  &
  \text{(Put $z:=y-(\beta-4\alpha)\xi$)}
  \\
  &=
  \exp\left[\frac{N\xi}{2}(\beta-4\alpha)\right]
  \int_{C_{\varphi}}
  \frac{\exp\left[N\left(-\frac{z^2}{2(\beta-4\alpha)\xi}\right)\right]}
       {\cosh\left(\frac{2x+z+(\beta-4\alpha)\xi}{2}\right)}
  \,dz
  \\
  &\quad
  +
  4\pi
  \sum_{j}(-1)^{j}
  \exp
  \left[
    N\left(-\frac{\bigl(-2x+(2j+1)\pi\sqrt{-1}\bigr)^2}{2(\beta-4\alpha)\xi}
    -2x+(2j+1)\pi\sqrt{-1}\right)
  \right],
\end{split}
\end{equation*}
where $C_{\varphi}+(\beta-4\alpha\xi)$ is the line $\left\{te^{\varphi\sqrt{-1}}+(\beta-4\alpha)\xi\mid t\in\R\right\}$, and $j$ runs over integers such that $-2x+(2j+1)\pi\sqrt{-1}$ is between $C_{\varphi}$ and $C_{\varphi}+(\beta-4\alpha)\xi$.
Note that since $x\in C_{\varphi}$ the range does not depend on $x$.
Therefore we have
\begin{equation*}
\begin{split}
  &I_3
  \\
  =&
  e^{\frac{N\xi}{2}(\beta-4\alpha)}
  \iint_{C_{\varphi}\times C_{\varphi}}
  \frac{\sinh\left(\frac{x}{\alpha}\right)
        \exp\left[N\left(-\frac{z^2}{2(\beta-4\alpha)\xi}-\frac{x^2}{2\alpha\xi}+2x\right)\right]}
       {\cosh\left(\frac{x}{2}\right)\cosh\left(\frac{2x+z+(\beta-4\alpha)\xi}{2}\right)}
  \,dz\,dx
  \\
  &+
  4\pi
  \sum_{j}(-1)^{j}
  e^{N(2j+1)\pi\sqrt{-1}}
  \int_{C_{\varphi}}
  \frac{\sinh\left(\frac{x}{\alpha}\right)}{\cosh\left(\frac{x}{2}\right)}
  \exp
  \left[
    N\left(-\frac{\bigl(-2x+(2j+1)\pi\sqrt{-1}\bigr)^2}{2(\beta-4\alpha)\xi}
    -\frac{x^2}{2\alpha\xi}\right)
  \right]
  \,dx
  \\
  =&
  e^{\frac{N\xi}{2}(\beta-4\alpha)}
  \int_{C_{\varphi}}
  \exp\left[N\left(-\frac{z^2}{2(\beta-4\alpha)\xi}\right)\right]
  \left(
    \int_{C_{\varphi}}
    \frac{\sinh\left(\frac{x}{\alpha}\right)
          \exp\left[N\left(-\frac{x^2}{2\alpha\xi}+2x\right)\right]}
         {\cosh\left(\frac{x}{2}\right)\cosh\left(\frac{2x+z+(\beta-4\alpha)\xi}{2}\right)}
    \,dx
  \right)
  \,dz
  \\
  &+
  4\pi
  \sum_{j}(-1)^{j}
  e^{N(2j+1)\pi\sqrt{-1}}
  \int_{C_{\varphi}}
  \frac{\sinh\left(\frac{x}{\alpha}\right)}{\cosh\left(\frac{x}{2}\right)}
  \exp
  \left[
    N\left(-\frac{\bigl(-2x+(2j+1)\pi\sqrt{-1}\bigr)^2}{2(\beta-4\alpha)\xi}
    -\frac{x^2}{2\alpha\xi}\right)
  \right]
  \,dx.
\end{split}
\end{equation*}
Put
\begin{equation*}
  K(z)
  :=
  \int_{C_{\varphi}}
  \frac{\sinh\left(\frac{x}{\alpha}\right)
        \exp\left[N\left(-\frac{x^2}{2\alpha\xi}+2x\right)\right]}
       {\cosh\left(\frac{x}{2}\right)\cosh\left(\frac{2x+z+(\beta-4\alpha)\xi}{2}\right)}
  \,dx.
\end{equation*}
Then we have
\begin{equation*}
\begin{split}
  &K(z)
  \\
  =&
  \int_{C_{\varphi}+2\alpha\xi}
  \frac{\sinh\left(\frac{x}{\alpha}\right)
        \exp\left[N\left(-\frac{x^2}{2\alpha\xi}+2x\right)\right]}
       {\cosh\left(\frac{x}{2}\right)\cosh\left(\frac{2x+z+(\beta-4\alpha)\xi}{2}\right)}
  \,dx
  \\
  &+
  2\pi\sqrt{-1}
  \sum_{k}
  \Res
    \left(
    \frac{\sinh\left(\frac{x}{\alpha}\right)
          \exp\left[N\left(-\frac{x^2}{2\alpha\xi}+2x\right)\right]}
         {\cosh\left(\frac{x}{2}\right)\cosh\left(\frac{2x+z+(\beta-4\alpha)\xi}{2}\right)}
    ;x=(2k+1)\pi\sqrt{-1}
  \right)
  \\
  &+
  2\pi\sqrt{-1}
  \sum_{l}
  \Res
    \left(
    \frac{\sinh\left(\frac{x}{\alpha}\right)
          \exp\left[N\left(-\frac{x^2}{2\alpha\xi}+2x\right)\right]}
         {\cosh\left(\frac{x}{2}\right)\cosh\left(\frac{2x+z+(\beta-4\alpha)\xi}{2}\right)}
    ;x=\frac{1}{2}\bigl((2l+1)\pi\sqrt{-1}-z-(\beta-4\alpha)\xi\bigr)
  \right)
  \\
  &\text{(Put $w:=x-2\alpha\xi$)}
  \\
  =&
  e^{2N\alpha\xi}
  \int_{C_{\varphi}}
  \frac{\sinh\left(\frac{w+2\alpha\xi}{\alpha}\right)
        \exp\left[N\left(-\frac{w^2}{2\alpha\xi}\right)\right]}
       {\cosh\left(\frac{w+2\alpha\xi}{2}\right)
        \cosh\left(\frac{2w+z+\beta\xi}{2}\right)}
  \,dw
  \\
  &+
  4\pi\sqrt{-1}
  \sum_{k}
  (-1)^{k}
  \frac{\sin\left(\frac{(2k+1)\pi}{\alpha}\right)
        \exp
        \left[N\left(\frac{(2k+1)^2\pi^2}{2\alpha\xi}+2(2k+1)\pi\sqrt{-1}\right)\right]}
       {\cosh\left(\frac{2(2k+1)\pi\sqrt{-1}+z+(\beta-4\alpha)\xi}{2}\right)}
  \\
  &+
  2\pi
  \sum_{l}
  (-1)^{l}
  \frac{\sinh\left(\frac{(2l+1)\pi\sqrt{-1}-z-(\beta-4\alpha)\xi}{2\alpha}\right)}
       {\cosh\left(\frac{(2l+1)\pi\sqrt{-1}-z-(\beta-4\alpha)\xi}{4}\right)}
  \\
  &\quad\times
  \exp
  \left[
    N
    \left(
      -\frac{\bigl((2l+1)\pi\sqrt{-1}-z-(\beta-4\alpha)\xi\bigr)^2}{8\alpha\xi}
      +(2l+1)\pi\sqrt{-1}-z-(\beta-4\alpha)\xi
    \right)
  \right],
\end{split}
\end{equation*}
where $k$ runs over integers such that $\alpha\nmid(2k+1)$ and that $(2k+1)\pi\sqrt{-1}$ is between $C_{\varphi}$ and $C_{\varphi}+2\alpha\xi$, and $l$ runs over integers such that $\frac{1}{2}\bigl((2l+1)\pi\sqrt{-1}-z-(\beta-4\alpha)\xi\bigr)$ is between $C_{\varphi}$ and $C_{\varphi}+2\alpha\xi$.
\par
Therefore we have
\begin{equation}\label{eq:I_3}
\begin{split}
  &I_3
  \\
  =&
  e^{\frac{N\xi}{2}\beta}
  \int_{C_{\varphi}\times C_{\varphi}}
  \frac{\sinh\left(\frac{w+2\alpha\xi}{\alpha}\right)
        \exp\left[N\left(-\frac{w^2}{2\alpha\xi}-\frac{z^2}{2(\beta-4\alpha)\xi}\right)\right]}
       {\cosh\left(\frac{w+2\alpha\xi}{2}\right)
        \cosh\left(\frac{2w+z+\beta\xi}{2}\right)}
  \,dw\,dz
  \\
  &+
  4\pi\sqrt{-1}
  e^{\frac{N\xi}{2}(\beta-4\alpha)}
  \sum_{k}
  (-1)^{k}
  \sin\left(\frac{(2k+1)\pi}{\alpha}\right)
  \exp\left[N\left(\frac{(2k+1)^2\pi^2}{2\alpha\xi}+2(2k+1)\pi\sqrt{-1}\right)\right]
  \\
  &\quad\times
  \int_{C_{\varphi}}
  \frac{\exp\left[N\left(-\frac{z^2}{2(\beta-4\alpha)\xi}\right)\right]}
       {\cosh\left(\frac{2(2k+1)\pi\sqrt{-1}+z+(\beta-4\alpha)\xi}{2}\right)}
  \,dz
  \\
  &+
  2\pi
  e^{\frac{N\xi}{2}(\beta-4\alpha)}
  \sum_{l}
  (-1)^{l}
  e^{N\bigl((2l+1)\pi\sqrt{-1}-(\beta-4\alpha)\xi\bigr)}
  \int_{C_{\varphi}}
  \frac{
        \sinh\left(\frac{(2l+1)\pi\sqrt{-1}-z-(\beta-4\alpha)\xi}{2\alpha}\right)}
       {\cosh\left(\frac{(2l+1)\pi\sqrt{-1}-z-(\beta-4\alpha)\xi}{4}\right)}
  \\
  &\quad\times
  \exp
  \left[
    N
    \left(
      -\frac{\bigl((2l+1)\pi\sqrt{-1}-z-(\beta-4\alpha)\xi\bigr)^2}{8\alpha\xi}
      -z-\frac{z^2}{2(\beta-4\alpha)\xi}
    \right)
  \right]
  \,dz
  \\
  &+
  4\pi
  \sum_{j}(-1)^{j}
  e^{N(2j+1)\pi\sqrt{-1}}
  \int_{C_{\varphi}}
  \frac{\sinh\left(\frac{x}{\alpha}\right)}{\cosh\left(\frac{x}{2}\right)}
  \exp
  \left[
    N\left(-\frac{\bigl(-2x+(2j+1)\pi\sqrt{-1}\bigr)^2}{2(\beta-4\alpha)\xi}
    -\frac{x^2}{2\alpha\xi}\right)
  \right]
  \,dx
  \\
  =&
  e^{\frac{N\xi}{2}\beta}
  \int_{C_{\varphi}\times C_{\varphi}}
  \frac{\sinh\left(\frac{w+2\alpha\xi}{\alpha}\right)
        \exp\left[N\left(-\frac{w^2}{2\alpha\xi}-\frac{z^2}{2(\beta-4\alpha)\xi}\right)\right]}
       {\cosh\left(\frac{w+2\alpha\xi}{2}\right)
        \cosh\left(\frac{2w+z+\beta\xi}{2}\right)}
  \,dw\,dz
  \\
  &+
  4\pi\sqrt{-1}
  e^{\frac{N\xi}{2}(\beta-4\alpha)}
  \sum_{k}
  (-1)^{k+1}
  \sin\left(\frac{(2k+1)\pi}{\alpha}\right)
  \exp\left[N\left(\frac{(2k+1)^2\pi^2}{2\alpha\xi}+2(2k+1)\pi\sqrt{-1}\right)\right]
  \\
  &\quad\times
  \int_{C_{\varphi}}
  \frac{\exp\left[N\left(-\frac{z^2}{2(\beta-4\alpha)\xi}\right)\right]}
       {\cosh\left(\frac{z+(\beta-4\alpha)\xi}{2}\right)}
  \,dz
  \\
  &+
  2\pi
  \sum_{l}
  (-1)^{l}
  e^{N\left(\frac{(2l+1)^2\pi^2}{2\beta\xi}+(2l+1)\pi\sqrt{-1}\right)}
  \int_{C_{\varphi}}
  \frac{\sinh\left(\frac{(2l+1)\pi\sqrt{-1}-z-(\beta-4\alpha)\xi}{2\alpha}\right)}
       {\cosh\left(\frac{(2l+1)\pi\sqrt{-1}-z-(\beta-4\alpha)\xi}{4}\right)}
  \\
  &\quad\times
  \exp
  \left[
    N
    \left(
      -\frac{\beta}{8\alpha(\beta-4\alpha)\xi}
      \left(z-\frac{(\beta-4\alpha)\bigl((2l+1)\pi\sqrt{-1}-\beta\xi\bigr)}{\beta}\right)^2
    \right)
  \right]
  \,dz
  \\
  &+
  4\pi
  \sum_{j}(-1)^{j}
  e^{N\left(\frac{(2j+1)^2\pi^2}{2\beta\xi}+(2j+1)\pi\sqrt{-1}\right)}
  \int_{C_{\varphi}}
  \frac{\sinh\left(\frac{x}{\alpha}\right)}{\cosh\left(\frac{x}{2}\right)}
  \\
  &\quad\times
  \exp
  \left[
    N
    \left(
      -\frac{\beta}{2\alpha(\beta-4\alpha)\xi}
      \left(x-\frac{2(2j+1)\alpha\pi\sqrt{-1}}{\beta}\right)^2
    \right)
  \right]
  \,dx.
\end{split}
\end{equation}
By the residue theorem, the integral in the third term of \eqref{eq:I_3} becomes
\begin{equation*}
\begin{split}
  &\int_{C_{\varphi}+\frac{(\beta-4\alpha)\bigl((2l+1)\pi\sqrt{-1}-\beta\xi\bigr)}{\beta}}
  f(z)\,dz
  \\
  -&
  2\pi\sqrt{-1}
  \sum_{m}
  \Res\bigl(f(z);z=(2l+1)\pi\sqrt{-1}-2(2m+1)\pi\sqrt{-1}-(\beta-4\alpha)\xi\bigr),
\end{split}
\end{equation*}
where
\begin{multline*}
  f(z)
  :=
  \frac{\sinh\left(\frac{(2l+1)\pi\sqrt{-1}-z-(\beta-4\alpha)\xi}{2\alpha}\right)}
       {\cosh\left(\frac{(2l+1)\pi\sqrt{-1}-z-(\beta-4\alpha)\xi}{4}\right)}
  \\
  \times
  \exp
  \left[
    N
    \left(
      -\frac{\beta}{8\alpha(\beta-4\alpha)\xi}
      \left(z-\frac{(\beta-4\alpha)\bigl((2l+1)\pi\sqrt{-1}-\beta\xi\bigr)}{\beta}\right)^2
    \right)
  \right]
\end{multline*}
and $m$ runs over integers such that $\alpha\nmid(2m+1)$ and that $(2l+1)\pi\sqrt{-1}-2(2m+1)\pi\sqrt{-1}-(\beta-4\alpha)\xi$ is between $C_{\varphi}+\frac{(\beta-4\alpha)\bigl((2l+1)\pi\sqrt{-1}-\beta\xi\bigr)}{\beta}$ and $C_{\varphi}$.
Note the minus sign in front of the second term.
This is because the line $C_{\varphi}+\frac{(\beta-4\alpha)\bigl((2l+1)\pi\sqrt{-1}-\beta\xi\bigr)}{\beta}$ is below the line $C_{\varphi}$, since $\frac{1}{2}\bigl((2l+1)\sqrt{-1}-z-(\beta-4\alpha)\xi\bigr)$ is between $C_{\varphi}$ and $C_{\varphi}+2\alpha\xi$.
Putting $w:=z-\frac{(\beta-4\alpha)\bigl((2l+1)\pi\sqrt{-1}-\beta\xi\bigr)}{\beta}$, the integral above becomes
\begin{equation*}
\begin{split}
  &\int_{C_{\varphi}}
  \frac{\sinh\left(\frac{2(2l+1)\pi\sqrt{-1}}{\beta}-\frac{w}{2\alpha}\right)}
       {\cosh\left(\frac{(2l+1)\alpha\pi\sqrt{-1}}{\beta}-\frac{w}{4}\right)}
  \exp\left[N\left(-\frac{\beta}{8\alpha(\beta-4\alpha)\xi}w^2\right)\right]
  \,dw
  \\
  +&
  8\pi
  \sum_{m}
  (-1)^{m}
  \sin\left(\frac{(2m+1)\pi}{\alpha}\right)
  \exp
  \left[
    N\times\frac{\bigl(2(2l+1)\alpha-(2m+1)\beta\bigr)^2\pi^2}{2\alpha\beta(\beta-4\alpha)\xi}
  \right].
\end{split}
\end{equation*}
The fourth term of \eqref{eq:I_3} becomes
\begin{equation*}
  \int_{C_{\varphi}+\frac{2(2j+1)\alpha\pi\sqrt{-1}}{\beta}}
  g(x)\,dx
  +
  2\pi\sqrt{-1}
  \sum_{n}
  \Res(g(x);x=(2n+1)\pi\sqrt{-1}),
\end{equation*}
where
\begin{equation*}
  g(x)
  :=
  \frac{\sinh\left(\frac{x}{\alpha}\right)}{\cosh\left(\frac{x}{2}\right)}
  \exp
  \left[
    N
    \left(
      -\frac{\beta}{2\alpha(\beta-4\alpha)\xi}
      \left(x-\frac{2(2j+1)\alpha\pi\sqrt{-1}}{\beta}\right)^2
    \right)
  \right]
\end{equation*}
and $n$ runs over integers such that $\alpha\nmid(2n+1)$ and that $(2n+1)\pi\sqrt{-1}$ is between $C_{\varphi}$ and $C_{\varphi}+\frac{2(2j+1)\alpha\pi\sqrt{-1}}{\beta}$.
Putting $z:=x-\frac{2(2j+1)\alpha\pi\sqrt{-1}}{\beta}$, the integral becomes
\begin{equation*}
\begin{split}
  &\int_{C_{\varphi}}
  \frac{\sinh\left(\frac{z}{\alpha}+\frac{2(2j+1)\pi\sqrt{-1}}{\beta}\right)}
       {\cosh\left(\frac{z}{2}+\frac{(2j+1)\alpha\pi\sqrt{-1}}{\beta}\right)}
  \exp
  \left[
    N
    \left(
      -\frac{\beta}{2\alpha(\beta-4\alpha)\xi}z^2
    \right)
  \right]
  \,dz
  \\
  +&
  4\pi
  \sum_{n}
  (-1)^{n}
  \sin\left(\frac{(2n+1)\pi}{\alpha}\right)
  \exp
  \left[
    N
    \left(
      \frac{\bigl((2n+1)\beta-2(2j+1)\alpha\bigr)^2\pi^2}{2\alpha\beta(\beta-4\alpha)\xi}
    \right)
  \right].
\end{split}
\end{equation*}
\par
So we have
\begin{equation}\label{eq:I_3_2}
\begin{split}
  &I_3
  \\
  =&
  e^{\frac{N\xi}{2}\beta}
  \int_{C_{\varphi}\times C_{\varphi}}
  \frac{\sinh\left(\frac{w+2\alpha\xi}{\alpha}\right)
        e^{N\left(-\frac{w^2}{2\alpha\xi}-\frac{z^2}{2(\beta-4\alpha)\xi}\right)}}
       {\cosh\left(\frac{w+2\alpha\xi}{2}\right)
        \cosh\left(\frac{2w+z+\beta\xi}{2}\right)}
  \,dw\,dz
  \\
  &+
  4\pi\sqrt{-1}
  e^{\frac{N\xi}{2}(\beta-4\alpha)}
  \sum_{k}
  (-1)^{k+1}
  \sin\left(\frac{(2k+1)\pi}{\alpha}\right)
  e^{N\left(\frac{(2k+1)^2\pi^2}{2\alpha\xi}+2(2k+1)\pi\sqrt{-1}\right)}
  \\
  &\quad\times
  \int_{C_{\varphi}}
  \frac{e^{N\left(-\frac{z^2}{2(\beta-4\alpha)\xi}\right)}}{\cosh\left(\frac{z+(\beta-4\alpha)\xi}{2}\right)}
  \,dz
  \\
  &+
  2\pi
  \sum_{l}
  (-1)^{l}
  e^{N\left(\frac{(2l+1)^2\pi^2}{2\beta\xi}+(2l+1)\pi\sqrt{-1}\right)}
  \int_{C_{\varphi}}
  \frac{\sinh\left(\frac{2(2l+1)\pi\sqrt{-1}}{\beta}-\frac{w}{2\alpha}\right)}
       {\cosh\left(\frac{(2l+1)\alpha\pi\sqrt{-1}}{\beta}-\frac{w}{4}\right)}
  e^{N\left(-\frac{\beta}{8\alpha(\beta-4\alpha)\xi}w^2\right)}
  \,dw
  \\
  &+
  16\pi^2
  \sum_{l,m}
  (-1)^{l+m}
  e^{
    N
    \left(
      \frac{\pi^2}{2\alpha(\beta-4\alpha)\xi}
      \bigl((2l+1)^2\alpha+(2m+1)^2\beta-4(2l+1)(2m+1)\alpha\bigr)
      +(2l+1)\pi\sqrt{-1}
    \right)
  }
  \\
  &\quad\times
  \sin\left(\frac{(2m+1)\pi}{\alpha}\right)
  \\
  &+
  4\pi
  \sum_{j}(-1)^{j}
  e^{N\left(\frac{(2j+1)^2\pi^2}{2\beta\xi}+(2j+1)\pi\sqrt{-1}\right)}
  \int_{C_{\varphi}}
  \frac{\sinh\left(\frac{z}{\alpha}+\frac{2(2j+1)\pi\sqrt{-1}}{\beta}\right)}
       {\cosh\left(\frac{z}{2}+\frac{(2j+1)\alpha\pi\sqrt{-1}}{\beta}\right)}
  e^{
    N
    \left(
      -\frac{\beta}{2\alpha(\beta-4\alpha)\xi}z^2
    \right)
  }
  \,dz
  \\
  &+
  16\pi^2
  \sum_{j,n}
  (-1)^{j+n}
  e^{N
    \left(
      \frac{\pi^2}{2\alpha(\beta-4\alpha)\xi}
      \bigl(
        (2j+1)^2\alpha+(2n+1)^2\beta-4(2j+1)(2n+1)\alpha
      \bigr)
      +(2j+1)\pi\sqrt{-1}
    \right)
  }
  \\
  &\quad\times
  \sin\left(\frac{(2n+1)\pi}{\alpha}\right),
\end{split}
\end{equation}
where
\begin{itemize}
\item
$j$ runs over integers such that $-2x+(2j+1)\pi\sqrt{-1}$ is between $C_{\varphi}$ and $C_{\varphi}+(\beta-4\alpha)\xi$ for $x\in C_{\varphi}$,
\item
$k$ runs over integers such that $\alpha\nmid(2k+1)$ and that $(2k+1)\pi\sqrt{-1}$ is between $C_{\varphi}$ and $C_{\varphi}+2\alpha\xi$,
\item
$l$ runs over integers such that $\frac{1}{2}\bigl((2l+1)\pi\sqrt{-1}-z-(\beta-4\alpha)\xi\bigr)$ is between $C_{\varphi}$ and $C_{\varphi}+2\alpha\xi$ for $z\in C_{\varphi}$,
\item
$m$ runs over integers such that $\alpha\nmid(2m+1)$ and that $(2l+1)\pi\sqrt{-1}-2(2m+1)\pi\sqrt{-1}-(\beta-4\alpha)\xi$ is between $C_{\varphi}$ and $C_{\varphi}+\frac{(\beta-4\alpha)\bigl((2l+1)\pi\sqrt{-1}-\beta\xi\bigr)}{\beta}$,
\item
$n$ runs over integers such that $\alpha\nmid(2n+1)$ and that $(2n+1)\pi\sqrt{-1}$ is between $C_{\varphi}$ and $C_{\varphi}+\frac{2(2j+1)\alpha\pi\sqrt{-1}}{\beta}$.
\end{itemize}
The ranges of indices can be simplified as follows:
\begin{itemize}
\item
$j$ runs over integers such that $(2j+1)\pi\sqrt{-1}$ is between $C_{\varphi}$ and $C_{\varphi}+(\beta-4\alpha)\xi$,
\item
$k$ runs over integers such that $\alpha\nmid(2k+1)$ and that $(2k+1)\pi\sqrt{-1}$ is between $C_{\varphi}$ and $C_{\varphi}+2\alpha\xi$,
\item
$l$ runs over integers such that $(2l+1)\pi\sqrt{-1}$ is between $C_{\varphi}+(\beta-4\alpha)\xi$ and $C_{\varphi}+\beta\xi$,
\item
$m$ runs over integers such that $\alpha\nmid(2m+1)$ and that $(2m+1)\pi\sqrt{-1}$ is between $C_{\varphi}-\frac{(\beta-4\alpha)\xi}{2}+\frac{(2l+1)\pi\sqrt{-1}}{2}$ and $C_{\varphi}+\frac{2(2l+1)\alpha\pi\sqrt{-1}}{\beta}$,
\item
$n$ runs over integers such that $\alpha\nmid(2n+1)$ and that $(2n+1)\pi\sqrt{-1}$ is between $C_{\varphi}$ and $C_{\varphi}+\frac{2(2j+1)\alpha\pi\sqrt{-1}}{\beta}$.
\end{itemize}
Replacing $w$ with $-2z$ in the third term of \eqref{eq:I_3_2}, we can combine the third and the fifth terms into one:
\begin{equation*}
  4\pi
  \sum_{j}(-1)^{j}
  e^{N\left(\frac{(2j+1)^2\pi^2}{2\beta\xi}+(2j+1)\pi\sqrt{-1}\right)}
  \int_{C_{\varphi}}
  \frac{\sinh\left(\frac{z}{\alpha}+\frac{2(2j+1)\pi\sqrt{-1}}{\beta}\right)}
       {\cosh\left(\frac{z}{2}+\frac{(2j+1)\alpha\pi\sqrt{-1}}{\beta}\right)}
  e^{
    N
    \left(
      -\frac{\beta}{2\alpha(\beta-4\alpha)\xi}z^2
    \right)
  }
  \,dz,
\end{equation*}
where $j$ runs over integers such that $(2j+1)\pi\sqrt{-1}$ is between $C_{\varphi}$ and $C_{\varphi}+\beta\xi$.
We can also combine the fourth and the sixth terms into one:
\begin{equation*}
\begin{split}
  &16\pi^2
  \sum_{l,m}
  (-1)^{l+m}
  e^{N
    \left(
      \frac{N\pi^2}{2\alpha(\beta-4\alpha)\xi}
      \bigl((2l+1)^2\alpha+(2m+1)^2\beta-4(2l+1)(2m+1)\alpha\bigr)
      +(2l+1)\pi\sqrt{-1}
    \right)
  }
  \\
  &\times
  \sin\left(\frac{(2m+1)\pi}{\alpha}\right),
\end{split}
\end{equation*}
where $(l,m)$ runs over pairs of integers such that
\begin{itemize}
\item
$\alpha\nmid(2m+1)$,
\item
$(2l+1)\pi\sqrt{-1}$ is between $C_{\varphi}$ and $C_{\varphi}+\beta\xi$,
\item
$(2m+1)\pi\sqrt{-1}$ is between $C_{\varphi}$ and $C_{\varphi}+\frac{2(2l+1)\alpha\pi\sqrt{-1}}{\beta}$ when $(2l+1)\pi\sqrt{-1}$ is between $C_{\varphi}$ and $C_{\varphi}+(\beta-4\alpha)\xi$,
\item
$(2m+1)\pi\sqrt{-1}$ is between $C_{\varphi}-\frac{(\beta-4\alpha)\xi}{2}+\frac{(2l+1)\pi\sqrt{-1}}{2}$ and $C_{\varphi}+\frac{2(2l+1)\alpha\pi\sqrt{-1}}{\beta}$ when $(2l+1)\pi\sqrt{-1}$ is between $C_{\varphi}+(\beta-4\alpha)\xi$ and $C_{\varphi}+\beta\xi$.
\end{itemize}
\par
Therefore we have
\begin{equation}\label{eq:I_3_3}
\begin{split}
  &I_3
  \\
  =&
  e^{\frac{N\xi}{2}\beta}
  \int_{C_{\varphi}\times C_{\varphi}}
  \frac{\sinh\left(\frac{w+2\alpha\xi}{\alpha}\right)
        e^{N\left(-\frac{w^2}{2\alpha\xi}-\frac{z^2}{2(\beta-4\alpha)\xi}\right)}}
       {\cosh\left(\frac{w+2\alpha\xi}{2}\right)
        \cosh\left(\frac{2w+z+\beta\xi}{2}\right)}
  \,dw\,dz
  \\
  &+
  4\pi
  \sum_{j}(-1)^{j}
  e^{N\left(\frac{(2j+1)^2\pi^2}{2\beta\xi}+(2j+1)\pi\sqrt{-1}\right)}
  \int_{C_{\varphi}}
  \frac{\sinh\left(\frac{z}{\alpha}+\frac{2(2j+1)\pi\sqrt{-1}}{\beta}\right)}
       {\cosh\left(\frac{z}{2}+\frac{(2j+1)\alpha\pi\sqrt{-1}}{\beta}\right)}
  e^{
    N
    \left(
      -\frac{\beta}{2\alpha(\beta-4\alpha)\xi}z^2
    \right)
  }
  \,dz
  \\
  &+
  4\pi\sqrt{-1}
  e^{\frac{N\xi}{2}(\beta-4\alpha)}
  \sum_{k}
  (-1)^{k+1}
  \sin\left(\frac{(2k+1)\pi}{\alpha}\right)
  e^{N\left(\frac{(2k+1)^2\pi^2}{2\alpha\xi}+2(2k+1)\pi\sqrt{-1}\right)}
  \\
  &\quad\times
  \int_{C_{\varphi}}
  \frac{e^{N\left(-\frac{z^2}{2(\beta-4\alpha)\xi}\right)}}{\cosh\left(\frac{z+(\beta-4\alpha)\xi}{2}\right)}
  \,dz
  \\
  &+
  16\pi^2
  \sum_{l,m}
  (-1)^{l+m}
  e^{N
    \left(
      \frac{\pi^2}{2\alpha(\beta-4\alpha)\xi}
      \bigl((2l+1)^2\alpha+(2m+1)^2\beta-4(2l+1)(2m+1)\alpha\bigr)
      +(2l+1)\pi\sqrt{-1}
    \right)
  }
  \\
  &\quad\times
  \sin\left(\frac{(2m+1)\pi}{\alpha}\right)
\end{split}
\end{equation}
where
\begin{itemize}
\item
$j$ runs over integers such that $(2j+1)\pi\sqrt{-1}$ is between $C_{\varphi}$ and $C_{\varphi}+\beta\xi$.
\item
$k$ runs over integers such that $\alpha\nmid(2k+1)$ and that $(2k+1)\pi\sqrt{-1}$ is between $C_{\varphi}$ and $C_{\varphi}+2\alpha\xi$,
\item
$(l,m)$ runs over pairs of integers such that
\begin{itemize}
\item
$\alpha\nmid(2m+1)$,
\item
$(2l+1)\pi\sqrt{-1}$ is between $C_{\varphi}$ and $C_{\varphi}+\beta\xi$,
\item
$(2m+1)\pi\sqrt{-1}$ is between $C_{\varphi}$ and $C_{\varphi}+\frac{2(2l+1)\alpha\pi\sqrt{-1}}{\beta}$ when $(2l+1)\pi\sqrt{-1}$ is between $C_{\varphi}$ and $C_{\varphi}+(\beta-4\alpha)\xi$,
\item
$(2m+1)\pi\sqrt{-1}$ is between $C_{\varphi}-\frac{(\beta-4\alpha)\xi}{2}+\frac{(2l+1)\pi\sqrt{-1}}{2}$ and $C_{\varphi}+\frac{2(2l+1)\alpha\pi\sqrt{-1}}{\beta}$ when $(2l+1)\pi\sqrt{-1}$ is between $C_{\varphi}+(\beta-4\alpha)\xi$ and $C_{\varphi}+\beta\xi$.
\end{itemize}
\end{itemize}
\par
By using the saddle point method (Theorem~\ref{thm:saddle}), the integrals in \eqref{eq:I_3_3} are approximated as follows:
\begin{gather*}
  \int_{C_{\varphi}\times C_{\varphi}}
  \frac{\sinh\left(\frac{w+2\alpha\xi}{\alpha}\right)
        e^{N\left(-\frac{w^2}{2\alpha\xi}-\frac{z^2}{2(\beta-4\alpha)\xi}\right)}}
       {\cosh\left(\frac{w+2\alpha\xi}{2}\right)
        \cosh\left(\frac{2w+z+\beta\xi}{2}\right)}
  \,dw\,dz
  \underset{N\to\infty}{\sim}
  \frac{2\pi\xi\sqrt{\alpha(\beta-4\alpha)}}{N}
  \frac{\sinh(2\xi)}{\cosh(\alpha\xi)\cosh\left(\frac{\beta\xi}{2}\right)},
  \\
  \int_{C_{\varphi}}
  \frac{\sinh\left(\frac{z}{\alpha}+\frac{2(2j+1)\pi\sqrt{-1}}{\beta}\right)}
       {\cosh\left(\frac{z}{2}+\frac{(2j+1)\alpha\pi\sqrt{-1}}{\beta}\right)}
  e^{N\left(-\frac{\beta}{2\alpha(\beta-4\alpha)\xi}z^2\right)}
  \,dz
  \underset{N\to\infty}{\sim}
  \sqrt{\frac{2\pi\alpha(\beta-4\alpha)\xi}{N\beta}}
  \frac{\sqrt{-1}\sin\left(\frac{2(2j+1)\pi}{\beta}\right)}
       {\cos\left(\frac{(2j+1)\alpha\pi}{\beta}\right)},
  \\
  \int_{C_{\varphi}}
  \frac{e^{N\left(-\frac{z^2}{2(\beta-4\alpha)\xi}\right)}}
       {\cosh\left(\frac{z+(\beta-4\alpha)\xi}{2}\right)}
  \,dz
  \underset{N\to\infty}{\sim}
  \sqrt{\frac{2\pi(\beta-4\alpha)\xi}{N}}
  \frac{1}{\cosh\left(\frac{(\beta-4\alpha)\xi}{2}\right)}.
\end{gather*}
So we have the following asymptotic equivalence.
\begin{equation*}
\begin{split}
  &I_3
  \\
  \underset{N\to\infty}{\sim}&
  e^{\frac{N\xi}{2}\beta}
  \frac{2\pi\xi\sqrt{\alpha(\beta-4\alpha)}}{N}
  \frac{\sinh(2\xi)}{\cosh(\alpha\xi)\cosh\left(\frac{\beta\xi}{2}\right)}
  \\
  &+
  4\pi\sqrt{-1}
  \sqrt{\frac{2\pi\alpha(\beta-4\alpha)\xi}{N\beta}}
  \sum_{j}(-1)^{j}
  e^{N\left(\frac{(2j+1)^2\pi^2}{2\beta\xi}+(2j+1)\pi\sqrt{-1}\right)}
  \frac{\sin\left(\frac{2(2j+1)\pi}{\beta}\right)}
       {\cos\left(\frac{(2j+1)\alpha\pi}{\beta}\right)}
  \\
  &+
  4\pi\sqrt{-1}
  e^{\frac{N\xi}{2}(\beta-4\alpha)}
  \sqrt{\frac{2\pi(\beta-4\alpha)\xi}{N}}
  \frac{1}{\cosh\left(\frac{(\beta-4\alpha)\xi}{2}\right)}
  \\
  &\quad\times
  \sum_{k}
  (-1)^{k+1}
  \sin\left(\frac{(2k+1)\pi}{\alpha}\right)
  e^{N\left(\frac{(2k+1)^2\pi^2}{2\alpha\xi}+2(2k+1)\pi\sqrt{-1}\right)}
  \\
  &+
  16\pi^2
  \sum_{l,m}
  (-1)^{l+m}
  e^{N
    \left(
      \frac{\pi^2}{2\alpha(\beta-4\alpha)\xi}
      \bigl((2l+1)^2\alpha+(2m+1)^2\beta-4(2l+1)(2m+1)\alpha\bigr)
      +(2l+1)\pi\sqrt{-1}
    \right)
  }
  \\
  &\quad\times
  \sin\left(\frac{(2m+1)\pi}{\alpha}\right).
\end{split}
\end{equation*}
Since it is easy to prove that $I_1\underset{N\to\infty}{\sim}0$ and $I_2\underset{N\to\infty}{\sim}0$, we finally have
\begin{equation*}
\begin{split}
  &J_{N}\left(T(2,\alpha)^{(2,\beta)};\exp(\xi/N)\right)
  \\
  \underset{N\to\infty}{\sim}&
  \frac{e^{\frac{-\beta(N^2-1)\xi}{2N}}}{2\sinh(\xi/2)}
  \frac{N}{4\xi\pi\sqrt{\alpha(\beta-4\alpha)}}
  e^{-\frac{\xi}{2N}\left(\frac{\beta-4\alpha}{4}+\frac{\alpha}{4}+\frac{1}{\alpha}\right)}
  I_3
  \\
  \underset{N\to\infty}{\sim}&
  \frac{\sinh(2\xi)}{4\sinh(\xi/2)\cosh(\alpha\xi)\cosh\left(\frac{\beta\xi}{2}\right)}
  \\
  &+
  \sqrt{-1}
  \frac{1}{2\sinh(\xi/2)}
  \sqrt{\frac{2\pi N}{\beta\xi}}
  \sum_{j}(-1)^{j}
  e^{N\left(\frac{-\beta\xi}{2}+\frac{(2j+1)^2\pi^2}{2\beta\xi}+(2j+1)\pi\sqrt{-1}\right)}
  \frac{\sin\left(\frac{2(2j+1)\pi}{\beta}\right)}
       {\cos\left(\frac{(2j+1)\alpha\pi}{\beta}\right)}
  \\
  &+
  \sqrt{-1}
  \frac{1}{2\sinh(\xi/2)}
  \sqrt{\frac{2\pi N}{\alpha\xi}}
  \sum_{k}
  (-1)^{k+1}
  e^{N\left(-2\alpha\xi+\frac{(2k+1)^2\pi^2}{2\alpha\xi}+2(2k+1)\pi\sqrt{-1}\right)}
  \frac{\sin\left(\frac{(2k+1)\pi}{\alpha}\right)}{\cosh\left(\frac{(\beta-4\alpha)\xi}{2}\right)}
  \\
  &+
  \frac{1}{2\sinh(\xi/2)}
  \frac{4\pi N}{\xi\sqrt{\alpha(\beta-4\alpha)}}
  \\
  &\quad\times
  \sum_{l,m}
  (-1)^{l+m}
  e^{N
    \left(
      \frac{-\beta\xi}{2}
      +
      \frac{\pi^2}{2\alpha(\beta-4\alpha)\xi}
      \bigl((2l+1)^2\alpha+(2m+1)^2\beta-4(2l+1)(2m+1)\alpha\bigr)
      +(2l+1)\pi\sqrt{-1}
    \right)}
  \sin\left(\frac{(2m+1)\pi}{\alpha}\right)
  \\
  =&
  \frac{1}{\Delta(T(2,2a+1)^{2b+1});\exp\xi}
  \\
  &+
  \frac{\sqrt{-\pi}}{2\sinh(\xi/2)}
  \sqrt{\frac{N}{\xi}}
  \sum_{j}
  \tau_1(\xi;j)
  \exp\left[\frac{N}{\xi}S_1(\xi;j)\right]
  \\
  &+
  \frac{\sqrt{-\pi}}{2\sinh(\xi/2)}
  \sqrt{\frac{N}{\xi}}
  \sum_{k}
  \tau_2(\xi;k)
  \exp\left[\frac{N}{\xi}S_2(\xi;k)\right]
  \\
  &+
  \frac{\pi}{2\sinh(\xi/2)}
  \frac{N}{\xi}
  \sum_{l,m}
  \tau_3(\xi;l,m)
  \exp\left[\frac{N}{\xi}S_3(\xi;l,m)\right],
\end{split}
\end{equation*}
where
\begin{align*}
  \tau_1(\xi;j)
  &:=
  (-1)^j
  \sqrt{\frac{2}{\beta}}
  \frac{\sin\left(\frac{2(2j+1)\pi}{\beta}\right)}
       {\cos\left(\frac{(2j+1)\alpha\pi}{\beta}\right)},
  \\
  S_1(\xi;j)
  &:=
  (-1)^{k+1}
  (2j+1)\xi\pi\sqrt{-1}
  -
  \frac{\beta\xi^2}{2}
  +
  \frac{(2j+1)^2\pi^2}{2\beta},
  \\
  \tau_2(\xi;k)
  &:=
  \sqrt{\frac{2}{\alpha}}
  \frac{\sin\left(\frac{(2k+1)\pi}{\alpha}\right)}{\cosh\left(\frac{(\beta-4\alpha)\xi}{2}\right)},
  \\
  S_2(\xi;k)
  &:=
  2(2k+1)\xi\pi\sqrt{-1}
  -
  2\alpha\xi^2
  +
  \frac{(2k+1)^2\pi^2}{2\alpha},
  \\
  \tau_3(\xi;l,m)
  &:=
  (-1)^{l+m}
  \frac{4}{\sqrt{\alpha(\beta-4\alpha)}}
  \sin\left(\frac{(2m+1)\pi}{\alpha}\right),
  \\
  S_3(\xi;l,m)
  &:=
  (2l+1)\xi\pi\sqrt{-1}
  -
  \frac{\beta\xi^2}{2}
  +
  \frac{\pi^2}{2\alpha(\beta-4\alpha)}\bigl((2l+1)^2\alpha+(2m+1)^2\beta-4(2l+1)(2m+1)\alpha\bigr).
\end{align*}
\section{Topological interpretation of the asymptotic expansion of the colored Jones polynomial}
\label{sec:topological}
In this section we study topological interpretation of the asymptotic expansion of the colored Jones polynomial.
Throughout this section we put $u:=\xi-2\pi\sqrt{-1}$.
\subsection{Figure-eight knot}
First of all we review the case of the figure-eight knot.
\par
The original volume conjecture (Conjecture~\ref{conj:VC}) for the figure-eight knot was prove by T.~Ekholm (see for example \cite{Murakami:Novosibirsk}).
Conjecture~\ref{conj:Gukov/Murakami} was proved by Yokota and the author \cite{Murakami/Yokota:JREIA2007} in the case of the figure-eight knot.
The following theorem appears in \cite{Murakami:JTOP2013}, proving Conjecture~\ref{conj:Dimofte/Gukov}.
\begin{thm}[{\cite[Theorem~1.4]{Murakami:JTOP2013}}]
Let $\FigEight$ be the figure-eight knot.
If $u$ is real and sufficiently small, then we have
\begin{multline*}
  J_N\Bigl(\FigEight;\exp\bigl((2\pi\sqrt{-1}+u)/N\bigr)\Bigr)
  \\
  =
  \frac{\sqrt{-\pi}}{2\sinh(u/2)}
  \mathbb{T}^{\FigEight}_{\mu}(u)^{-1/2}
  \left(\frac{N}{2\pi\sqrt{-1}+u}\right)^{1/2}
  \exp\left[\frac{S(u)N}{2\pi\sqrt{-1}+u}\right],
\end{multline*}
where
\begin{equation*}
  \mathbb{T}^{\FigEight}_\mu(u)
  =
  \frac{\sqrt{(2\cosh{u}+1)(2\cosh{u}-3)}}{2}
\end{equation*}
is the twisted Reidemeister torsion of the representation $\rho_{u,+}$ associated with the meridian $\mu$, and
\begin{equation*}
  S(u)
  :=
  \Li_2\left(e^{u-\varphi(u)}\right)-\Li_2\left(e^{u+\varphi(u)}\right)-u\varphi(u).
\end{equation*}
Here $\varphi(u):=\arccosh\bigl(\cosh{u}-1/2\bigr)$ and $\displaystyle\Li_2(z):=-\int_{0}^{z}\frac{\log(1-x)}{x}\,dx$ is the dilogarithm function.
Moreover if we define $v(u):=2\dfrac{d\,S(u)}{d\,u}-2\pi\sqrt{-1}$, then we have
\begin{equation*}
  \CS_{u,v(u)}(\rho_{u,+})
  =
  S(u)-u\pi\sqrt{-1}-\frac{uv(u)}{4},
\end{equation*}
where $\CS_{u,v(u)}$ is the Chern--Simons invariant of $\rho_{u,+}$ associated with $(u,v(u))$.
\end{thm}
Note that by taking the derivative, one can confirm \eqref{eq:figure8_CS} since $S(u)$ coincides with $\sqrt{-1}\Vol\left(S^3\setminus\FigEight\right)$.
Note also that \eqref{eq:fig8_Reidemeister} identifies $\mathbb{T}^{\FigEight}_{\mu}(u)$.
\subsection{Torus knot}
Let $S_k(\xi)$ and $\tau_k$ be as defined in \eqref{eq:asymptotic_torus_CS} and \eqref{eq:asymptotic_torus_Reidemeister}, respectively.
We put
\begin{equation*}
\begin{split}
  v_k(u)
  &:=
  2\frac{d\,S_k(u+2\pi\sqrt{-1})}{d\,u}-2\pi\sqrt{-1}
  \\
  &=
  -2(2a+1)u-4(2a+1-k)\pi\sqrt{-1}.
\end{split}
\end{equation*}
Then we have
\begin{equation*}
  S_k(u+2\pi\sqrt{-1})-\pi\sqrt{-1}u-\frac{1}{4}uv_{k}(u)
  =
  \frac{(2k+1)^2\pi^2}{2(2a+1)}
  +4(a-k)\pi^2
  -(2a+1-k)u\pi\sqrt{-1}.
\end{equation*}
So if we put $l:=-2(2a+1-k)$ and $\omega_1:=\exp\left(\frac{(2k+1)\pi\sqrt{-1}}{2a+1}\right)$, we have
\begin{equation*}
  \CS_{u,v}(\rho_{u,\omega_1})
  =
  S_k(u+2\pi\sqrt{-1})-\pi\sqrt{-1}u-\frac{1}{4}u v_{k}(u)
\end{equation*}
modulo $\pi^2\Z$ with $v:=-2(2a+1)u+2l\pi\sqrt{-1}=-2(2a+1)u-2(2a+1-k)\pi\sqrt{-1}$ from \eqref{eq:torus_knot_CS}.
\begin{rem}\label{rem:torus_knot_topological}
I have chosen $l$ to be an odd integer, but here I say $l$ is even.
Recall that $l$ should have been odd from \eqref{eq:rep_longitude_torus_knot}.
So if we regard $\rho_{u,\omega_1}$ as a representation to $\operatorname{PSL}(2;\C)$, not $\SL(2;\C)$, we can avoid this trouble.
\end{rem}
\par
We also have
\begin{equation*}
  \tau_{k}^{-2}
  =
  \mathbb{T}^{T(2,2a+1)}_{\mu}\left(\rho_{u,\omega_1}\right)
\end{equation*}
up to a sign from \eqref{eq:torus_knot_Reidemeister}.
\par
Therefore we have the following theorem, which is a special case of Theorem~\ref{thm:Hikami/Murakami_intro}.
\begin{thm}[\cite{Hikami/Murakami:Bonn}]
\label{thm:Hikami/Murakami}
For a complex number $\xi$ that is not purely imaginary and $\Im\xi\ge0$, we have the following asymptotic equivalence:
\begin{multline*}
  J_{N}\bigl(T(2,2a+1);\exp(\xi/N)\bigr)
  \\
  \underset{N\to\infty}{\sim}
  \frac{1}{\Delta\bigl(T(2,2a+1);\exp(\xi)\bigr)}
  +
  \frac{\sqrt{-\pi}}{2\sinh(\xi/2)}\sqrt{\frac{N}{\xi}}
  \sum_{k}\tau_k\exp\left[S_k(\xi)\frac{N}{\xi}\right],
\end{multline*}
where 
\begin{align*}
  S_k(\xi)
  &:=
  \frac{-\bigl((2k+1)\pi\sqrt{-1}-(2a+1)\xi\bigr)^2}{2(2a+1)},
  \\
  \intertext{and}
  \tau_k
  &:=
  (-1)^k\frac{4\sin\left(\frac{(2k+1)\pi}{2a+1}\right)}{\sqrt{2(2a+1)}}.
\end{align*}
Moreover $\tau_k^{-2}$ is the homological twisted Reidemeister torsion $\mathbb{T}^{T(2,2a+1)}_k$ of the irreducible representation $\rho_{u,\omega_1}\colon\pi_1(S^3\setminus{T(2,2a+1)})\to\SL(2;\C)$ \rm{(}$\omega_1:=\exp\left(\frac{(2k+1)\pi\sqrt{-1}}{2a+1}\right)$\rm{)} associated with the meridian, and $S_k(\xi)-\pi\sqrt{-1}u-\frac{1}{4}uv_k(u)$ is the $\SL(2;\C)$ Chern--Simons invariant of $\rho_{u,\omega_1}$ with respect to the pair $(u,v_k(u))$ with $u:=\xi-2\pi\sqrt{-1}$ and $v_k(u):=2\dfrac{d\,S_k(\xi)}{d\,\xi}\Big|_{\xi:=2\pi\sqrt{-1}+u}-2\pi\sqrt{-1}$.
\end{thm}
\subsection{Twice-iterated torus knot}
\label{subsec:topological_iterated_torus_knot}
In this subsection we will prove Theorem~\ref{thm:main}.
\par
In the following sub-subsections (Sub-subsection~\ref{subsubsec:AN_topological}--\ref{subsubsec:NN_topological}) I will relate $S_1(\xi;j)$, $S_2(\xi;k)$, $S_3(\xi;l,m)$, $\tau_1(\xi;j)$, $\tau_2(\xi;k)$, and $\tau_3(\xi;l,m)$ to the Chern--Simons invariants and the Reidemeister torsions.
\par
Recall that we put $\alpha=2a+1$ and $\beta=2b+1$.
\subsubsection{$\Im\rho_{C}$ is Abelian and $\Im\rho_{P}$ is non-Abelian}
\label{subsubsec:AN_topological}
We put
\begin{equation*}
\begin{split}
  v_{1,\omega_2}(u)
  &:=
  2\frac{d\,S_1(u+2\pi\sqrt{-1};j)}{d\,u}-2\pi\sqrt{-1}
  \\
  &=
  -2(2b+1)u-4(2b+1-j)\pi\sqrt{-1}.
\end{split}
\end{equation*}
Then we have
\begin{equation*}
  S_1(u+2\pi\sqrt{-1};j)-\pi\sqrt{-1}u-\frac{1}{4}u v_{1,\omega_2}(u)
  =
  \frac{(2j+1)^2\pi^2}{2(2b+1)}
  +4(b-j)\pi^2
  -(2b+1-j)u\pi\sqrt{-1}.
\end{equation*}
So if we put $m:=-2(2b+1-j)$ and $\omega_2:=\exp\left(\frac{(2j+1)\pi\sqrt{-1}}{2b+1}\right)$, we have
\begin{equation*}
  \CS_{u,v}\left(\rho^{\rm{AN}}_{u,\omega_2}\right)
  =
  S_1(u+2\pi\sqrt{-1};j)-\pi\sqrt{-1}u-\frac{1}{4}u v_{1,\omega_2}(u)
\end{equation*}
modulo $\pi^2\Z$ with $v:=-2(2b+1)u+2m\pi\sqrt{-1}$ from \eqref{eq:AN_CS}.
For the parity of $m$, see Remark~\ref{rem:torus_knot_topological}.
\par
We also have
\begin{equation*}
  \tau_{1}(\xi;j)^{-2}
  =
  \mathbb{T}^{T(2,2a+1)^{(2,2b+1)}}_{\mu}\left(\rho^{\rm{AN}}_{u,\omega_2}\right)
\end{equation*}
up to a sign from \eqref{eq:AN_Reidemeister}.
\subsubsection{$\Im\rho_{C}$ is non-Abelian and $\Im\rho_{P}$ is Abelian}
\label{subsubsec:NA_topological}
We put
\begin{equation*}
\begin{split}
  v_{2,\omega_1}(u)
  &:=
  2\frac{d\,S_2(u+2\pi\sqrt{-1};k)}{d\,u}-2\pi\sqrt{-1}
  \\
  &=
  -8(2a+1)u-2\bigl(8(2a+1)-4k-1\bigr)\pi\sqrt{-1}.
\end{split}
\end{equation*}
Then we have
\begin{multline*}
  S_2(u+2\pi\sqrt{-1};k)-\pi\sqrt{-1}u-\frac{1}{4}uv_{2,\omega_1}(u)
  \\
  =
  \frac{(2k+1)^2\pi^2}{2(2a+1)}
  +4\bigl(2(2a+1)-2k-1\bigr)\pi^2
  -\frac{1}{2}\bigl(8(2a+1)-4k-1\bigr)u\pi\sqrt{-1}.
\end{multline*}
So if we put $l:=-\dfrac{1}{2}\bigl(8(2a+1)-4k-1\bigr)$ and $\omega_1:=\exp\left(\frac{(2k+1)\pi\sqrt{-1}}{2a+1}\right)$, we have
\begin{equation*}
  \CS_{u,v}\left(\rho^{\rm{NA}}_{u,\omega_1}\right)
  =
  S_2(u+2\pi\sqrt{-1};k)-\pi\sqrt{-1}u-\frac{1}{4}u v_{2,\omega_1}(u)
\end{equation*}
modulo $\pi^2\Z$ with $v:=-8(2a+1)u+4l\pi\sqrt{-1}$ from \eqref{eq:NA_CS}.
For the integrality of $l$, see Remark~\ref{rem:torus_knot_topological}.
\par
We also have
\begin{equation*}
  \tau_{2}(\xi;k)^{-2}
  =
  \mathbb{T}^{T(2,2a+1)^{(2,2b+1)}}_{\mu}\left(\rho^{\rm{NA}}_{u,\omega_1}\right)
\end{equation*}
up to a sign from \eqref{eq:NA_Reidemeister}.
\subsubsection{Both $\Im\rho_{C}$ and $\Im\rho_{P}$ are non-Abelian}
\label{subsubsec:NN_topological}
We put
\begin{equation*}
\begin{split}
  v_{3,\omega_1,\omega_3}(u)
  &:=
  2\frac{d\,S_3(u+2\pi\sqrt{-1};l,m)}{d\,u}-2\pi\sqrt{-1}
  \\
  &=
  -2(2b+1)u-4(2b+1-l)\pi\sqrt{-1}.
\end{split}
\end{equation*}
Then we have
\begin{equation*}
\begin{split}
  &S_3(u+2\pi\sqrt{-1};k)-\pi\sqrt{-1}u-\frac{1}{4}uv_{3,\omega_1,\omega_3}(u)
  \\
  =&
  \frac{\pi^2}{2\alpha(\beta-4\alpha)}
  \left((2l+1)^2\alpha+(2m+1)^2\beta-4(2l+1)(2m+1)\alpha\right)
  \\
  &
  +4(b-l)\pi^2
  -(2b+1-l)u\pi\sqrt{-1}
  \\
  =&
  \frac{\pi^2}{2(\beta-4\alpha)}\bigl(2l+1-2(2m+1)\bigr)^2
  +
  \frac{\pi^2}{2\alpha}(2m+1)^2
  \\
  &
  +4(b-l)\pi^2
  -(2b+1-l)u\pi\sqrt{-1}.
\end{split}
\end{equation*}
So if we put $n:=-2(2b+1-l)$, $k:=m$, $h:=l-2m-1$, $\omega_1:=\exp\left(\frac{(2k+1)\pi\sqrt{-1}}{2a+1}\right)$ and $\omega_3:=\exp\left(\frac{(2h+1)\pi\sqrt{-1}}{2b+1-4(2a+1)}\right)$, we have
\begin{equation*}
  \CS_{u,v}\left(\rho^{\rm{NN}}_{u,\omega_1,\omega_3}\right)
  =
  S_3(u+2\pi\sqrt{-1};l,m)-\pi\sqrt{-1}u-\frac{1}{4}u v_{3,\omega_1,\omega_3}(u)
\end{equation*}
modulo $\pi^2\Z$ with $v:=-2(2b+1)u+2n\pi\sqrt{-1}$ from \eqref{eq:NN_CS}.
For the parity of $n$, see Remark~\ref{rem:torus_knot_topological}.
\par
Unfortunately, $\tau_3(\xi;l,m)^{-2}=\frac{(2a+1)\bigl(2b+1-4(2a+1)\bigr)}{16\sin^2\left(\frac{(2m+1)\pi}{2a+1}\right)}$ does not coincide with the right hand side of \eqref{eq:NN_Reidemeister}.
\bibliography{mrabbrev,hitoshi}
\bibliographystyle{amsplain}
\end{document}